\documentclass[11pt]{amsart}
\usepackage[top=25mm, bottom=25mm, left=25mm, right=25mm]{geometry}

\newcommand{\tikzannotsize}{\scriptsize}
\newcommand{\tikzlabelsize}{\small}
\newcommand{\tablesize}{\small}
\newcommand{\figlabelsize}{\normalsize}
\newcommand{\figtitlesize}{\large}

\usepackage{amsmath}
\usepackage{amssymb}
\usepackage{bm}

\usepackage{adjustbox}
\usepackage{graphicx}
\usepackage{amssymb}
\usepackage{bm}

\usepackage{caption}
\usepackage{subcaption}
\usepackage[normalem]{ulem}
\usepackage{setspace,booktabs,lipsum,changepage,makecell}

\usepackage{pgfplots}
\pgfplotsset{compat=1.16}          \usepgfplotslibrary{groupplots}

\usepackage{tikz,tikz-3dplot}
\usetikzlibrary{decorations.pathreplacing,calligraphy}
\usetikzlibrary{positioning}
\usetikzlibrary{calc,patterns}
\usetikzlibrary{fit,fpu}
\usepackage{ifthen}
\usepackage{algorithm}
\usepackage[noend]{algpseudocode}

\usepackage{multirow}

\newcommand{\msf}[1]{\mathsf{#1}}

\newcommand{\pvct}[1]{\bm{#1}}

\newcommand{\mtx}[1]{\bm{\mathsf{#1}}}

\newcommand{\pxx}{\pvct{x}}
\newcommand{\pyy}{\pvct{y}}

\newcommand{\tree}{\mathcal{T}}

\pgfplotsset{
  every axis/.append style={
    width=\linewidth,          height=\linewidth,     scale only axis,           font=\fontsize{11}{13}\selectfont,       tick label style={font=\fontsize{10}{12}\selectfont},
    label style={font=\fontsize{12}{14}\selectfont},
    title style={font=\fontsize{12}{14}\selectfont},
    legend style={font=\fontsize{9}{11}\selectfont},
  },
  every axis plot/.append style={line width=1pt, mark size=2.2pt} }

\definecolor{myblue}{HTML}{1F77B4}
\definecolor{myred}{HTML}{D62728}
\definecolor{mygreen}{HTML}{2CA02C}

\pgfplotsset{
  cycle list={
    {color=myblue, mark=*},
    {color=myred, mark=square*},
    {color=mygreen, mark=triangle*}
  }
}

\newcommand{\addONguideRel}[3]{\pgfmathsetmacro{\xstart}{#1}\pgfmathsetmacro{\ystart}{#2}\pgfmathsetmacro{\xend}{#1 + #3}\pgfmathsetmacro{\yend}{#2 + #3}\draw[dotted, very thick, gray]
    (rel axis cs:\xstart,\ystart) -- (rel axis cs:\xend,\yend);
\node[anchor=north, font=\figlabelsize,yshift=-1.0em]
    at (rel axis cs:\xend,\yend) {$\mathcal{O}(N)$};
}
 
\newtheorem{remark}{Remark}

\title[Randomized Strong Recursive Skeletonization]{Randomized Strong Recursive Skeletonization:\\ Compression and LU Factorization of Hierarchical Matrices using Matrix--Vector Products}

\author[Yesypenko]{Anna Yesypenko\textsuperscript{*}}
\thanks{\textsuperscript{*}Department of Mathematics, The Ohio State University,
Columbus, OH, USA. \texttt{yesypenko.1@osu.edu}}

\author[Martinsson]{Per-Gunnar Martinsson\textsuperscript{$\dagger$}}
\thanks{\textsuperscript{$\dagger$}Department of Mathematics, The University of Texas at Austin,
Austin, TX, USA. \texttt{pgm@oden.utexas.edu}}

\begin{document}

\begin{abstract}
The hierarchical matrix framework partitions matrices into subblocks that are either small or of low numerical rank, enabling linear storage complexity and efficient matrix--vector multiplication. This work focuses on the $\mathcal{H}^2$-matrix format constructed under the strong admissibility condition, which has two key properties: (1) a compressed representation that approximates far-field interactions with low-rank blocks while near-field interactions are stored densely, and (2) a nested basis structure that reuses basis matrices across hierarchy levels. Although these matrices support fast Cholesky and LU factorizations, implementing them---especially for 3D PDE discretizations---remains challenging due to the nested recursions and recompressions involved.

This paper introduces an algorithm that simultaneously compresses and factorizes a general $\mathcal{H}^{2}$ matrix, using only the action of the matrix and its adjoint on vectors. The number of required matrix--vector products is independent of the matrix size, and depends only on the problem geometry and a rank parameter. 
The resulting LU factorization is invertible and can serve as an approximate direct solver, with accuracy influenced by the spectral properties of the matrix.

To achieve competitive sample complexity, the method employs dense Gaussian test matrices without explicitly encoding structured sparsity. Samples are drawn only once at the start of the algorithm; as the factorization proceeds, structure is dynamically introduced into the test matrices through efficient linear algebraic operations.  Numerical experiments demonstrate robustness to indefiniteness and ill-conditioning, as well as the efficiency of the method for integral and differential equations in 2D and 3D. \end{abstract}

\maketitle

\section{Introduction} \label{sec:intro}

Dense matrices arising in mathematical physics often exhibit 
internal structures that enable the efficient solution of problems 
involving elliptic partial differential equations (PDEs). 
In many cases, a dense matrix can be partitioned into subblocks,
each of which is either small or admits a compressed representation. 
Early methods such as the Fast 
Multipole Method (FMM) \cite{rokhlin1987} leveraged these properties 
to accelerate matrix--vector operations for a specific class of 
matrices. Later, the $\mathcal{H}^2$-matrix methodology 
\cite{2008_bebendorf_book,2003_borm_introduction_H_matrix,hackbusch2015hierarchical}
reinterpreted 
the FMM in algebraic terms by representing off-diagonal subblocks as approximately low rank,
thereby extending fast algorithmic techniques to a broader class of matrices.
Suppose we have a dense matrix $\mtx A \in \mathbb{R}^{N \times N}$ that is compressible as an $\mathcal{H}^2$ matrix and that a fast matrix--vector product algorithm is available. In this case, the linear system 
\begin{equation*} \mtx A\,\mtx u = \mtx f 
\end{equation*} can be solved using iterative methods. However, these methods may exhibit unsatisfactory convergence or may be impractical when solving for multiple right-hand sides.

The $\mathcal{H}^2$-matrix methodology not only enables fast matrix--vector products but also supports operations such as matrix--matrix multiplication and invertible factorization in linear complexity. However, implementing these algorithms, especially for matrices from three-dimensional PDE discretizations, remains challenging. In this work, we focus on direct inversion techniques for $\mathcal{H}^2$ matrices under the strong admissibility condition.

Two key features characterize this format. First, the $\mathcal{H}^2$ notation refers to a dual-hierarchy structure that achieves linear complexity through recursive partitioning and the use of nested bases. The representation captures interactions across multiple scales via a tree-based geometric partitioning, with low-rank bases from finer levels successively incorporated into coarser levels. Second, under the strong admissibility condition, interactions between a box and sufficiently distant boxes are approximated as numerically low rank, whereas interactions with neighboring boxes are stored densely. This approach, inspired by the FMM \cite{1997_greengard_shortcourse,rokhlin1987}, is optimal for compressing matrices arising from surfaces or volumes in 3D. The storage cost scales linearly with the matrix size $N$, with a constant that depends on the geometry and the compression rank.

Despite the advantages this matrix format offers, 
the practical use of $\mathcal{H}^2$ inversion algorithms under the strong admissibility condition has been limited by two core challenges: 
\begin{enumerate} 
\item \textbf{Compression as $\mathcal{H}^2$}. While efficient techniques exist for obtaining an $\mathcal{H}^2$ representation of a matrix in certain settings where matrix entries are explicitly available---such as adaptive cross approximation \cite{bebendorf2000approximation,bebendorf2003adaptive} and proxy point techniques \cite{xing2020interpolative,ye2020analytical}---the case where $\mtx{A}$ is only available through its action on vectors remains largely open. Randomized techniques have been proposed to address this problem \cite{2011_lin_lu_ying}, but the prefactors involved tend to be very large.
\item \textbf{Invertible factorization of an $\mathcal{H}^2$ matrix}. 
The problem of inverting an $\mathcal{H}^{2}$ matrix remains highly challenging, with existing methods relying either on repeated recompressions in recursive structures \cite{2008_bebendorf_book,2010_borm_book},
or on highly storage-intensive data structures \cite{2017_ho_ying_strong_RS,sushnikova2023fmm}. There do exist much faster algorithms for specialized subclasses of matrices \cite{2012_greengard_ho_recursive_skeletonization,2005_martinsson_fastdirect,2010_xia}, but these tend to not be suitable for problems arising from general geometries in three dimensions.
\end{enumerate}

This manuscript introduces Randomized Strong Recursive Skeletonization (\texttt{RSRS}), 
an algorithm that \textit{simultaneously} compresses and inverts $\mathcal{H}^2$ matrices under the strong admissibility condition.
The algorithm produces the ``Strong Recursive Skeletonization'' factorization introduced in \cite{2017_ho_ying_strong_RS} and further improved in \cite{sushnikova2023fmm} (cf.~also \cite{2014_darve_IFMM}).
\texttt{RSRS} broadens the applicability of $\mathcal{H}^2$ inversion to a wide range of dense matrices for which  fast matrix--vector products with the matrix and its adjoint are available (in contrast \cite{2017_ho_ying_strong_RS,sushnikova2023fmm} which assume that matrix entries of $\mtx{A}$ are readily available).

Specifically, assume that $\mtx{A} \in \mathbb{R}^{N \times N}$ is an $\mathcal{H}^2$ matrix equipped with geometric information for its rows and columns,
and that fast routines exist to apply both $\mtx{A}$ and its adjoint $\mtx{A}^*$ to arbitrary vectors. 
The \texttt{RSRS} algorithm requires two tall, thin matrices $\mtx{\Omega}$ and $\mtx{\Psi}$ with entries drawn from the standard normal distribution, which define the random sketches:
\begin{equation*}
\underset{N \times s}{\mtx{Y}}\ =\ \underset{N \times N}{\mtx{A}}\  \underset{N \times s}{\mtx{\Omega}}
\qquad\mbox{and}\qquad
\underset{N \times s}{\mtx{Z}}\ =\ \underset{N \times N}{\mtx{A}^*} \underset{N \times s}{\mtx{\Psi}},
\end{equation*}
where $\mtx \Omega, \mtx \Psi \sim \mathcal N(\mtx 0,\mtx I)$ are standard Gaussian matrices.
The method then reconstructs an approximate invertible factorization of $\mtx{A}$ by post-processing
$$\left\{\underset{N\times s}{\mtx Y},\ \underset{N\times s}{\mtx Z},\ \underset{N\times s}{\mtx \Omega},\ \underset{N\times s}{\mtx \Psi}\right\}$$
without directly accessing individual entries of $\mtx{A}$. 
The number of samples $s$ needed depends linearly on the maximal rank of the off-diagonal blocks, but is independent of the matrix size $N$.

\texttt{RSRS} is immediately applicable in a range of important environments. As a solver for boundary integral equations, \texttt{RSRS} can be used to compute an approximate inverse of any integral operator $\mtx{A}$ for which a fast matrix–vector multiplication algorithm, such as the Fast Fourier Transform (FFT) \cite{cooley1965algorithm} or the FMM~\cite{rokhlin1987,rokhlin1997}, is available. As a solver for PDE discretizations, \texttt{RSRS} can substantially accelerate and simplify the treatment of dense blocks that arise in the course of sparse LU factorization. The rank structure in these dense matrices can be exploited to achieve competitive complexity and high practical performance in sparse direct solvers \cite{2017_amestoy_BLR_sparse_direct,2009_martinsson_FEM,2019_martinsson_fast_direct_solvers,2009_xia_superfast,yesypenko2024slablu}. In uncertainty quantification, \texttt{RSRS} may be useful for factorizing the Hessian or dense matrices associated with the Jacobian that arise in PDE-constrained optimization problems \cite{alger2019data,ambartsumyan2020hierarchical,isaac2015scalable}.

\subsection{Key Insights and Contributions}

The task of recovering hierarchical matrices from matrix--vector products is an active research area, with wider implications to the recovery of structured representations of PDE models from data \cite{boulle2022data,boulle2023mathematical}. 
A key challenge in black-box randomized algorithms is designing test vectors that
allow for approximate matrix recovery with competitive sample complexity---that is, using as few matrix--vector products as possible.
Some previous works \cite{2011_lin_lu_ying,2016_martinsson_hudson2,2011_martinsson_randomhudson} rely on accessing matrix entries or on designing test matrices with carefully placed zeros to (1) sketch low-rank subblocks and (2) extract small subblocks. However, for LU factorization of $\mathcal{H}^2$ matrices, such techniques can increase sample complexity and complicate the practical implementation. \texttt{RSRS} avoids these issues by employing generic dense Gaussian test matrices and leveraging specialized randomized sketching techniques introduced in \cite{2022_levitt_dissertation,levitt2024linear}. Any required structure is imposed via linear transformations applied to the test matrices, rather than by explicitly encoding sparsity patterns or fixed zeros.

The algorithm \texttt{RSRS} simultaneously compresses and factorizes the matrix, representing a significant departure from traditional methods, which treat these stages separately. 
Previous work on black-box randomized algorithms for rank-structured matrices \cite{levitt2024linear,2011_lin_lu_ying,2016_martinsson_hudson2,2011_martinsson_randomhudson,2013_xia_randomized} 
has primarily focused on weak admissibility, due to its simplicity and the availability of exact inversion algorithms \cite{2022_martinsson_gpu_hodlr,2012_martinsson_FDS_survey,2010_xia}. 
In contrast, strong admissibility poses additional challenges, such as the need for repeated recompression of off-diagonal blocks during inversion, which has historically hindered its application to large 3D problems \cite{2014_darve_IFMM,2008_bebendorf_book,2010_borm_book}. Compared to previous algorithms for $\mathcal{H}^2$ inversion, \texttt{RSRS} introduces key innovations by combining the compression and inversion steps, leading to improved computational efficiency and reduced storage requirements, as the updated off-diagonal blocks are maintained in a compressed form, rather than stored explicitly.

A core idea motivating our approach is how randomized sketching behaves under multiplicative transformations.
Suppose we are given a randomized sample pair $(\mtx{Y}, \mtx{\Omega})$ such that $\mtx{Y} = \mtx{A} \mtx{\Omega}$.
Now consider an updated matrix $\hat{\mtx{A}}$ obtained by applying invertible left and right transformations:
\[
\hat{\mtx{A}} = \mtx{L} \mtx{A} \mtx{U}.
\]
We can then define a new sample pair $(\hat{\mtx{Y}}, \hat{\mtx{\Omega}})$ for $\hat{\mtx{A}}$ as
\[
\hat{\mtx{Y}} = \mtx{L} \mtx{Y}, \qquad \hat{\mtx{\Omega}} = \mtx{U}^{-1} \mtx{\Omega}.
\]
That $(\hat{\mtx{Y}}, \hat{\mtx{\Omega}})$ is indeed a sample pair for $\hat{\mtx{A}}$ 
follows directly from the calculation
\[
\hat{\mtx{Y}} = \mtx{L} \mtx{A} \mtx{\Omega} = \mtx{L} \mtx{A} \mtx{U} \mtx{U}^{-1} \mtx{\Omega} = \hat{\mtx{A}} \hat{\mtx{\Omega}}.
\]
Observe that if $\mtx{L}$ and $\mtx{U}$ are well-conditioned and independent of the test matrix $\mtx{\Omega}$, then the new sample pair $(\hat{\mtx{Y}}, \hat{\mtx{\Omega}})$ provides a high-quality sketch of $\hat{\mtx{A}}$. 
In the method described, the matrices $\mtx{L}$ and $\mtx{U}$ are not entirely independent of $\mtx{\Omega}$, but extensive numerical experiments indicate that this dependence is sufficiently weak that sketches remain accurate throughout the compression and factorization process.
This principle underlies the \texttt{RSRS} algorithm: as we recursively apply structured factorizations to compress and invert $\mtx{A}$, we update sketch matrices in tandem using transformations derived from the factorization itself, without ever needing direct access to matrix entries.

\subsection{Outline of the Paper}

The manuscript is organized as follows: Section \ref{sec:linalg_preliminaries} reviews the mathematical preliminaries of randomized linear algebra and sparse LU factorization. Section \ref{sec:skel_algos} describes LU factorizations of $\mathcal{H}^2$ matrices under the strong admissibility condition, focusing on the case where matrix entries are directly available. Section \ref{sec:randomized_compression_and_lu} presents the computation of an invertible factorization when the matrix is accessible only through its action on vectors. Two key tools are introduced for efficiently computing low-rank factorizations of admissible subblocks and extracting dense subblocks using matrix--vector products. Finally, Section \ref{sec:numerical_results} presents numerical results for various discretizations and complex geometries in two and three dimensions, demonstrating the performance of \texttt{RSRS} in terms of speed, sample requirements, and accuracy. 
\section{Preliminaries}
\label{sec:linalg_preliminaries}

We briefly summarize the notations used throughout the paper.
Let $\msf{I},\msf{J}$ denote ordered index sets.
The notation $\mtx A_{\msf I \msf J}$ denotes the subblock of matrix $\mtx A$ corresponding to the set of row indices $\msf I$ and the column indices $\msf J$.
The Euclidean norm of a vector $\mtx{x}$ is $\|\mtx{x}\|$ and for a given matrix $\mtx{A}$, the induced operator norm (spectral) norm is written as $\|\mtx{A}\|$.

We also introduce some shorthands for common linear algebraic constructions.
A matrix $\mtx Q \in \mathbb{R}^{m \times n}$ is to have orthonormal columns if $\mtx Q^* \mtx Q = \mtx I$.
For $m\geq n$, the operation
\begin{equation}
\texttt{orth}(\mtx A) := \mtx Q \in \mathbb{R}^{m \times n}
\label{eq:orth}
\end{equation}
applied to a matrix $\mtx A \in \mathbb{R}^{m \times n}$ returns an orthonormal basis for the column space of $\mtx A$ such that $\mtx Q \mtx Q^* \mtx A = \mtx A$. The notation
\begin{equation}
\texttt{orth}_k(\mtx A) := \mtx Q \in \mathbb{R}^{m \times k}
\label{eq:orthk}
\end{equation}
gives an orthonormal basis for the space spanned by the dominant $k$ left singular vectors of $\mtx A$. 

For $m < n$, the operation
\begin{equation}
\texttt{null}(\mtx A) := \mtx Q \in \mathbb{R}^{n \times (n-m)}
\label{eq:null}
\end{equation}
gives an orthonormal basis for the nullspace of $\mtx A$.
The orthogonal bases in (\ref{eq:orth})–(\ref{eq:null}) are not unique—they are defined up to unitary transformations—but this non-uniqueness has no effect on the proposed methods.
The dagger notation denotes the Moore–Penrose pseudoinverse, which is uniquely defined for any matrix; for example, if $\mtx A \in \mathbb{R}^{n \times m}$ with $m < n$, then $\mtx A\mtx A^{\dagger} = \mtx I_m$.

\subsection{Randomized Low Rank Approximation}
\label{sec:randomized_lr}
Suppose we would like to compute a rank-$k$ approximation to the matrix $\mtx A \in \mathbb{R}^{m \times n}$, i.e.\ to find matrix $\mtx Q \in \mathbb{R}^{m \times k}$ with orthonormal columns 
and some matrix $\mtx B \in \mathbb{R}^{k \times n}$ such that
\begin{equation*}
\left \| \underset{m \times n}{\mtx A\vphantom{Q}}  - \underset{m \times k}{\mtx Q}\ \underset{k \times n}{\vphantom{\mtx Q}\mtx B}\right \|\qquad \text{is small}.
\end{equation*}
This task can be accomplished with $(k+p)$ matrix--vector products of $\mtx A$ and its adjoint, where $p$ is a small parameter, e.g. $p=5$.

First, we generate a randomized sketch of the matrix $\mtx A$ as 
\begin{equation}
\underset{m \times (k+p)}{\mtx Y}\ =\  \underset{m \times n\vphantom{()}}{\mtx A}\ \underset{n \times (k+p)}{\mtx \Omega}, \qquad \mtx \Omega \sim \mathcal N(\mtx 0,\mtx I)
\label{eq:Asketchlr}
\end{equation}
where $\mtx \Omega$ is a standard Gaussian matrix. With high probability, the columns of $\mtx Y$ span the dominant column space of $\mtx A$ \cite[Sec.~10]{2020_acta_martinsson_tropp}.
Then an approximate factorization of $\mtx A$ can be computed as
\begin{equation}
\underset{m \times k}{\mtx Q} = \texttt{orth}_k(\mtx Y), \qquad \mtx B := \mtx Q^*\ \mtx A,
\label{eq:Qorth}
\end{equation}
where $\texttt{orth}_k$ is defined in equation (\ref{eq:orthk}). The computation for $\mtx B$ requires the action of the adjoint of $\mtx A$. 

Although the computation involves randomization, the produced basis in (\ref{eq:Qorth}) using random sketch (\ref{eq:Asketchlr}) is within a factor that is polynomial in $k$ and $p$ of the optimal error. For $k \ge 2, p \ge 4$ and $k+p \le \min(m,n)$, the probability that 
\begin{equation*}\left\| \mtx A - \mtx Q \mtx Q^* \mtx A \right\| \le \left( 1 + 6 \sqrt {(k+p) p \log p} \right) \sigma_k + 3 \sqrt{(k+p) \sum_{j>k} \sigma_j^2}
\end{equation*} approaches
1 at a rate faster than any exponential function with increasing $p$ \cite[Corollary 10.9]{2011_martinsson_randomsurvey}. For many PDE operators, where the singular values decay exponentially, the gap from optimality is essentially bounded by a polynomial factor that depends only on $k$ and $p$.
Randomized sketching methods can be used to 
construct a wide range of low-rank decompositions, 
including the interpolative decomposition of Section \ref{sec:ID}. 
These methods are especially useful in the black-box setting because the matrix $\mtx A$ is only accessed through its action on vectors.

\subsection{Block elimination matrices}
\label{sec:block_elim}

Consider a block matrix of the form
\begin{equation}
\mtx{A}
=
\left(\begin{array}{rrr}
\mtx{A}_{11} & \mtx{A}_{12} &       \\
\mtx{A}_{21} & \mtx{A}_{22} & \mtx{A}_{23} \\
      & \mtx{A}_{32} & \mtx{A}_{33}
\end{array}\right).
\label{eq:Ablockmat}
\end{equation}
If $\mtx{A}_{11}$ is nonsingular, we can ``decouple'' it from the other blocks via
one step of block Gaussian elimination, by multiplying $\mtx A$ in (\ref{eq:Ablockmat}) on the left and right with matrices $\mtx L$ and $\mtx U$ as
\begin{equation}
\label{eq:LAU}
\mtx L\ \mtx A\ \mtx U =
\left(\begin{array}{ccc}
\mtx{A}_{11} &  &       \\
& \mtx{S}_{22} & \mtx{A}_{23} \\
  & \mtx{A}_{32} & \mtx{A}_{33}
  \end{array}\right),
\end{equation}
where the matrices $\mtx L$ and $\mtx U$ are unit-triangular matrices
\begin{equation}
\mtx{L}\ = \
\left(\begin{array}{ccc}
\mtx{I} &  &   \\
-\mtx{A}_{21}\mtx{A}_{11}^{-1} & \mtx{I} & \hphantom{\mtx{0}_{23}}  \\
 &  & \mtx{I}
\end{array}\right)
\qquad\mbox{and}\qquad
\mtx{U}\ =\
\left(\begin{array}{ccc}
\mtx{I} & -\mtx{A}_{11}^{-1}\mtx{A}_{12} &  \hphantom{\mtx{0}_{13}} \\
 & \mtx{I} &   \\
& & \mtx{I}
\end{array}\right),
\label{eq:LUdef}
\end{equation}
and the submatrix $\mtx{S}_{22}$ in (\ref{eq:LAU}) is the 
Schur complement
\begin{equation*}
\mtx{S}_{22} = \mtx{A}_{22} - \mtx{A}_{21}\mtx{A}_{11}^{-1}\mtx{A}_{12}.
\end{equation*}
Block-elimination matrices of the form in (\ref{eq:LUdef}) are simple to invert by toggling the sign of the off-diagonal block.
When $\mtx A$ is symmetric, the factor $\mtx U = \mtx L^*$ and symmetry is preserved in (\ref{eq:LAU}).

\subsection{The interpolative decomposition}
\label{sec:ID}
Let $\mtx{A}_{\msf I \msf J}$ be a matrix subblock of size $m\times n$
and approximate rank $k$. 
The \textit{interpolative decomposition (ID)} of $\mtx{A}_{\msf I \msf J}$
 is a low-rank factorization where a subset
of $k$ rows (or columns) is used to span the row space (or column space) of $\mtx A$. To be precise, for the \textit{row ID}, we find a partition
$\msf I = \msf R \cup \msf S$, such that the matrix admits a low-rank decomposition
\begin{equation}\label{eq:id_factorization}
\mtx A_{\msf I \msf J}\ = \begin{pmatrix} \mtx A_{\msf R \msf J}\\ \mtx A_{\msf S \msf J} \end{pmatrix}
\approx
\begin{pmatrix}\underset{(m-k) \times k}{\mtx T}\\ \underset{k \times k}{\mtx I}
\end{pmatrix}
\underset{k \times n}{\mtx A_{\msf S \msf J}}.
\end{equation}
We use the compact notation
\begin{equation}
\texttt{id}(\mtx A_{\msf I \msf J}) =
\left[
\underset{m - k}{\msf R} \cup\ \underset{k}{\msf S},\
\underset{(m-k) \times k}{\mtx T}
\right],
\label{eq:id_def}
\end{equation}
for an index partition and interpolation matrix $\mtx T$ that satisfy (\ref{eq:id_factorization}).

Finding the optimal $k$ rows is a combinatorially hard problem. However, the strong rank-revealing QR factorization \cite{gu1996} is guaranteed to produce a near-optimal factorization. In practice, the standard pivoted QR with a greedy approach performs well.
Although the error in an approximate low-rank $k$ interpolative decomposition can, in theory, be significantly larger than that obtained by truncating a singular value decomposition, the practical error is usually modest when the singular values of the input matrix decay at a reasonable rate, as is often the case for PDE problems.
For numerical stability, it is desirable that the matrix \(\mtx{T}\) be well-conditioned, which in practice means keeping its entries small. It has been demonstrated that one can always choose the set \(\mathsf{S}\) such that every entry of \(\mtx{T}\) has a modulus bounded by one, and practical algorithms exist to ensure that these entries remain modest \cite{2005_martinsson_skel,gu1996,martinsson2008id}.

In scenarios where the matrix is not easily accessible, randomized methods \cite{2021_martinsson_dong_CUR,2007_martinsson_PNAS,2011_martinsson_random1,2014_martinsson_CUR} provide an efficient means of computing the skeleton set and the interpolation matrix \(\mtx{T}\). Suppose that one computes the ID of the sketch
\begin{equation}
\underset{ m \times (k+p)}{\mtx Y_{\msf I}}\ =\ 
\underset{ \vphantom{()} m \times n}{\mtx A_{\msf I\msf J}}\ 
\underset{ n \times (k+p)}{\mtx \Omega_{\vphantom {B}}},
\qquad
\texttt{id}( \mtx Y_{\msf I}) = \left[\underset{m - k}{\msf R} \cup\ \underset{k}{\msf S},\
\underset{(m-k) \times k}{\mtx T} \right].
\label{eq:YIsketch}
\end{equation}
Then, by some simple observations \cite[Section 13.3]{2020_acta_martinsson_tropp}, the information needed for the ID of \(\mtx{A}_{\mathsf{IJ}}\) is contained in the ID of the sketch $\mtx Y_{\msf I}$. Because the sketch approximately spans the dominant column space of $\mtx A$, the matrix admits
the low-rank decomposition 
$$
\underset{m \times n\vphantom{)}}{\mtx A_{\msf I \msf J}}\
 \approx \ \underset{m \times (k+p)}{\mtx Y_{\msf I}}\
  \underset{(k+p) \times n}{\mtx B_{\vphantom{B}}}
$$ 
for some matrix $\mtx B$. Consequently,
\begin{equation*}
\begin{pmatrix} \mtx A_{\msf R \msf J}\\ \mtx A_{\msf S \msf J}\end{pmatrix}\
\approx\ \begin{pmatrix} \mtx Y_{\msf R}\\ \mtx Y_{\msf S} \end{pmatrix} \mtx B\ 
\approx\ \begin{pmatrix} \mtx T\\\mtx I \end{pmatrix}\ \mtx Y_{\msf S}\ \mtx B\
=\ \begin{pmatrix} \mtx T\\ \mtx I \end{pmatrix} \mtx A_{\msf S \msf J}.
\end{equation*}
Thus, the key information needed for the ID of $\mtx A_{\msf {IJ}}$ is encapsulated in the ID of the sketch $\mtx Y_{\msf I}$ defined in (\ref{eq:YIsketch}). Using deterministic methods, the complexity of the ID is $\mathcal O(m n^2)$. Randomized methods reduce the complexity to $\mathcal O(m n k + mk^2)$ when using Gaussian random matrices. Likewise, to compute a representative set of columns, the algorithm described
can be executed by sketching the transpose of the matrix. 
\section{LU Factorization of Hierarchical Matrices using Strong Admissibility}
\label{sec:skel_algos}

In this section, we describe the construction of an invertible factorization for $\mathcal{H}^2$ matrices under the strong admissibility condition, largely following the presentation in \cite{2017_ho_ying_strong_RS}. As an example, consider the matrix $\mtx A$ defined by the Green’s function $G$ for the Laplace equation:
\begin{align*}
\begin{split}
\mtx A_{ij} &= G(\pxx_i, \pxx_j),\ i \neq j\\ 
&\text{where}\ G(\pxx_i, \pxx_j) = \begin{cases}
\log(\|\pxx_i - \pxx_j\|),& \pxx \in \mathbb{R}^2 \\
{(\| \pxx_i - \pxx_j\|)}^{-1},& \pxx \in \mathbb{R}^3.
\end{cases}
\end{split}
\end{align*}
An appropriate quadrature correction is applied on the diagonal to ensure that $\mtx A$ is invertible. Matrices of this type may arise from the discretization of an integral equation or from covariance matrices in statistics. For simplicity in this section, we assume that the matrix entries are readily accessible and that the matrix is stored densely. However, when $\mtx A$ originates from the discretization of an integral equation with an explicit formula for evaluating its entries, efficient algorithms that avoid forming the matrix densely are detailed in \cite{2017_ho_ying_strong_RS,sushnikova2023fmm}.

The algorithms rely on organizing a set of points $\{\pxx_j\}_{j=1}^N$ into a hierarchical structure, either a {quadtree} in two dimensions or an {octree} in three dimensions, depending on whether the points lie in 2D or 3D space.
Formally, we construct a tree $\tree$ in which each node, or {box} $B$, contains a subset of the points. Initially, all points are contained in a single box called the {root}. The root box is recursively subdivided into $2^d$ child boxes, where $d$ is the dimension of the space. This subdivision continues until each box contains no more than $m$ points, where $m$ is a user-specified threshold.
We refer to a box that has children as a {tree box}, and a box with no children as a {leaf}. The {depth} of a box is its distance (in number of edges) from the root box. The collection of all boxes at depth $\ell$ is called {level} $\ell$ of the tree. Thus, level 0 consists of just the root box, level 1 contains its $2^d$ children, and higher levels correspond to progressively finer partitions of the domain.
The {depth} of the tree, denoted $L$, is the maximum depth of any box, and is approximately given by
$L \approx \log_2 \left( \frac{N}{m} \right)$.

Two boxes in the tree are said to be {adjacent} if they share a face, edge, or corner. In cases where the point distribution is non-uniform, only those boxes containing more than $m$ points are further subdivided, resulting in an {adaptive} tree structure. We assume that such adaptive trees satisfy a {2:1 balance condition}, meaning that any two adjacent leaf boxes differ in depth by at most one. This constraint limits the number of adjacent boxes and helps maintain computational efficiency.
For a given box $B$, we distinguish between its {neighbor} boxes---those that are adjacent---and its {far-field} boxes---those that are well-separated. In a slight abuse of notation, we also use $\msf B$ to denote the set of indices corresponding to points within box $B$. Similarly, $\msf N$ refers to the indices of points in the neighboring boxes of $B$, while $\msf F$ refers to the indices of points in its far-field boxes.

\subsection{Strong Skeletonization}
\label{sec:strong_skel}

For a dense matrix $\mtx A$, the procedure described in this section introduces and exploits sparsity in a modified system via a sequence of linear transformations applied to $\mtx A$. Assume that the matrix is tessellated according to the geometry, following the index order $[\msf B, \msf N, \msf F]$, so that
\begin{equation}
\mtx A = 
\begin{pmatrix}
\mtx A_{\msf B \msf B} & \mtx A_{\msf B \msf N} & \mtx A_{\msf B \msf F} \\
\mtx A_{\msf N \msf B} & \mtx A_{\msf N \msf N} & \mtx A_{\msf N \msf F} \\
\mtx A_{\msf F \msf B} & \mtx A_{\msf F \msf N} & \mtx A_{\msf F \msf F} \\
\end{pmatrix},\quad
\text{where}\quad
\begin{minipage}{0.4\textwidth}
\begin{tabular}{l}
$\msf B$: target box indices for box $B$,\\
$\msf N$: near field indices,\\
$\msf F$: far field indices.
\end{tabular}
\end{minipage}
\label{eq:labels_matrix}
\end{equation}
Consider the subblock $\mtx A_{\msf F \msf B}$ corresponding to the interaction between a box and its far field points. Because the points are well-separated, the block is numerically low-rank and admits the interpolative decomposition of Section \ref{sec:ID} as
\begin{equation}
\underset{m\vphantom{()}}{\msf B} = \underset{m-k}{\msf R} \cup\ \underset{k\vphantom{()}}{\msf S}\ , 
\qquad 
\begin{pmatrix} \mtx A_{\msf F \msf R}& \mtx A_{\msf F \msf S} \end{pmatrix}
\approx 
 \mtx A_{\msf F \msf S}\begin{pmatrix}\underset{(m-k) \times k}{\mtx T}& \underset{k \times k}{\mtx I} \end{pmatrix}.
\label{eq:skel_lowrank}
\end{equation}
Likewise, an analogous statement holds for $\mtx A_{\msf B \msf F}$. It is often convenient to choose a skeleton set and corresponding interpolation matrix which is applicable for both $\mtx A_{\msf B \msf F}$ and $\mtx A_{\msf F \msf B}$. This can be done by computing the ID of the concatenation so that
\begin{equation}
\begin{pmatrix}
\mtx A_{\msf F \msf R} & \mtx A_{\msf F \msf S}\\
\mtx A^*_{\msf R \msf F} & \mtx A^*_{\msf S \msf F}
\end{pmatrix} \approx \begin{pmatrix}
\mtx A_{\msf F \msf S}\\
\mtx A^*_{\msf S \msf F}
\end{pmatrix}
\begin{pmatrix}
\mtx T & \mtx I
\end{pmatrix}.
\label{eq:skel_lowrank_simult}
\end{equation}
Instead of using low-rank decompositions, which project onto a $k$-dimensional subspace, we aim to remain in the full $m$-dimensional space while introducing sparsity into the system. In this approach, an equivalent formulation of (\ref{eq:skel_lowrank}) is given by
\begin{equation*}
\left( \begin{array}{cc}{\mtx A}_{\msf F \msf R} & {\mtx A}_{\msf F \msf S}\end{array}\right) \left(\begin{array}{cc}\mtx I & \\[1mm]-\mtx T & \mtx I \end{array} \right)
=
\left( \begin{array}{cc}\mtx A_{\msf F \msf R} -  \mtx A_{\msf F \msf S} \mtx T& {\mtx A}_{\msf F \msf S}\end{array}\right)
\approx
\left( \begin{array}{cc}\mtx 0 & {\mtx A}_{\msf F \msf S}\end{array}\right).
\end{equation*}
To apply this local transformation globally, we embed it into the full matrix using sparsifying matrices $\mtx E$ and $\mtx F$ acting on all indices. For an appropriate permutation 
$[\msf R, \msf S, \msf N, \msf F]$ of \eqref{eq:labels_matrix}, we define
\begin{equation}
\mtx E = 
\left( \begin{array}{cc|cc} \mtx I & -\mtx T^* \\& \mtx I \\\hline && \mtx I\\ &&&\mtx I
\end{array}\right),\qquad
\mtx F = 
\left( \begin{array}{cc|cc} \mtx I & \\- \mtx T& \mtx I \\\hline && \mtx I\\ &&&\mtx I
\end{array}\right).
\label{eq:orth_sparsifying}
\end{equation}
Applying $\mtx E$ and $\mtx F$ on the left and right of $\mtx A$, yields the sparsified system
\begin{equation}
\mtx E\ \mtx A\ \mtx F \approx 
\left(\begin{array}{cc|cc}
{\mtx X}_{\msf R \msf R} & {\mtx X}_{\msf R \msf S} & {\mtx X}_{\msf R \msf N} & \\[1mm]
{\mtx X}_{\msf S \msf R} & {\mtx A}_{\msf S \msf S} & {\mtx A}_{\msf S \msf N} &  {\mtx A}_{\msf S \msf F}\\ \hline
{\mtx X}_{\msf N \msf R} & {\mtx A}_{\msf N \msf S} & {\mtx A}_{\msf N \msf N} & {\mtx A}_{\msf N \msf F}\\[1mm]
& {\mtx A}_{\msf F \msf S} & {\mtx A}_{\msf F \msf N} & {\mtx A}_{\msf F \msf F}
\end{array}\right).
\label{eq:EAFskelstrong}
\end{equation}
The interactions between $\msf R$ and the far field $\msf F$ are approximately zero due to \eqref{eq:skel_lowrank_simult}. 
Importantly, the use of the ID does not modify any entries associated with the skeleton subset $\msf S$, thereby preserving a physical interpretation of the system. This retention of structure enables the application of analytic compression techniques, such as proxy surfaces \cite{xing2020interpolative,ye2020analytical} or adaptive cross-approximation \cite{bebendorf2000approximation,2005_borm_grasedyck_hybrid_ACA,grasedyck2005adaptive} for discretized boundary integral equations.

\begin{figure}[!htb]
\resizebox{\textwidth}{!}{\newcommand{\annotatematrix}[1]{\begin{scope}[
      shift={(#1.south west)},      x={(#1.south east)},          y={(#1.north west)}]          \foreach \i/\lbl in {1/F, 2/N, 3/B, 4/N, 5/F, 6/F}{\node[font=\tikzannotsize,overlay] at ({(\i-0.5)/6}, 1.01) {$\msf \lbl$};
\node[font=\tikzannotsize,overlay] at (-0.01, {1-(\i-0.5)/6}) {$\msf \lbl$};
    }
  \end{scope}}

\newcommand{\annotategeom}[2]{\begin{scope}[shift={(#1.south west)},
                x={(#1.south east)}, y={(#1.north west)}]
    \foreach \i/\lbl in {1/{\msf F}, 2/{\msf N}, 3/{#2}, 4/{\msf N}, 5/{\msf F}, 6/{\msf F}}{\node[font=\tikzannotsize,overlay] at ({(\i-.5)/6}, +0.01) {$\lbl$};
    }
  \end{scope}}

\begin{tikzpicture}

\node (m1) {\includegraphics[width=35mm]{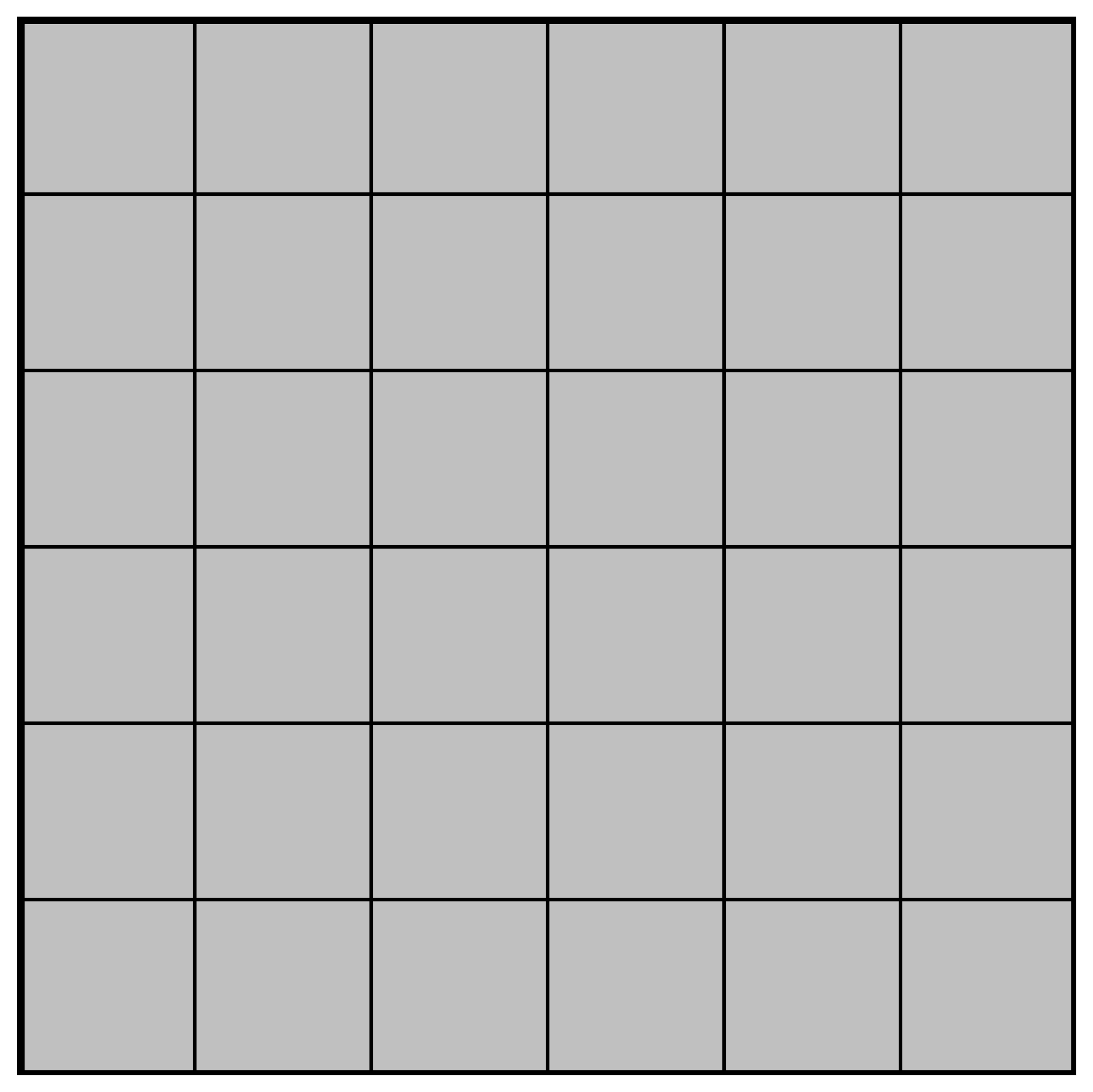}};

\node (m2) [right= 15mm of m1] {\includegraphics[width=35mm]{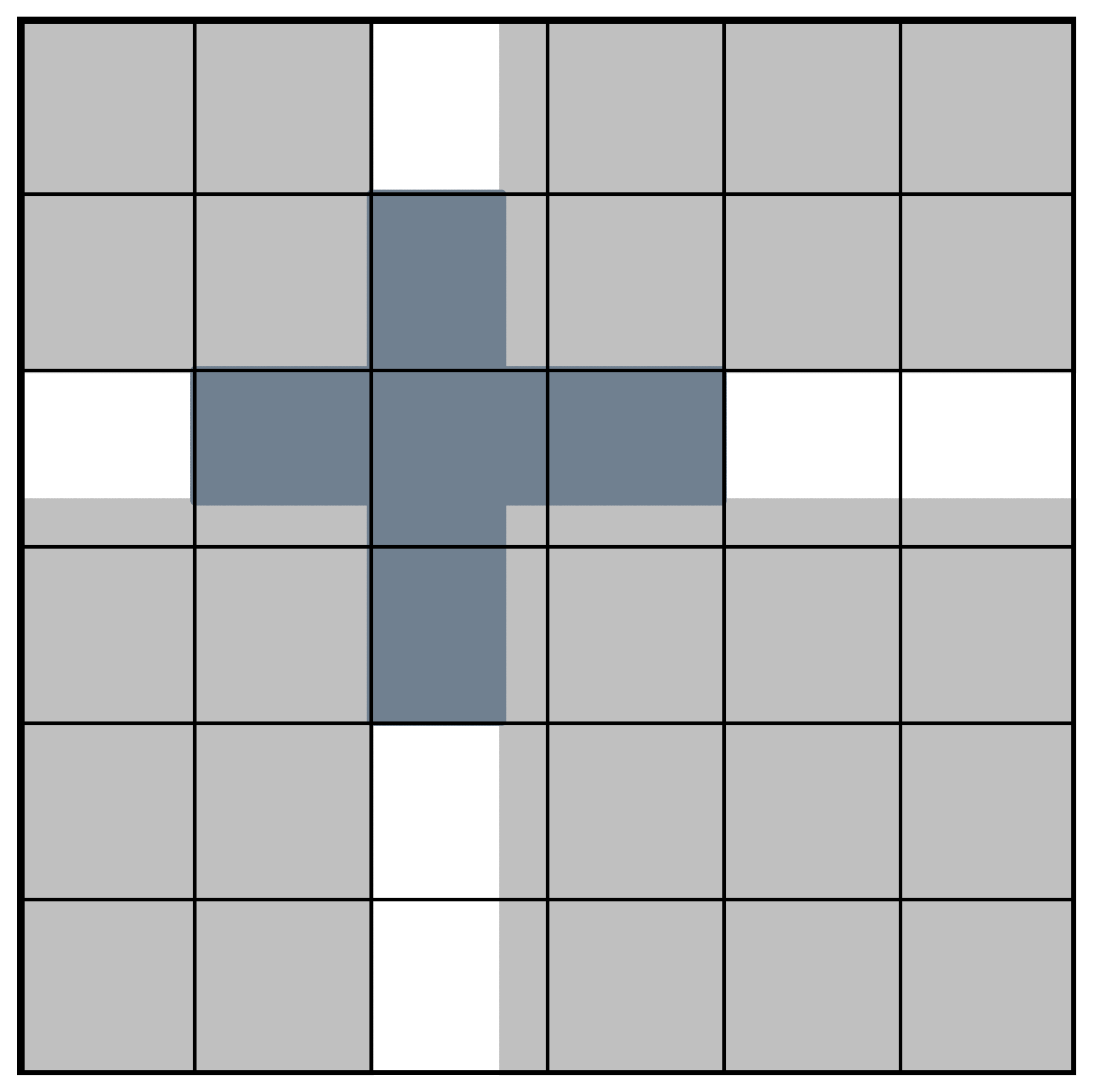}};

\node (m3) [right= 15mm of m2] {\includegraphics[width=35mm]{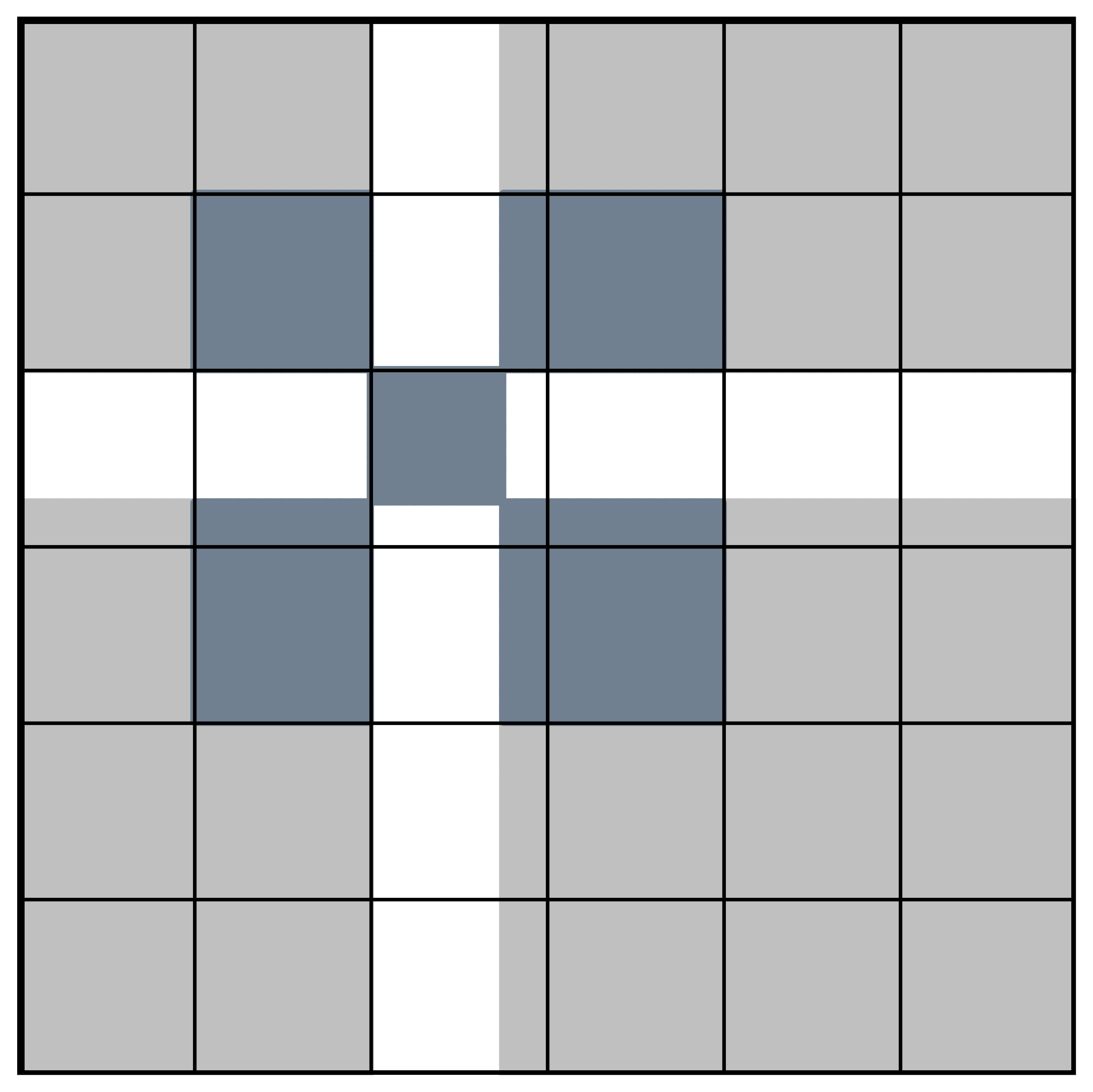}};

\node (d1) [below=0mm of m1] {\includegraphics[width=35mm]{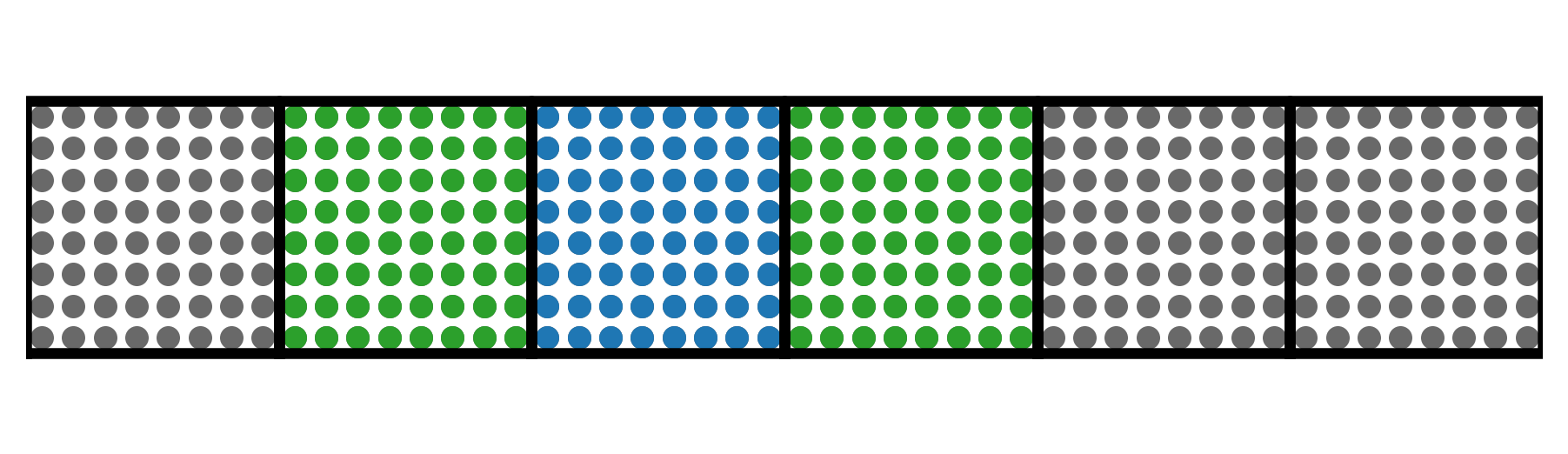}};
\node (d2) [below=0mm of m2] {\includegraphics[width=35mm]{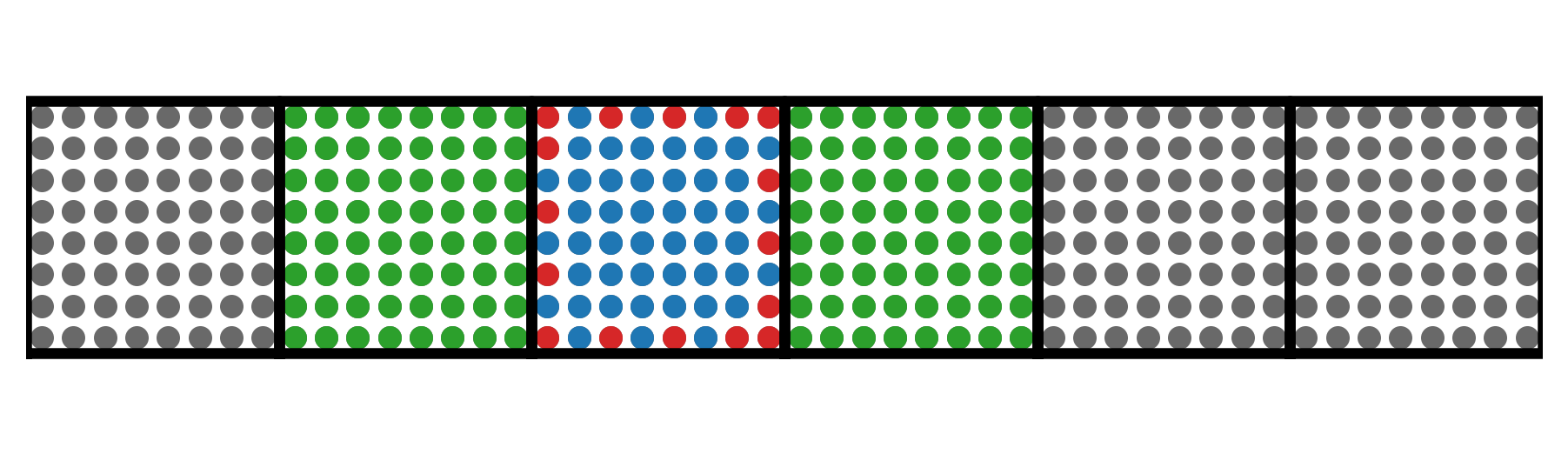}};
\node (d3) [below=0mm of m3] {\includegraphics[width=35mm]{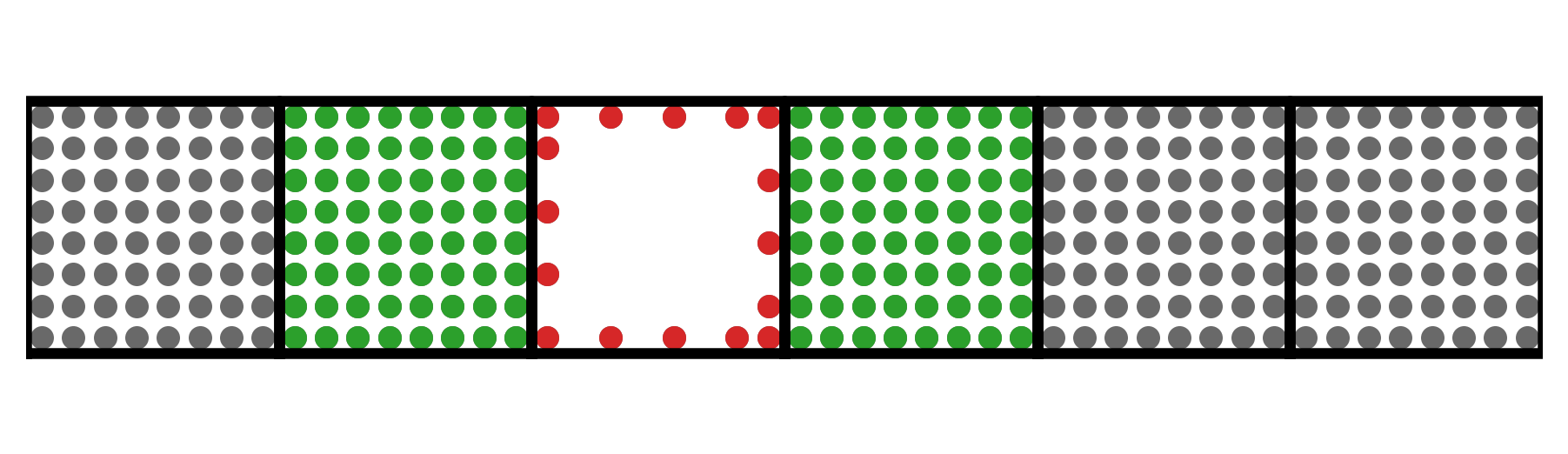}};

\draw[->, thick] (m1) -- node[below, midway, font=\tikzlabelsize] {$\msf{sparsify}$}    (m2);
\draw[->, thick] (m2) -- node[below, midway, font=\tikzlabelsize] {$\msf{diagonalize}$} (m3);

\draw[->, thick] (d1) -- (d2);
\draw[->, thick] (d2) -- (d3);

\annotatematrix{m1}

\annotategeom{d1}{\msf B}
\annotategeom{d2}{\msf R \cup \msf S}
\annotategeom{d3}{\msf S}

\end{tikzpicture}
 }
\caption{Illustration of strong skeletonization for a single box $B$. Top row: the matrix structure at each stage (with modified entries in blue). Bottom row: corresponding pseudo-1D geometry and index ordering. 
``Sparsify'' applies the ID to far-field interactions and identifies redundant 
indices $\mathsf R$. ``Diagonalize'' then eliminates the remaining interactions 
involving $\mathsf R$ via block elimination.
}
\label{fig:sparsify_then_elim}
\end{figure}

Now that the matrix is sparser, the objective is to \textit{diagonalize} the interactions between the redundant box indices $\msf R$ and the remaining indices of the matrix (i.e., to eliminate the residual sparse interactions). This is achieved using standard sparse block-elimination techniques, as described in Section \ref{sec:block_elim}. We next introduce matrices $\mtx L, \mtx U$ as
\begin{equation}
\def\arraystretch{1.2}
\mtx L = \left( \begin{array}{cc|cc}  
\mtx I & \hphantom{\mtx XY} & \hphantom{\mtx XY} & \hphantom{\mtx X}\\[1mm]
-\mtx X_{\msf S \msf R} \mtx X^{-1}_{\msf R\msf R } & \mtx I\\ \hline
-\mtx X_{\msf N \msf R} \mtx X^{-1}_{\msf R \msf R} &        & \mtx I\\[1mm]
&  &  & \mtx I
\end{array}\right),\quad
\mtx U = \left( \begin{array}{cc|cc} 
\mtx I & -\mtx X^{-1}_{\msf R \msf R} \mtx X_{\msf R \msf S}  & -\mtx X^{-1}_{\msf R \msf R }  \mtx X_{\msf R \msf N} & \hphantom{\mtx XY}\\[1mm]
& \mtx I\\ \hline
 & & \mtx I\\[1mm]
& &  & \mtx I
\end{array}\right),
\label{eq:LUstrong}
\end{equation}
where $\mtx X = \mtx E \mtx A \mtx F$ and the submatrices $\mtx X_{\msf R\msf R }, \mtx X_{\msf R\msf S }$, etc., are defined in \eqref{eq:EAFskelstrong}.
Then applying $\mtx L$ and $\mtx U$ to the left and right, respectively, yields the diagonalized matrix $\widetilde{\mtx A}$
\begin{equation}
\widetilde{\mtx A}(\mtx A; B) = \mtx L\ \left(\mtx E\ \mtx A\ \mtx F \right)\ \mtx U
 \approx 
\left(\begin{array}{cc|cc}
{\mtx X}_{\msf R \msf R} &  \\[1mm]
 & {\mtx X}_{\msf S \msf S} & {\mtx X}_{\msf S \msf N} &  {\mtx A}_{\msf S \msf F}\\ \hline
& {\mtx X}_{\msf N \msf S} & {\mtx X}_{\msf N \msf N} & {\mtx A}_{\msf N \msf F}\\[1mm]
& {\mtx A}_{\msf F \msf S}  & {\mtx A}_{\msf F \msf N} & {\mtx A}_{\msf F \msf F}
\end{array}\right),\label{eq:Atildedef}
\end{equation}
where modifications to the original matrix entries are introduced to the near-neighbor interactions.

We define \textit{diagonalization} matrices for a box $B$ as $\mtx V, \mtx W$, which are compositions of the skeletonization and block-elimination transformations in \eqref{eq:orth_sparsifying} and \eqref{eq:LUstrong}:
\begin{align}
\begin{split}
\mtx V^{-1}\ &=\ \mtx L\ \mtx E,\\
\mtx W^{-1}\ &=\ \mtx F\ \mtx U,\\
\widetilde{\mtx A}(\mtx A; B) \ &\approx\ \mtx V^{-1}\ \mtx A\ \mtx W^{-1},
\end{split}
\label{eq:tildeA_diag}
\end{align}
where $\widetilde{\mtx A}$ is diagonalized with respect to the redundant indices $\msf R \subseteq \msf B$ of box $B$.
Because the matrices $\mtx V, \mtx W$ are invertible, we can also express the decomposition as a sparse factorization of $\mtx A$:
\begin{equation*}
\mtx A\ \approx\ \mtx V\ \widetilde{\mtx A}\ \mtx W,
\end{equation*}
with $\widetilde{\mtx A}$ defined in \eqref{eq:Atildedef} and its corresponding sparsity pattern in \eqref{eq:tildeA_diag}.

The basic idea of the subsequent algorithmic steps is to repeat the sparsify-then-eliminate procedure for each box in the hierarchical tree decomposition. Because the interactions with respect to $\msf R$ have been diagonalized, these subsequent linear operations do not introduce additional nonzero interactions for $\msf R$. We make two important remarks regarding the diagonalized matrix. 

\begin{remark}[Well-conditioning of the sparsifying matrices] Applying the sparsifying matrices $\mtx E$ and $\mtx F$ on the left and right does modify the condition number of the subblock $\mtx X_{\msf R \msf R}$, however, the effect is modest. In fact, owing to the bounded entries of $\mtx T$ (as discussed in Section \ref{sec:ID}), both $\mtx E$ and $\mtx F$ remain well-conditioned, and the condition number of $\mtx A_{\msf R \msf R}$ and that of the modified subblock exhibit only minor differences in practice. 
\end{remark}

\begin{remark}[Modified near-field interactions] The diagonalized matrix $\widetilde{\mtx A}$ introduces an additive term of dense interactions between the near-neighbors $\msf N$, namely,
\[
\mtx X_{\msf N \msf N} = \mtx A_{\msf N \msf N} - \mtx X_{\msf N \msf R} \mtx X_{\msf R \msf R}^{-1} \mtx X_{\msf R \msf N}.
\]
Consider two boxes, $\alpha, \beta$ which neighbor $B$. Although $\alpha, \beta$ have a common neighbor, they may in fact be well-separated from each other. For a concrete example, see Figure \ref{fig:sparsify_then_elim}.  Therefore, this additional term could in principle impact the existing low-rank spectrum. In practice, however, the modified interactions remain compressible as low rank. In settings where the matrix entries are explicitly accessible, these modified entries are stored in an intermediate representation and then recompressed in later stages of the algorithm.
\label{remark:additive_lowrank}
\end{remark}

\subsection{Strong Recursive Skeletonization}
\label{sec:strong_rec_skel}

In this section, we describe a recursive algorithm obtained by composing the diagonalization matrices
defined in Section~\ref{sec:strong_skel} for the multilevel tree decomposition $\tree$. At a high level, the
algorithm proceeds level by level from the finest to the coarsest. At each level, we apply the procedure
in Section~\ref{sec:strong_skel} to each box and update the remaining active degrees of freedom.

Recall that the root box containing all points is on level~0, and that levels are labeled by their depth
(i.e., their distance from the root). The finest level is labeled $L$, and necessarily contains only leaf boxes.
We first describe how the diagonalization matrices are composed for all boxes on level $L$.
We then generalize the procedure to all levels. Following previous works
\cite{2017_ho_ying_strong_RS,sushnikova2023fmm}, the boxes are labeled $B_1, B_2, \dots$ according to the
order in which they are diagonalized.

Suppose that for the first box $B_{1}$ on level $L$, we compute diagonalization matrices $\mtx V_1, \mtx W_1$
so that a subset of indices $\msf R_{1} \subseteq \msf B_{1}$ is diagonalized and decoupled from the rest of the dense system as
\begin{equation*}\widetilde{\mtx A}(\mtx A; B_{1})\ \approx \ \mtx V_{1}^{-1}\ \mtx A\ \mtx W_1^{-1}.\end{equation*}
As noted in Remark \ref{remark:additive_lowrank}, these transformations introduce additive terms into the system, and the next step must diagonalize the modified system $\widetilde{\mtx A}$. Continuing with the second box $B_2$, we update the modified system as
\begin{equation*}
\widetilde{\mtx A}(\mtx A; B_1, B_2)\ =\ \mtx V_{2}^{-1}\ \widetilde{\mtx A}(\mtx A; B_1)\ \mtx W_{2}^{-1}.
\end{equation*}
Importantly, the matrices $\mtx V_2^{-1}, \mtx W_2^{-1}$ do not affect previously diagonalized indices, and
the redundant indices $\msf R_1 \subseteq \msf B_1$ remain diagonalized.

\begin{figure}[!htb!]
\resizebox{0.975\textwidth}{!}{

\newcommand{\annotatematrix}[1]{\begin{scope}[
      shift={(#1.south west)},      x={(#1.south east)},          y={(#1.north west)}]          \foreach \i in {1, 2, 3, 4, 5, 6, 7, 8}{\node[font=\tikzannotsize,overlay] at ({(\i-0.5)/8}, 1.0) {$\msf B_\i$};
\node[font=\tikzannotsize,overlay] at (-0.0, {1-(\i-0.5)/8}) {$\msf B_\i$};
    }
  \end{scope}}

\newcommand{\annotatematrixcoarsen}[5]{\begin{scope}[shift={(#1.south west)},
                x={(#1.south east)},   y={(#1.north west)}]   

    \foreach \i in {1,...,#2}{\node[font=\tikzannotsize,overlay]
        at ({#4 + (\i - 0.5)*((1 - #4)/#2)- 0.01}, {1 - #4 + #5 + 0.01})
        {$\msf B_{\number\numexpr\i+#3\relax}$};

      \node[font=\tikzannotsize,overlay]
        at ({#4 - #5 - 0.02}, {(1 - #4) - (\i - 0.5)*((1 - #4)/#2)})
        {$\msf B_{\number\numexpr\i+#3\relax}$};
    }

  \end{scope}}
\newcommand{\annotategeom}[2]{\begin{scope}[shift={(#1.south west)},
                x={(#1.south east)}, y={(#1.north west)}]

    \foreach \i in {1,...,8}{\node[font=\tikzannotsize,overlay] at ({(\i-.5)/8}, 0.2)
        {$\msf{\ifnum\i>\numexpr#2\relax B\else S\fi}_{\i}$};
    }\end{scope}}

\newcommand{\annotategeomcoarsen}[3]{\begin{scope}[shift={(#1.south west)},
                x={(#1.south east)}, y={(#1.north west)}]

    \foreach \i in {1,...,#2}{\node[font=\tikzannotsize,overlay] at ({(\i-.5)/#2}, 0.2)
        {$\msf B_{\number\numexpr\i+#3\relax}$};
    }\end{scope}}

\begin{tikzpicture}[node distance=0.5cm and 0.5cm, every node/.style={align=center}]

\node (mat12) {\includegraphics[width=60mm]{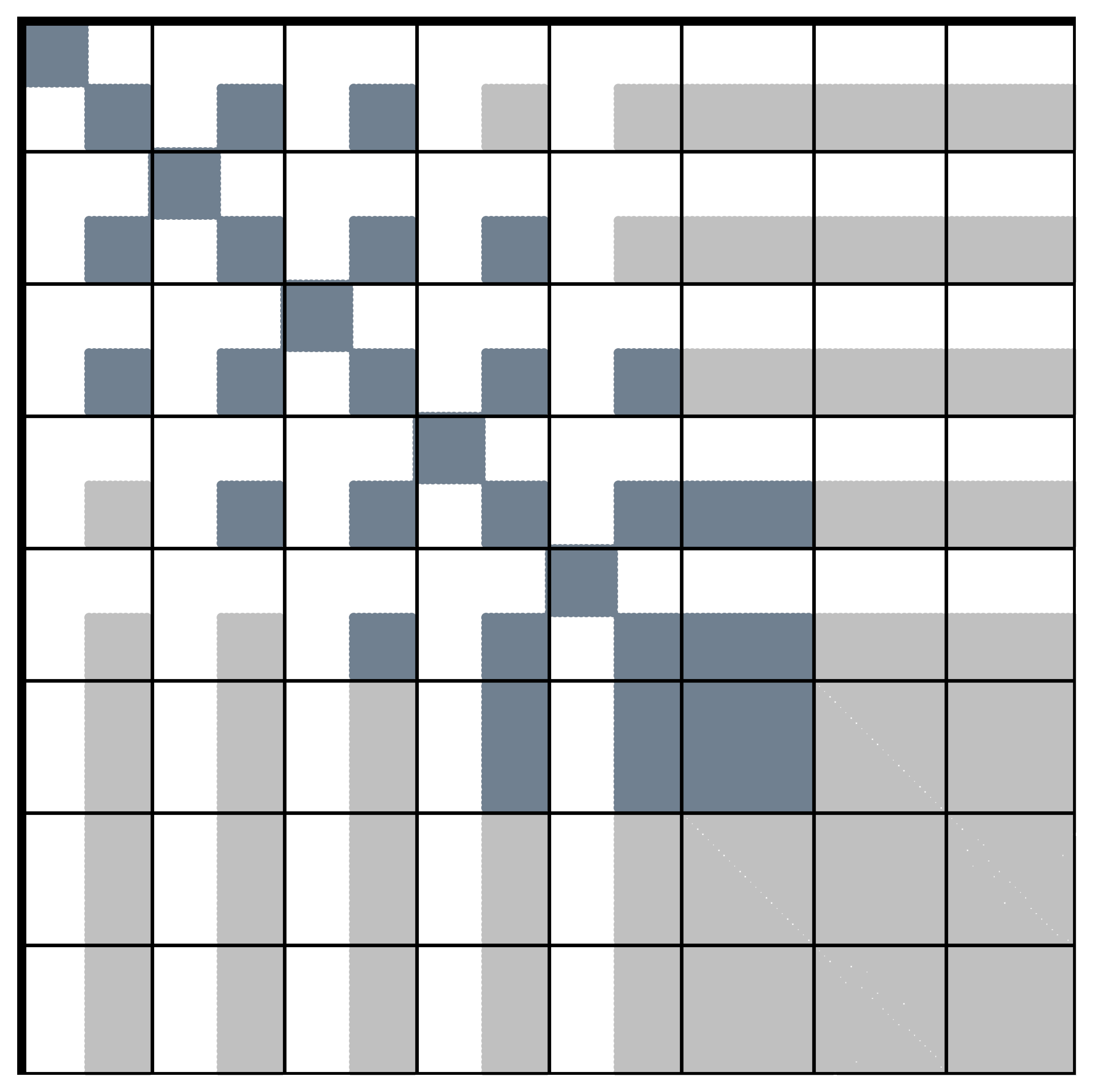}};
  \node[right=1.5cmof mat12] (level3mat) {\includegraphics[width=60mm]{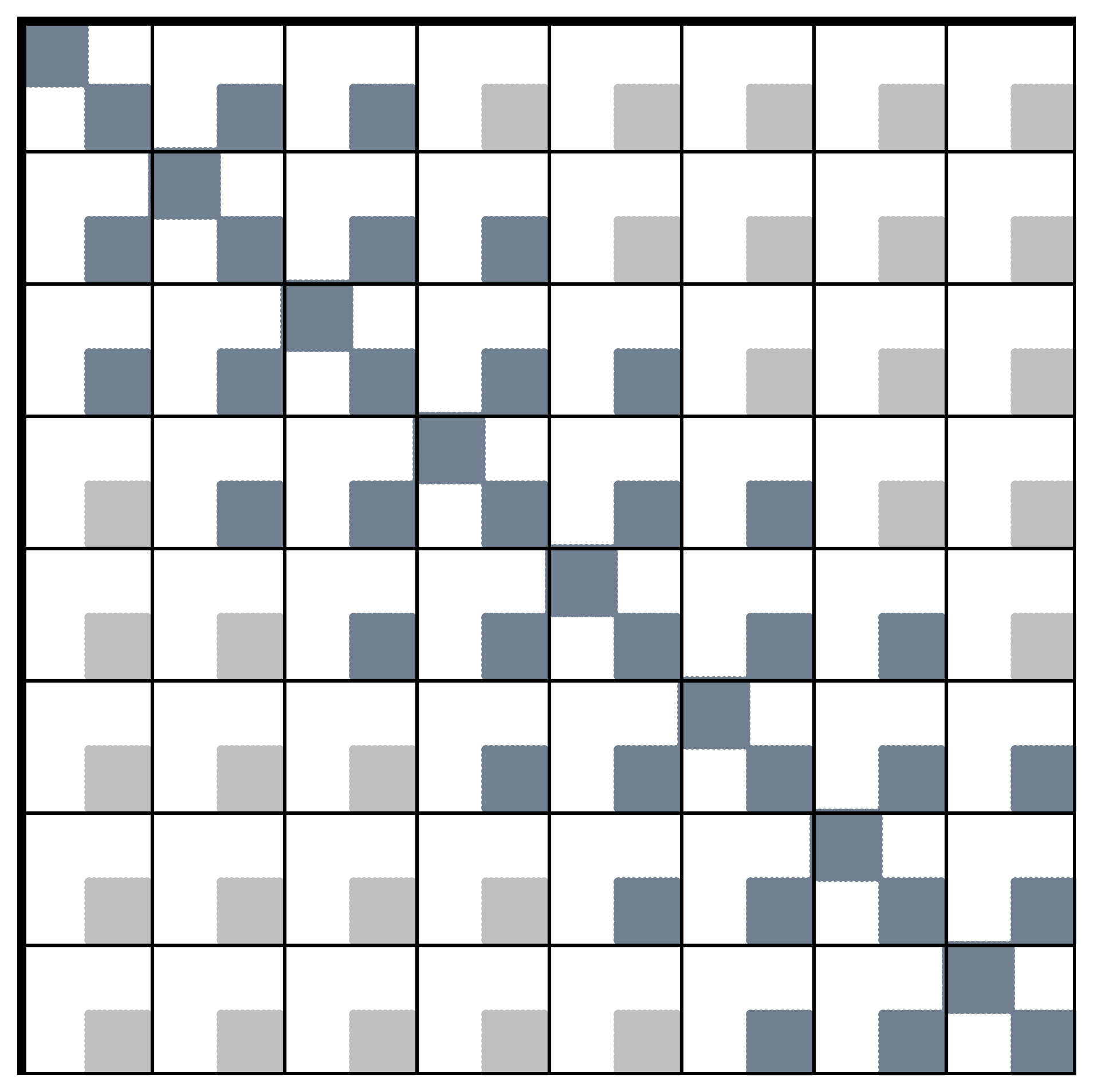}};
  \draw[->, thick, font=\tikzlabelsize] (mat12.east) -- node[below] {$\msf{diagonalize}$} (level3mat.west);
  
\node[below=-0.5cm of mat12] (geom12) {\includegraphics[width=60mm]{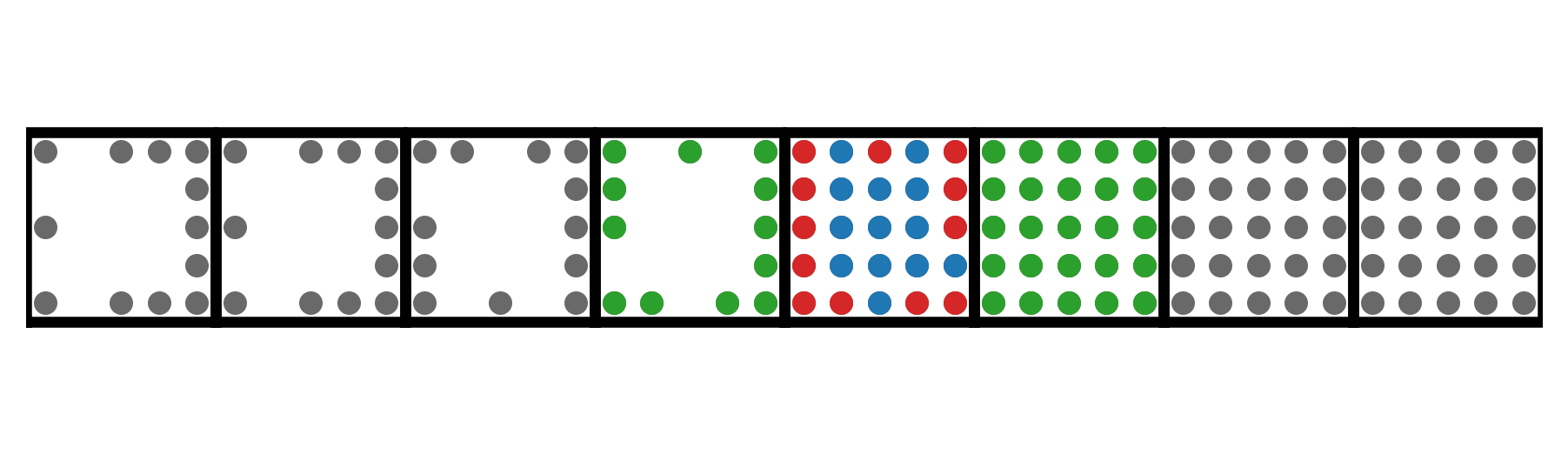}};
  \node[below=-0.5cm of level3mat] (level3geom) {\includegraphics[width=60mm]{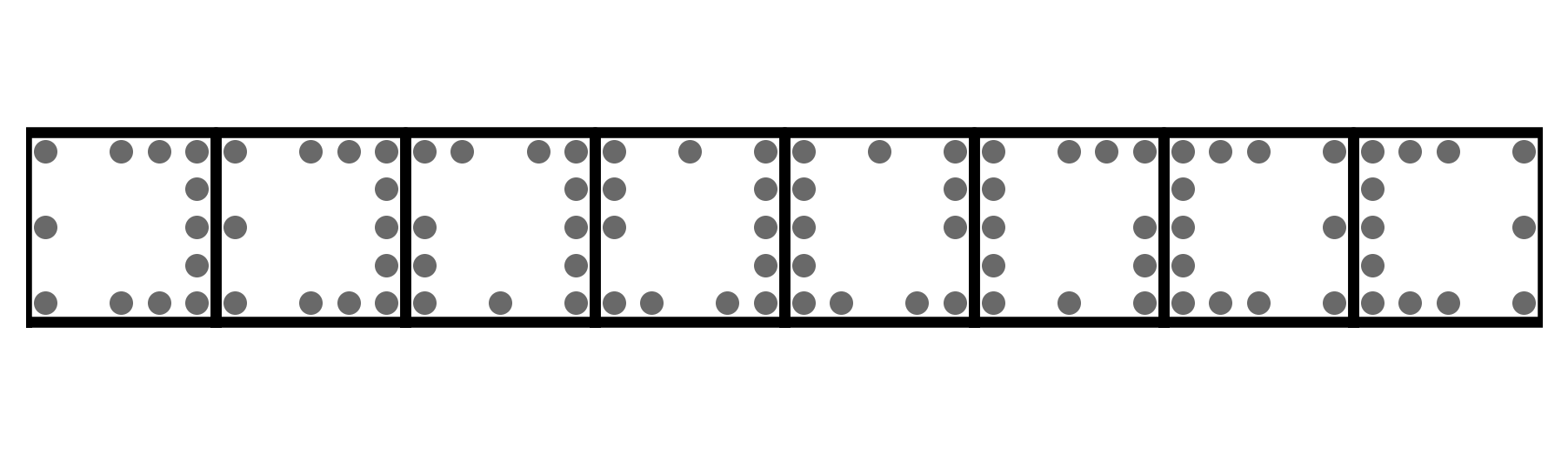}};
  \draw[->, thick, font=\tikzlabelsize] (geom12.east) -- (level3geom.west);

  \annotatematrix{mat12}
  \annotategeom{geom12}{4}
  \annotategeom{level3geom}{8}

\node[below= 2.5cm of mat12] (level2mat) {\includegraphics[width=60mm]{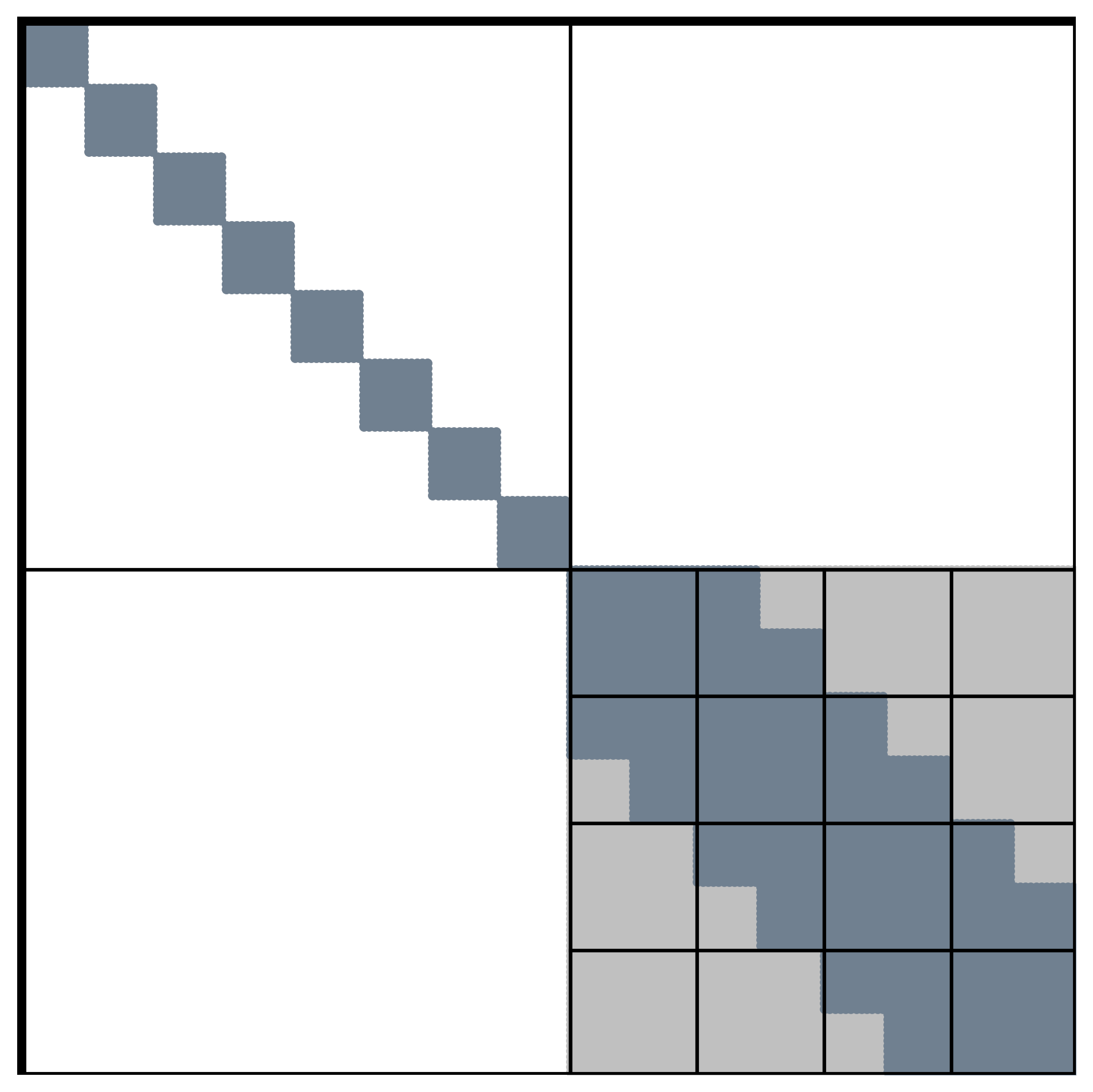}};
  \node[below= 2.5cm of level3mat] (level1mat) {\includegraphics[width=60mm]{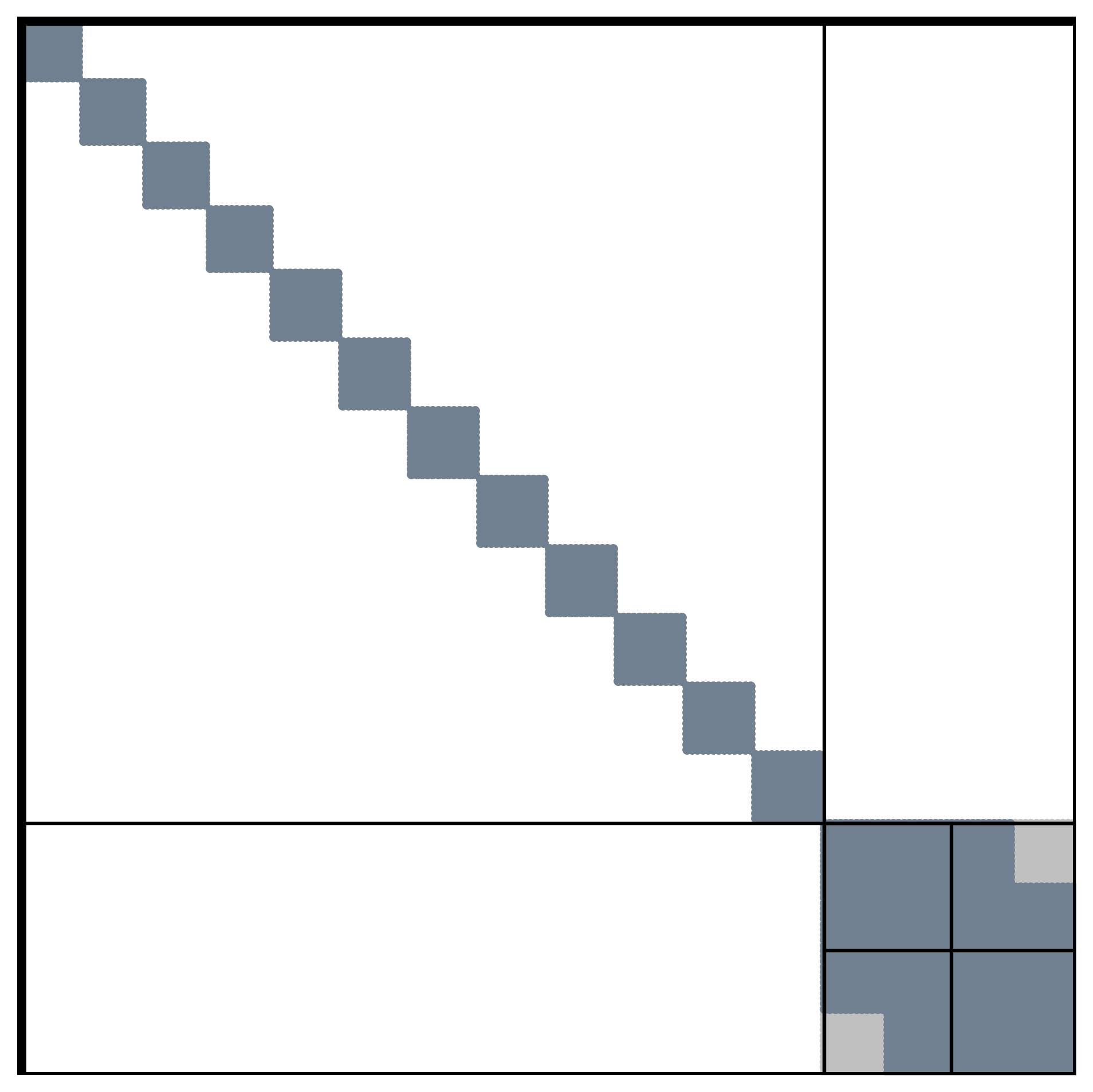}};
  \draw[->, thick, font=\tikzlabelsize] (level2mat.east) -- node[below]{\shortstack{$\msf{coarsen\ \&}$\\$\msf{diagonalize}$}} (level1mat.west);
  
\node[below=-0.5cm of level2mat] (level2geom) {\includegraphics[width=60mm]{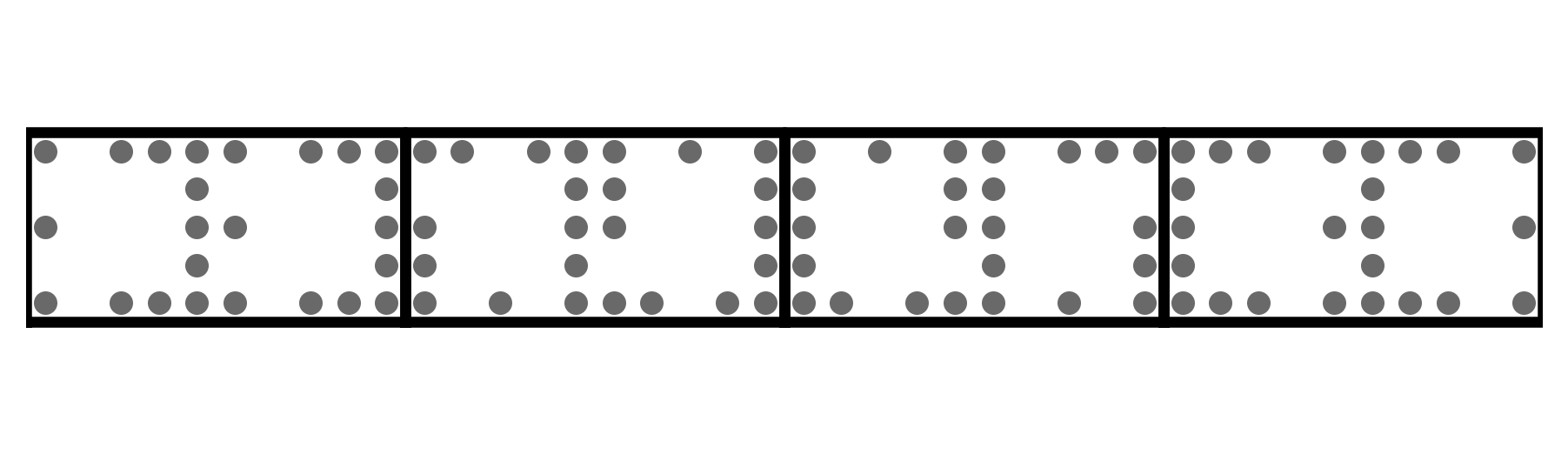}};
  \node[below=-0.5cm of level1mat] (level1geom) {\includegraphics[width=60mm]{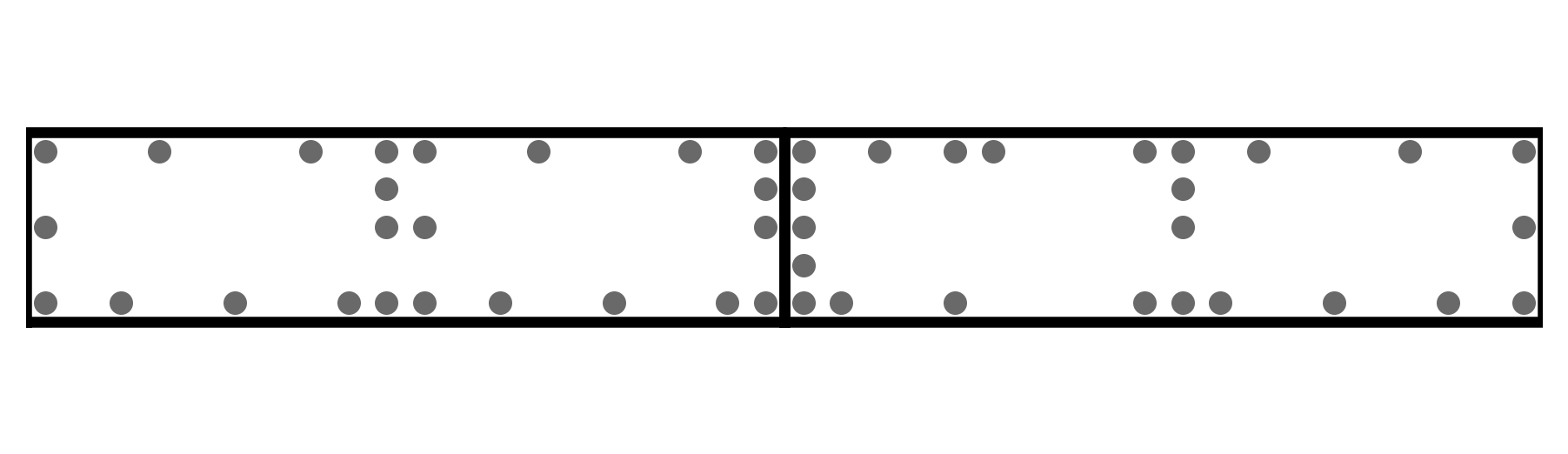}};
  \draw[->, thick,font=\tikzlabelsize] (level2geom.east) -- (level1geom.west);

  \annotategeomcoarsen{level2geom}{4}{8}
  \annotategeomcoarsen{level1geom}{2}{12}

  \annotatematrixcoarsen{level2mat}{4}{8}{0.5}{0}
  \annotatematrixcoarsen{level1mat}{2}{12}{0.75}{0.03}

\draw[overlay,->, thick] (level3geom.south east)
    .. controls +(1.5cm,-1.5cm) and +(-1.5cm,1.5cm) ..
    node[midway, sloped, below,font=\tikzlabelsize] {$\msf{coarsen}$}
    (level2mat.north west);
\end{tikzpicture}
 }
\caption{For a pseudo-1D geometry, the active points as well as the corresponding matrix are shown
at various stages of the computation. First, the redundant points of the boxes on the finest level are diagonalized.
To continue the computation, the remaining active points are regrouped according to the next coarse level of the tree,
which introduces low-rank subblocks that can be further diagonalized.}
\label{fig:SRS1d}
\end{figure}

The same procedure is applied to each subsequent box on level $L$, with each step operating on the current modified system while preserving the indices already diagonalized.
Let $n_L$ denote the number of boxes on the finest level $L$ (all of which are leaves). 
After diagonalizing all boxes on that level, the modified matrix has the form
\[
\widetilde{\mtx A}(\mtx A; B_1, \dots, B_{n_L})
= \mtx V_{n_L}^{-1} \cdots \mtx V_{1}^{-1}\,
\mtx A\, 
\mtx W_{1}^{-1} \cdots \mtx W_{n_L}^{-1},
\]
where $\msf R_1, \dots, \msf R_{n_L}$ are diagonalized.
The top row of Figure~\ref{fig:SRS1d} illustrates this process for a pseudo-1D
domain with eight leaf boxes.

To continue strong recursive skeletonization, we regroup the remaining skeleton indices according to a coarser
level of the tree. 
The interactions between boxes and their far-field at this coarser scale are
again numerically low rank, and can be sparsified and diagonalized using the same procedure
(cf. Remark~\ref{remark:nested_bases}).

We introduce notation for the remaining \textit{active} set of points, which is updated after diagonalizing indices $\msf R_1, \dots,  \msf R_j$,
\begin{equation}
\msf{active}\ =\ [1,\dots, N] \setminus \left (\ \msf R_1\  \cup\ \dots\ \cup\ \msf R_{j}\  \right),
\label{eq:active}\end{equation}
as well as corresponding notation for the remaining \textit{active} box and neighbor points for a box $B$
\begin{equation}
\msf B_{\msf {active}}\ =\ \msf B\ \cap\ \msf {active}, \qquad \msf N_{\msf {active}}\ =\ \msf N\ \cap\ \msf {active}.
\label{eq:BNactive}
\end{equation}
Only these active indices are used when computing the diagonalization matrices
$\mtx V^{-1}$ and $\mtx W^{-1}$ for subsequent steps.

\begin{figure}[!htb]
\resizebox{\textwidth}{!}{\begin{tikzpicture}[node distance=0.55cm and 0.55cm, auto]
\node (m1) {\includegraphics[width=45mm]{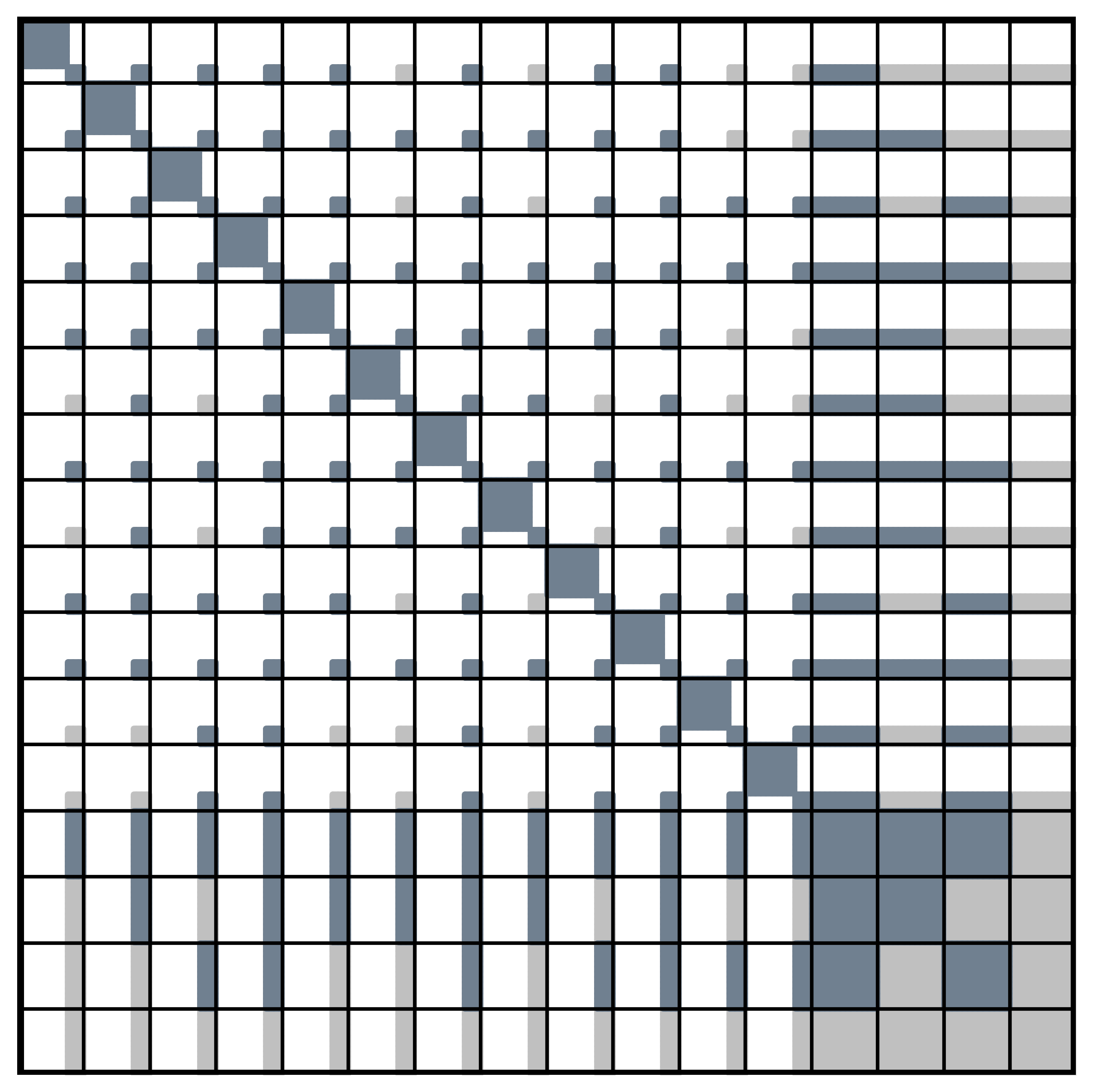}};
  \node (m2) [right=1.5cm of m1] {\includegraphics[width=45mm]{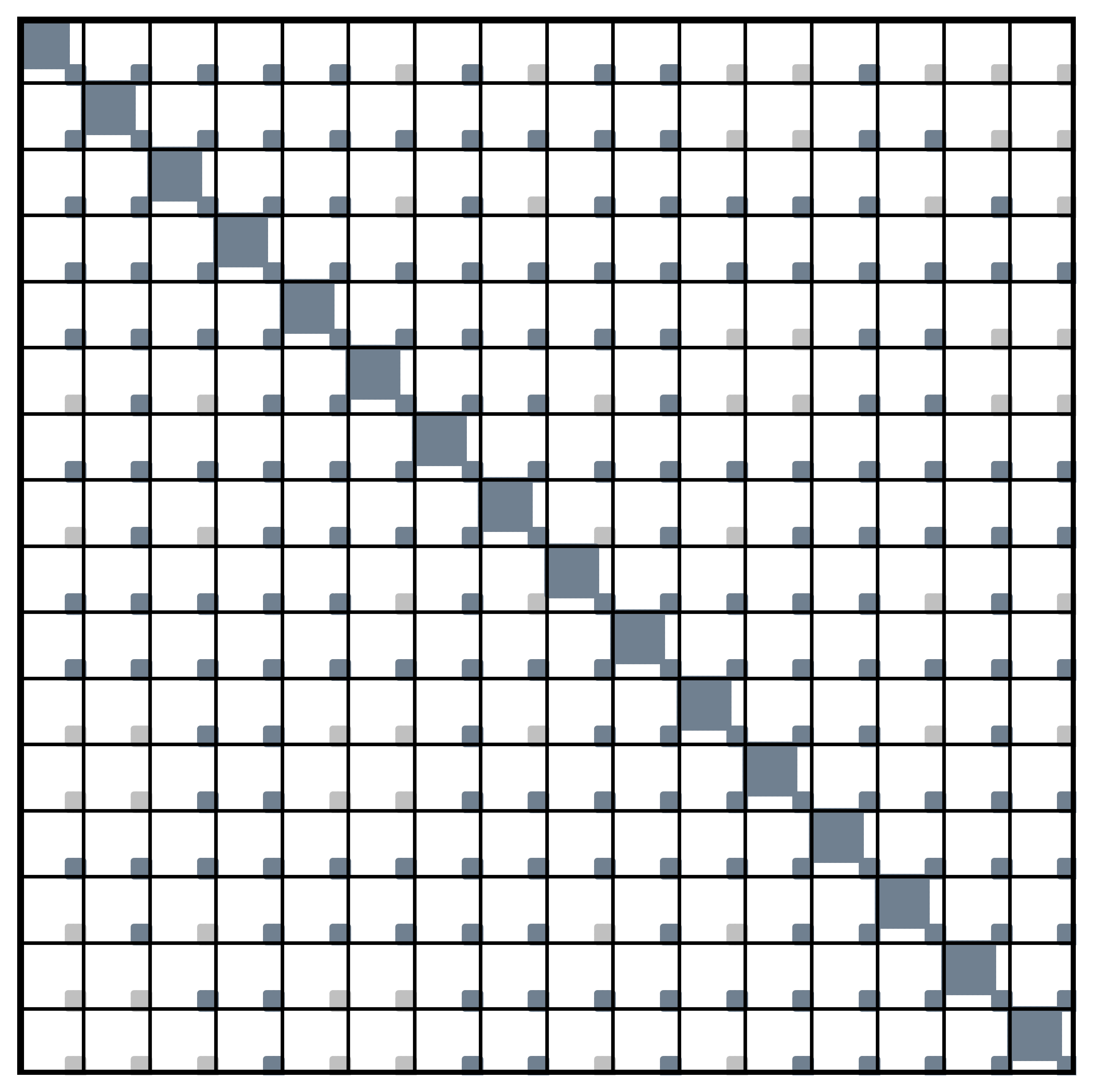}};
  \node (m3) [right=1.5cm of m2] {\includegraphics[width=45mm]{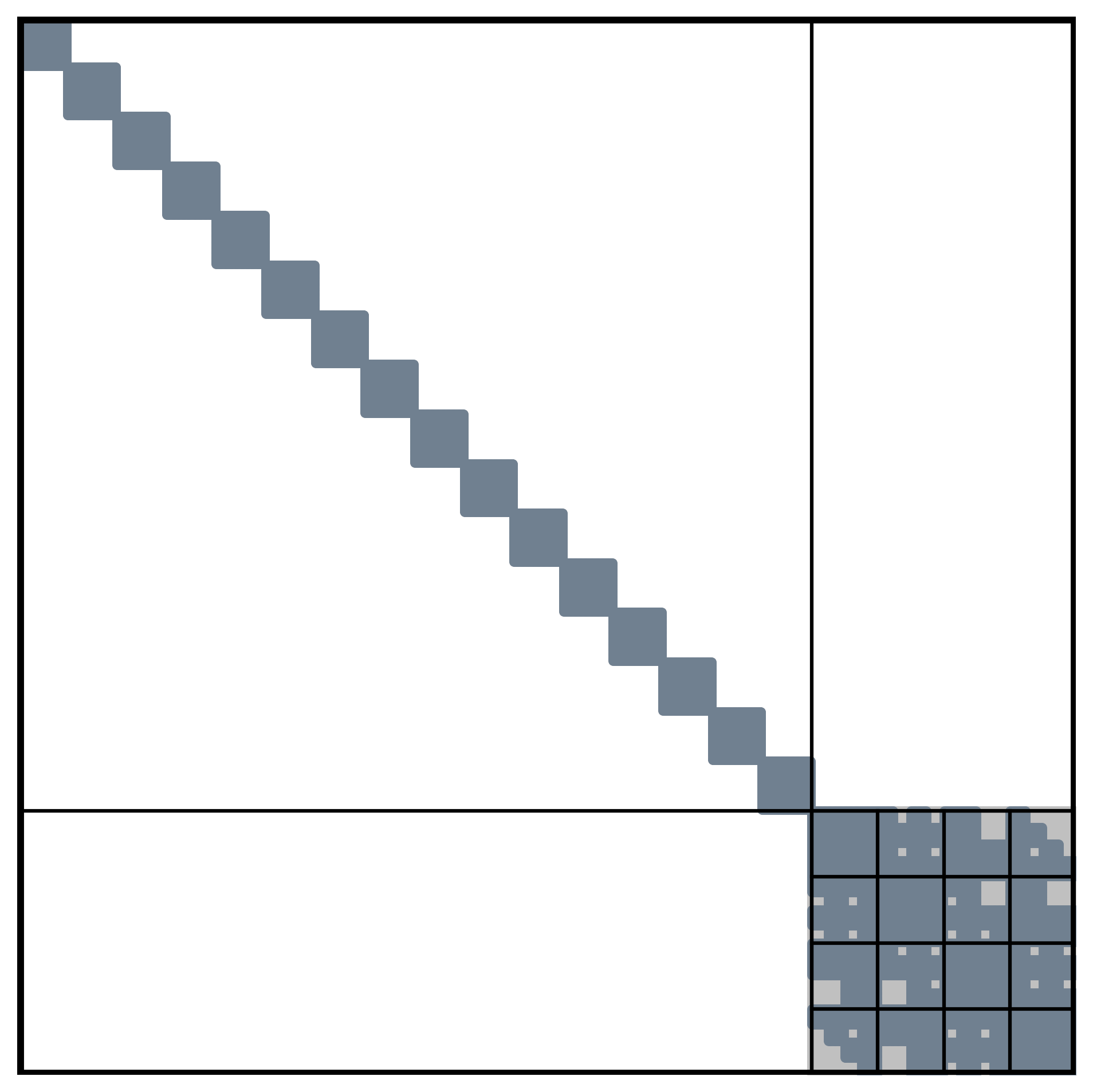}};

\node (d1) [below=0.00cm of m1,xshift=0.0cm] {\includegraphics[width=35mm]{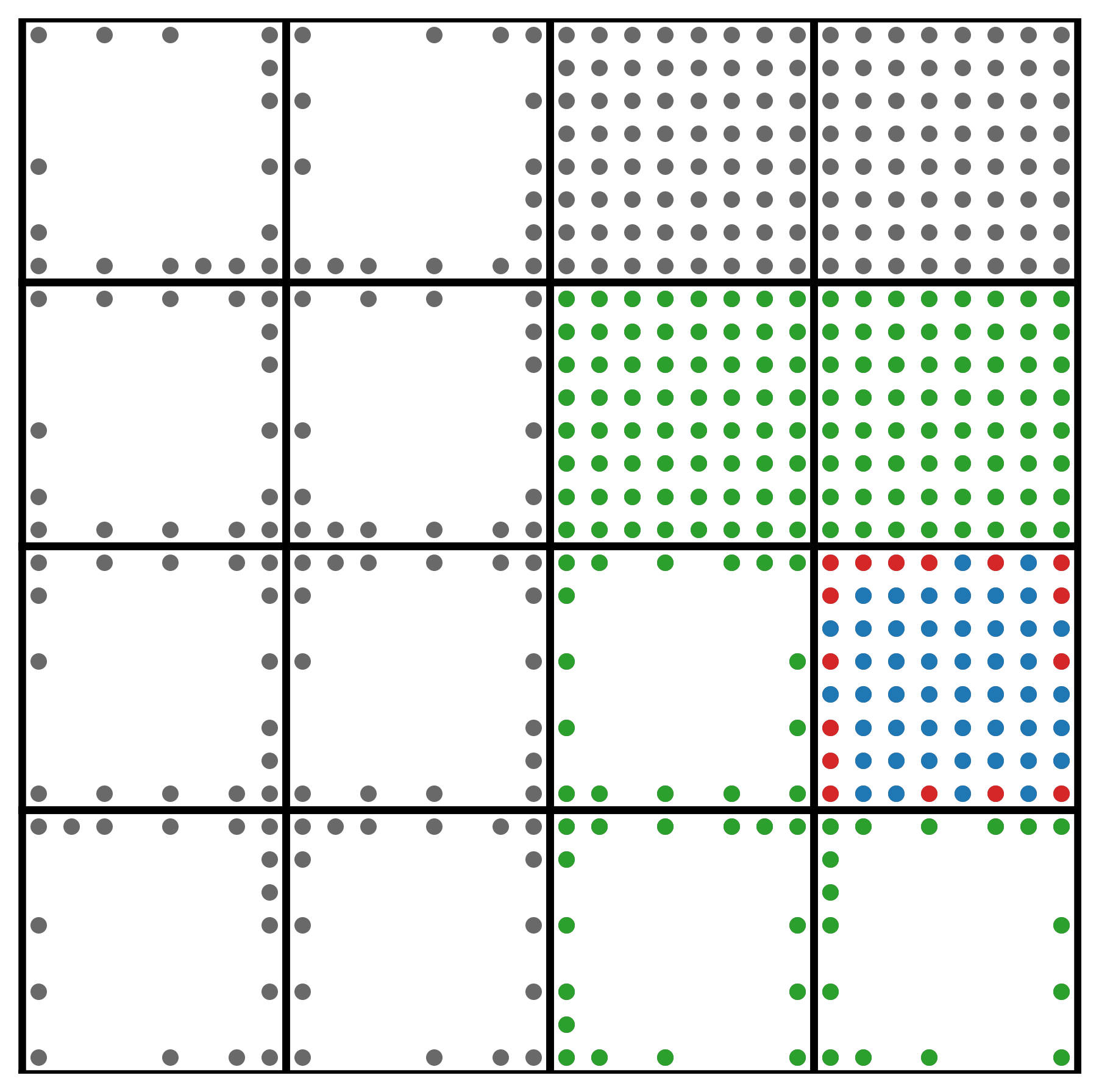}};
  \node (d2) [below=0.00cm of m2,xshift=0.0cm] {\includegraphics[width=35mm]{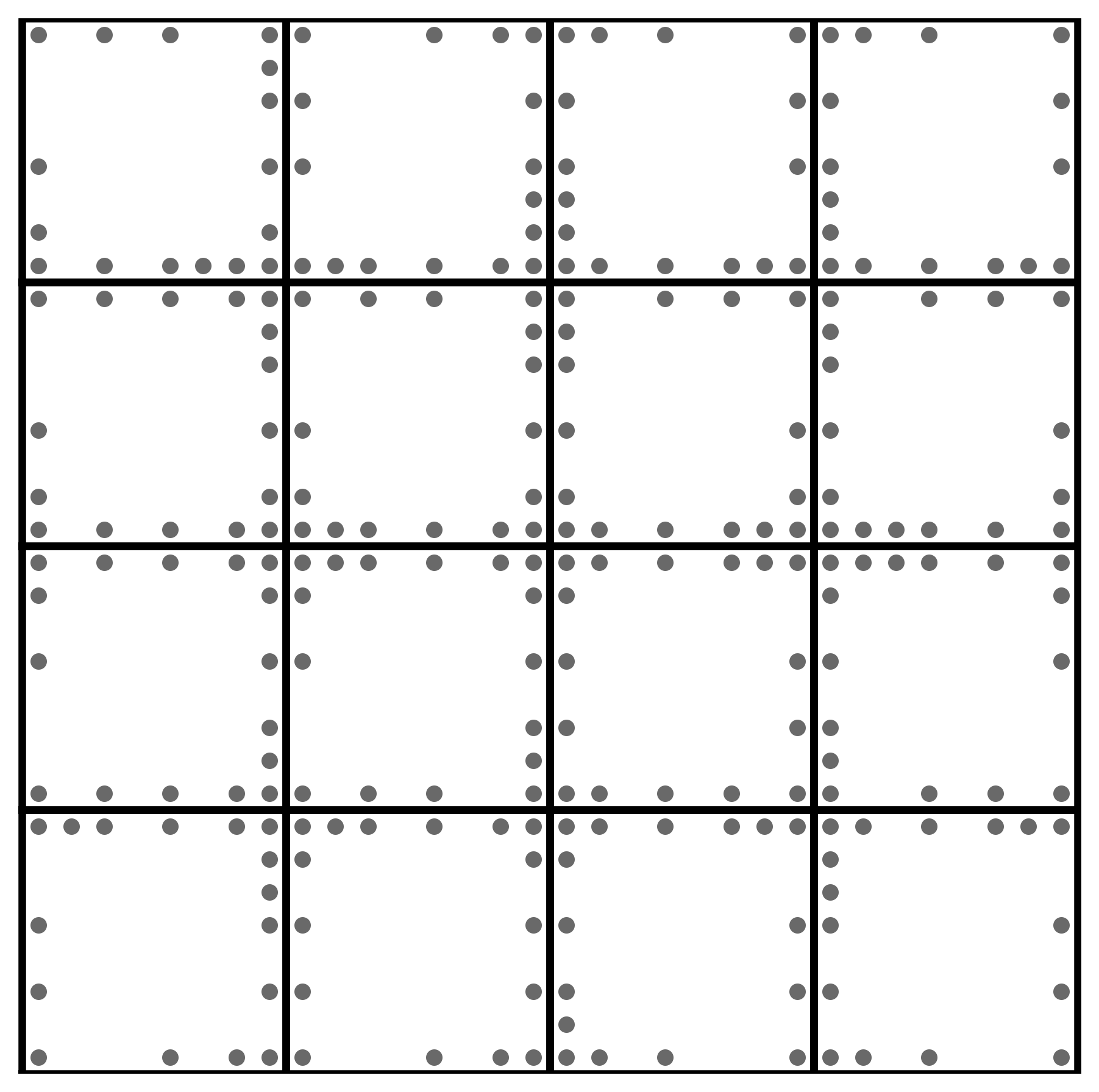}};
\node (d3) [below=0.00cm of m3,xshift=0.0cm] {\includegraphics[width=35mm]{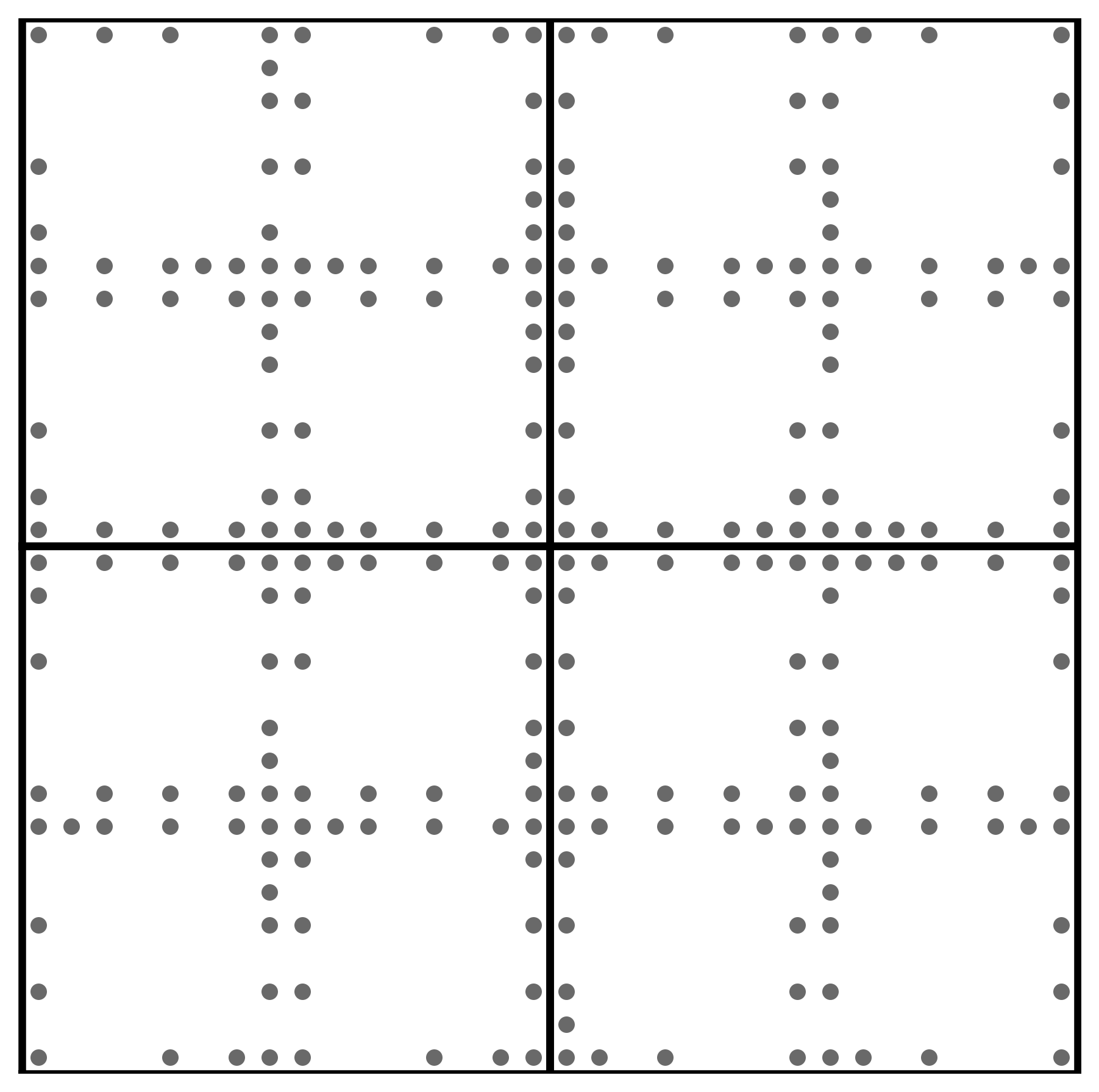}};
  
\draw[->, thick] (m1) -- node[below, midway,font=\tikzlabelsize] {$ \msf{diagonalize}$} (m2);
  \draw[->, thick] (m2) -- node[below, midway,font=\tikzlabelsize] {$ \msf{coarsen}$} (m3);
  
\draw[->, thick] (d1) -- node[above, midway] {} (d2);
  \draw[->, thick] (d2) -- node[above, midway] {} (d3);
\end{tikzpicture}
 }
\caption{The analog of Figure \ref{fig:SRS1d} when the computational domain is a unit box. Observe that many more blocks get updated for a true two-dimensional domain. For every diagonalization step, up to 35 pairwise interactions are introduced between non-neighboring boxes, which must be stored and later recompressed.}
\label{fig:SRS2d}
\end{figure}

When matrix entries are directly available, the modified entries of the diagonalized matrix are stored explicitly, and the unmodified entries can be accessed directly as needed. 
As we noted in Remark \ref{remark:additive_lowrank}, strong recursive skeletonization introduces additive
terms between boxes that may be well-separated. These interactions must be recompressed at later stages of the algorithm.
In computing the ID of the subblock $\widetilde{\mtx A}_{\msf B_{\rm active}\msf F}$,
the modified entries are handled algebraically, whereas the unmodified entries can be compressed using analytic methods
such as proxy surfaces or adaptive cross approximation. The modified entries in the far-field remain spatially localized, and efficient methods for their compression are detailed in \cite{sushnikova2023fmm}.
In 2D and 3D, significantly more modified interactions arise---see Figure \ref{fig:SRS2d}---and these must be stored and recompressed.

\begin{remark}[Nested bases]
A defining property of the $\mathcal{H}^2$ representation is that the bases are
\emph{nested} across levels of the hierarchical tree. In the present factorization,
this is achieved by constructing skeleton sets hierarchically: for each box $B$, a single
subset of indices—the \emph{skeleton set}—represents both the column space of
$\mtx A_{\msf B \msf F}$ and the row space of $\mtx A_{\msf F \msf B}$. The skeletons of a
parent box are selected from the union of its children’s skeletons so that they span the
corresponding interactions $\mtx A_{\msf B \msf F}$ and $\mtx A_{\msf F \msf B}$. 
This hierarchical construction is consistent with the geometric fact that any box
well-separated from the parent is also well-separated from each of its children.
\label{remark:nested_bases}
\end{remark}

To formalize the description of recursive skeletonization with strong admissibility, consider that the boxes are diagonalized in an upward traversal through the tree in order 
$1, \dots, M$, terminating at box $B_M$. 
Let the index vector $\msf B_t$ denote the remaining active points in the domain at the time the algorithm terminates,
and let the permutation vector $\msf P_t$ record the order in which points are eliminated.
We define the remaining dense active submatrix by extracting subindices
\begin{equation*}
\widetilde{\mtx A}_{\msf B_t \msf B_t}, \qquad \text{where}\ \quad \widetilde {\mtx A}\ =\ \widetilde {\mtx A} ( \mtx A; B_1, \dots B_M).
\end{equation*}
The full diagonalized matrix takes the form
\begin{align}
\begin{split}
\mtx P_t\ \mtx D\ \mtx P_t^T\ &=\
\mtx P_t\ 
\begin{pmatrix}
\mtx X_{\msf R_1 \msf R_1}\\
& \ddots\\
&& \mtx X_{\msf R_M \msf R_M}\\
&&&\widetilde{\mtx A}_{\msf B_t \msf B_t}
\end{pmatrix}\ \mtx P_t^T\\ &=\ \widetilde {\mtx A}( \mtx A; B_1, \dots B_M)\ \approx\  \mtx V_M^{-1}\  \dots \mtx V_1^{-1}\ \mtx A\ \mtx W_1^{-1}\ \dots\ \mtx W_M^{-1}
\end{split}
\label{eq:diag_sys}
\end{align}
The diagonalized system (\ref{eq:diag_sys})
leads to a sparse factorization of $\mtx A$,
and, by inversion, an expression for $\mtx A^{-1}$, since each of the diagonalization matrices is a composition of sparse,
invertible matrices
\begin{alignat}{3}
\mtx{A} &\approx\ & \mtx{V}_1 \cdots \mtx{V}_M\,& \mtx P_t\ \mtx{D}\ \mtx P_t^T \,& \mtx{W}_M \cdots \mtx{W}_1, \label{eq:Afact}\\
\mtx{A}^{-1} &\approx\ & \underbrace{\mtx{W}_1^{-1} \cdots \mtx{W}_M^{-1}}_{\text{downward pass}}\,& \mtx P_t\ \mtx{D}^{-1}\ \mtx P_t^T\ & \underbrace{\mtx{V}_M^{-1} \cdots \mtx{V}_1^{-1}}_{\text{upward pass}}.
\label{eq:Ainvfact}
\end{alignat}
For symmetric positive definite matrices, the diagonalization matrices $\mtx V, \mtx W$
are symmetric as well, and the factorization (\ref{eq:Afact}) can be used to compute the matrix square root $\sqrt{\mtx A}$.
\section{Randomized Compression and LU Factorization using Sketching}
\label{sec:randomized_compression_and_lu}

In this section, we describe randomized methods for recovering an invertible factorization of $\mtx A$, as defined in Section~\ref{sec:skel_algos}, in settings where direct access to matrix entries is prohibitively expensive. Instead, we assume that $\mtx A$ is accessible only through matrix-vector products with $\mtx A$ and its adjoint $\mtx A^*$. The factorization is reconstructed from random sketches of the form
\begin{equation} \label{eq:rand_sample_setting}
\underset{N \times s}{\mtx Y} = \underset{N \times N}{\mtx A} \underset{N \times s}{\mtx \Omega}, \qquad 
\underset{N \times s}{\mtx Z} = \underset{N \times N}{\mtx A^*} \underset{N \times s}{\mtx \Psi}, \qquad \mtx \Omega, \mtx \Psi \sim \mathcal{N}(\mtx 0, \mtx I),
\end{equation}
by post-processing the set $\{\mtx Y, \mtx \Omega, \mtx Z, \mtx \Psi\}$.
We defer discussion of the required sample size $s$ to later sections, but for now assume $s \ll N$.
Recall from Section~\ref{sec:strong_skel} that diagonalizing a set of indices $\msf R \subseteq \msf B$ proceeds by first compressing interactions between a target box and its far field, followed by extracting the remaining near-field (inadmissible) subblocks between $\msf R$ and neighboring indices $\msf N$.

A natural strategy for performing these operations within the randomized sketching framework is to design \textit{structured} test matrices $\mtx \Omega, \mtx \Psi$ containing zero or identity subblocks. For example, zeroing out the near-field subblocks isolates far-field interactions $\mtx A_{\msf B \msf F}$, while placing identity blocks on the near field allows direct extraction of sparsified near-field interactions.
However, constructing such structured test matrices poses significant challenges. First, they typically require substantially larger sample sizes to ensure accurate recovery. Second, and more fundamentally, the skeleton indices are computed dynamically during factorization (due to recompressions discussed in Section~\ref{sec:skel_algos}), making it impractical to design structured test matrices in advance.

Instead, we adopt a more flexible approach: we begin with dense Gaussian test matrices $\mtx \Omega, \mtx \Psi$ and impose the necessary structure \textit{as needed} during the algorithm via linear transformations. This strategy preserves the generality of randomized sketching while enabling both far-field compression and near-field extraction. Our approach builds on techniques introduced in~\cite{2022_levitt_dissertation}, trading modest post-processing overhead for a substantial reduction in the number of required samples.

\vspace{0.5em}

\begin{figure}[!htb]
\centering
\begin{subfigure}{0.45\textwidth}
\centering
\resizebox{!}{5.5cm}{\begin{tikzpicture}[>=latex, every node/.style={font=\figtitlesize}]

\node (Y) {\includegraphics[height=3.5cm]{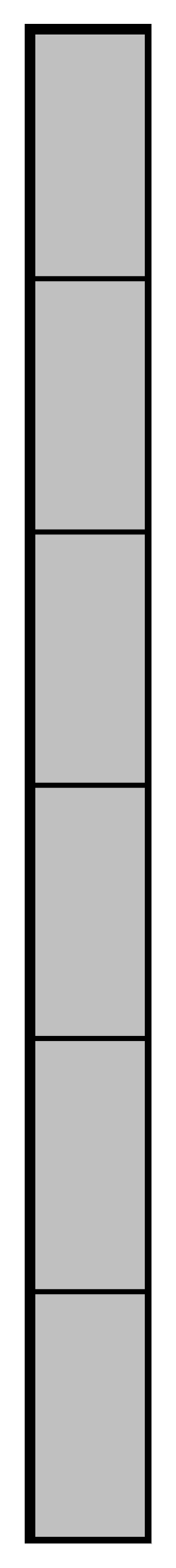}};
\node (Ylabel) [above=2mm of Y] {$\mtx Y'$};

\node (eqtop) [right=2mm of Ylabel] {$=$};

\node (eqmid) at (eqtop |- Y) {$=$};

\node (A) [right= of Y] {\includegraphics[height=3.5cm]{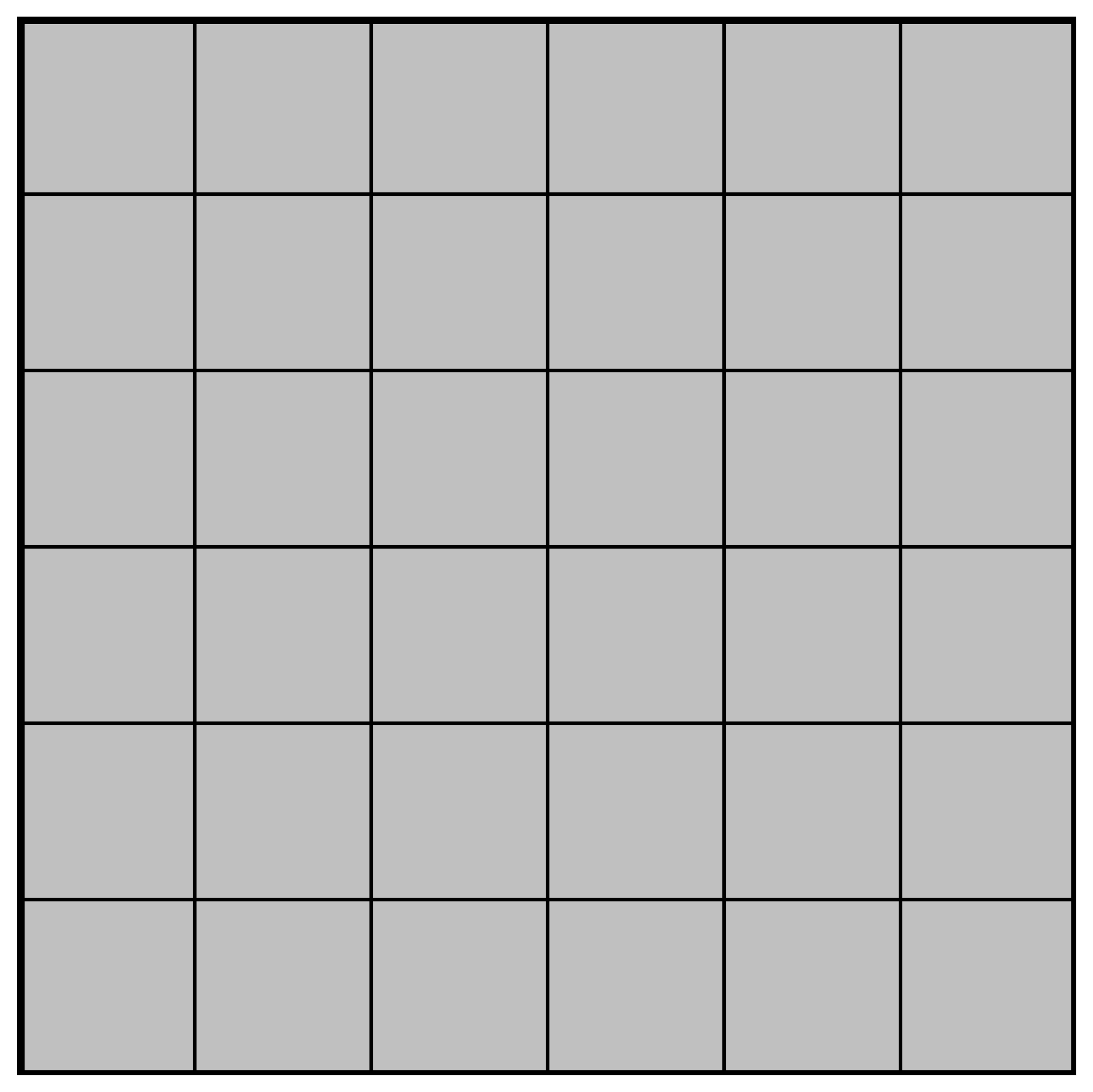}};
\node [overlay,above=2mm of A] {$\mtx A$};

\node (Omega) [right=0.05cm of A] {\includegraphics[height=3.5cm]{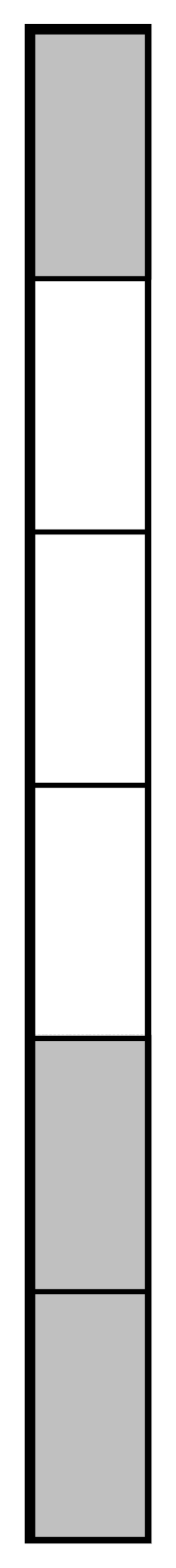}};
\node [above=2mm of Omega] {$\mtx \Omega'$};

\node (geom) [below=-0.4cm of A] {\includegraphics[width=3.5cm]{fig_SRS_demosketchtmp_step9_skel_geom.png}};

\end{tikzpicture}
 }
\caption{Setting for block nullification.}
\label{fig:setting_blocknull}
\end{subfigure}\hfill
\begin{subfigure}{0.55\textwidth}
\centering
\resizebox{!}{5.5cm}{\begin{tikzpicture}[>=latex, every node/.style={font=\figtitlesize}]

\node (Y) {\includegraphics[height=3.5cm]{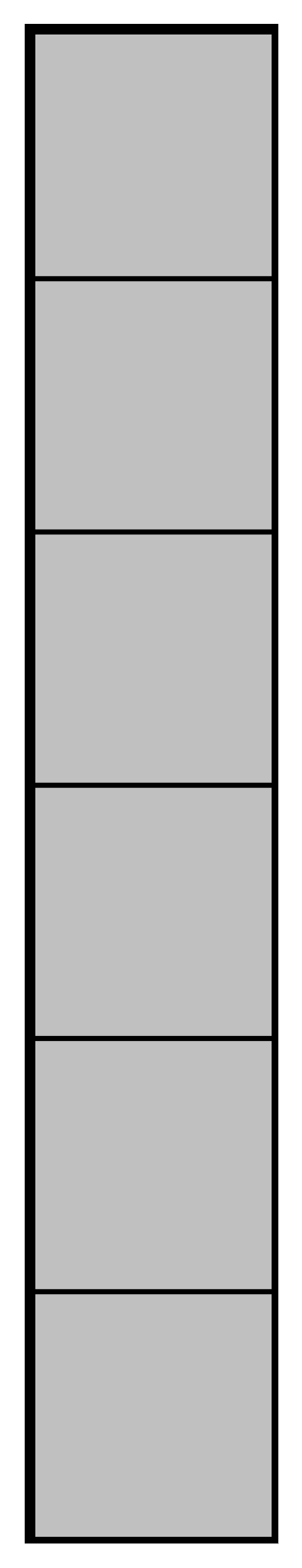}};
\node (Ylabel) [above=2mm of Y] {$\mtx Y'$};

\node (eqtop) [right=2mm of Ylabel] {$=$};
\node (eqmid) at (eqtop |- Y) {$=$};   

\node (A) [right= of Y] {\includegraphics[height=3.5cm]{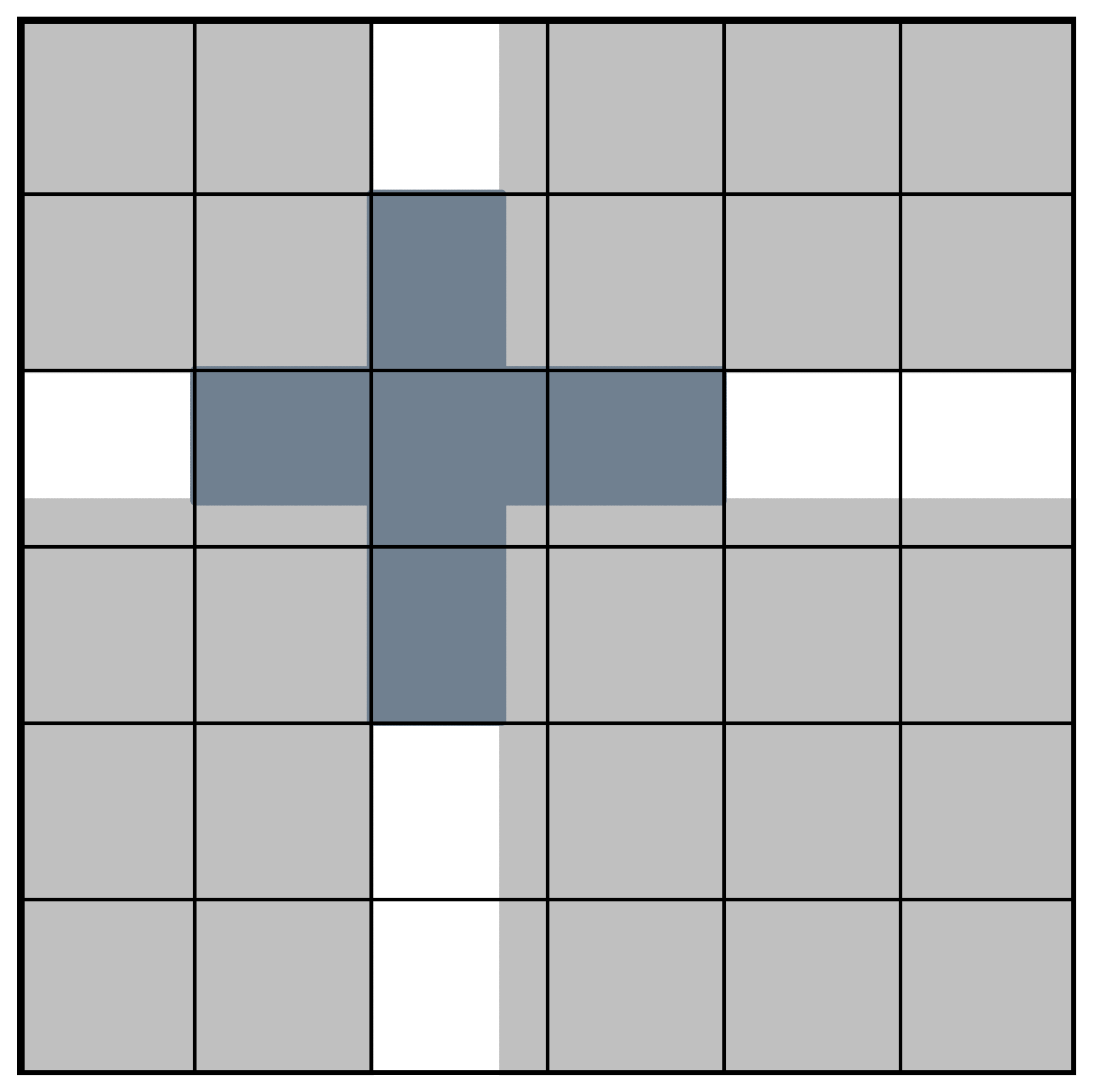}};
\node [overlay,above=2mm of A] {$\hat{\mtx A}$};

\node (Omega) [right=0.05cm of A, inner sep=0]
      {\includegraphics[height=3.5cm]{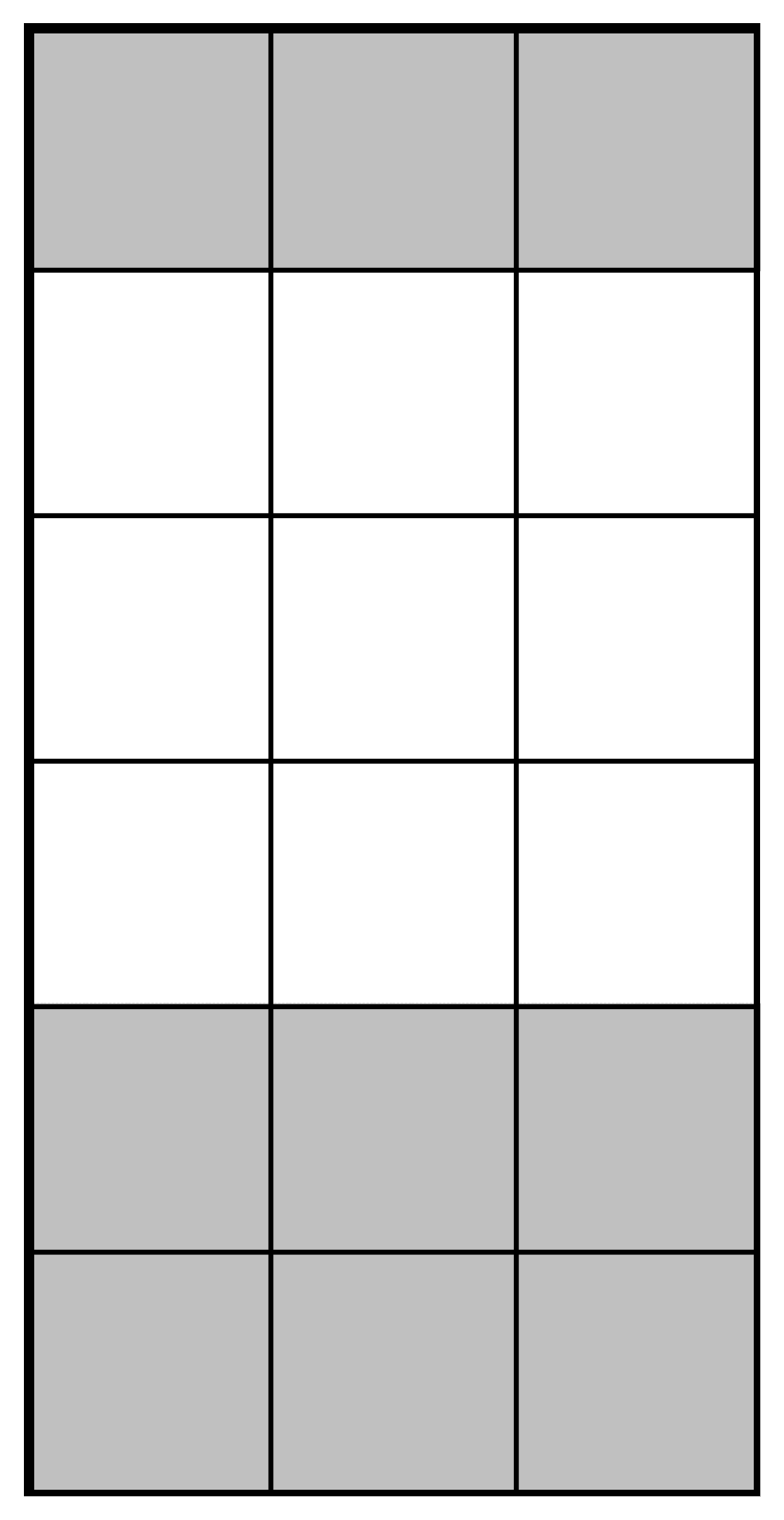}};
\node [above=2mm of Omega] {$\mtx \Omega'$};

\def\ncols{3}\def\nrows{6}
\begin{scope}[shift={(Omega.south west)},
              x={(Omega.south east)}, y={(Omega.north west)}]
  \tikzset{cell/.style={anchor=center}}
  \foreach \r/\c in {2/1, 3/2, 4/3}{
    \node[cell] at ({(\c-.5)/\ncols}, {1-(\r-.5)/\nrows}) {$\mtx I$};
  }
\end{scope}

\node (geom) [below=-0.4cm of A]
      {\includegraphics[width=3.5cm]{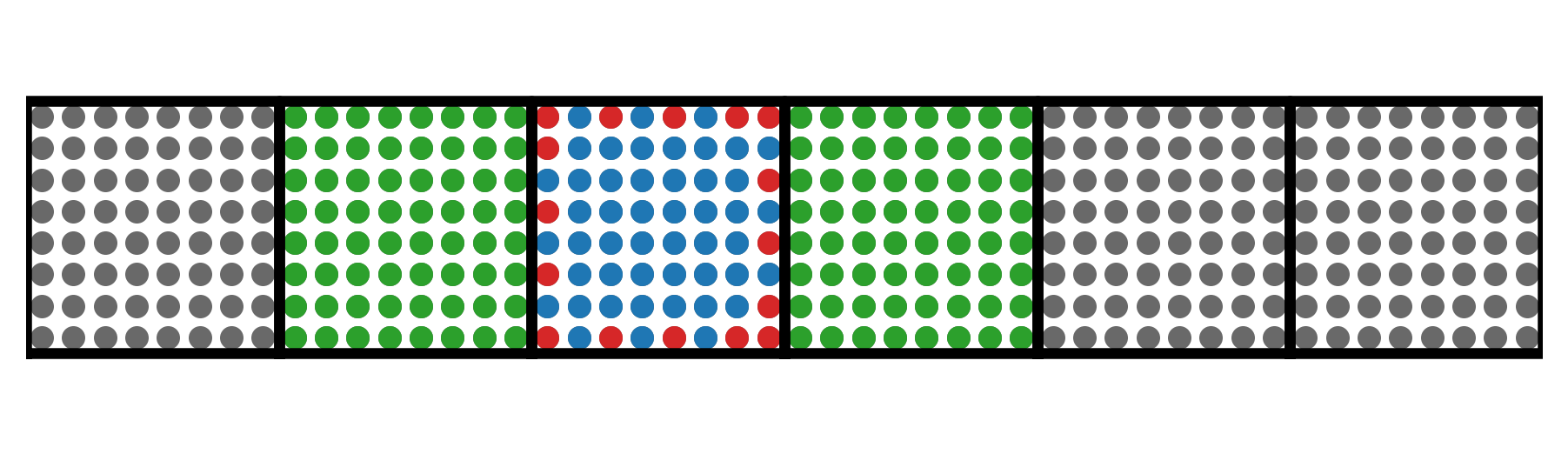}};

\end{tikzpicture}
 }
\caption{Setting for block extraction.}
\label{fig:setting_blockextract}
\end{subfigure}
\caption{Two types of test matrices are needed for \texttt{RSRS}.
(\subref{fig:setting_blocknull}) To sketch far-field interactions between a box and its distant neighbors.
(\subref{fig:setting_blockextract}) To extract near-field interactions between a subset of box points and neighboring points.
Rather than designing structured test matrices directly, we apply \textit{block nullification} and \textit{block extraction}, which use linear transformations of dense Gaussian matrices to introduce the desired structure.}
\label{fig:separable}
\end{figure}

\noindent \textit{Block Nullification:} Applies a linear transformation to a Gaussian test matrix $\mtx \Omega$ to produce a modified matrix $\mtx \Omega'$, where the contribution of the near field has been ``nullified.'' This enables efficient sampling of far-field interactions $\mtx A_{\msf B \msf F}$. See Figure \ref{fig:setting_blocknull} for an illustration of $\mtx \Omega'$.

\vspace{0.5em}

\noindent \textit{Block Extraction:} Applies a linear transformation to a Gaussian test matrix $\mtx \Omega$ to produce a modified matrix $\mtx \Omega'$ that ``extracts'' specific subblocks of a sparse matrix. This allows near-field subblocks to be extracted without prior knowledge of their locations. See Figure \ref{fig:setting_blockextract} for an illustration of $\mtx \Omega'$.

\subsection{Block Nullification}\label{sec:block_null}
Suppose we would like to compute the interpolative decomposition of the interactions between box indices $\msf B$ and far field indices $\msf F$ so that
(\ref{eq:skel_lowrank_simult}) holds. To accomplish this using
the randomized sketching (as discussed in Section \ref{sec:ID}), we need to generate
and postprocess the sketches
\begin{equation}
\mtx Y_{\msf B}'\  = \mtx A_{\msf B \msf F}\  \mtx \Omega'_{\msf F}, \qquad \mtx Z_{\msf B}' = \mtx A_{\msf F \msf B}^*\ \mtx \Psi'_{\msf F},
\label{eq:sketch_needed_blocknull_setup}
\end{equation}
where $\mtx \Omega'_{\msf F}$ and $\mtx \Psi'_{\msf F}$ are Gaussian random matrices.

For concreteness, suppose that each block in the tessellation has size at most $m$, and the test and sketch matrices are tessellated according to the decomposition in Figure \ref{fig:setting_blocknull}.
The matrix $\mtx A$ has full rank interactions between the target box and near-field boxes, complicating the straightforward use of randomized sketching. 
Ideally, the test matrices should reflect the sparsity pattern of the low-rank blocks we aim to sample.
Consider a \textit{structured} test matrix that, under an appropriate permutation $[\msf B, \msf N, \msf F]$, is designed to isolate the far-field interactions $\mtx A_{\msf B \msf F}$:
\begin{equation}
\underset{N \times (k+p)}{\mtx Y'}\ =\ \underset{N \times N}{\mtx A}\ \underset{N \times (k+p)}{\mtx \Omega'}, \qquad \text{where}\
{\mtx \Omega'}\ =\ \begin{pmatrix} \mtx 0_{\msf B}\\ \mtx 0_{\msf N}\\ \mtx \Omega'_{\msf F}\end{pmatrix}
\label{eq:Yprime_null}
\end{equation}
Then, extracting the subblock $\mtx Y'_{\msf B}$ yields the required sketch in (\ref{eq:sketch_needed_blocknull_setup}).

Block nullification achieves the same objective without explicitly constructing a sparse test matrix. 
Instead, we apply a linear transformation to a dense Gaussian matrix $\mtx \Omega$ to produce a modified test matrix $\mtx \Omega'$ with desired zero subblocks.
Specifically, given $\mtx \Omega \in \mathbb{R}^{N \times s}$ drawn from a Gaussian distribution, we compute a nullspace basis that annihilates the contributions of $\msf B$ and $\msf N$, resulting in
the desired test matrix
\begin{equation}
\underset{N \times (k+p)}{\mtx \Omega'} = \underset{N \times s\vphantom{()}}{\mtx \Omega}\ \underset{s \times (k+p)}{\mtx N'}, \qquad \text{where} \quad \mtx N' = \texttt{null} \begin{pmatrix}\mtx \Omega_{\msf B}\\ \mtx \Omega_{\msf N}\end{pmatrix},
\label{eq:Omegaprime_null}
\end{equation}
for a sufficient number of samples $s \ge |\msf B| + |\msf N| + k + p$.

Since we already have a sketch of $\mtx A$ in (\ref{eq:rand_sample_setting}), we do not need to explicitly form $\mtx \Omega'$. Instead, we can apply the nullspace basis $\mtx N'$ directly to the extracted block:
\begin{equation}
\mtx Y_{\msf B}' =
\underbrace{\mtx A_{\msf B,:} \mtx \Omega'}_{\text{extract } \msf B \text{ from } (\ref{eq:Yprime_null})}
=
\mtx A_{\msf B,:}
\underbrace{\mtx \Omega \mtx N'}_{\text{using } (\ref{eq:Omegaprime_null})}
=
\underbrace{\mtx A_{\msf B,:} \mtx \Omega}_{\text{extract } \msf B \text{ from } (\ref{eq:rand_sample_setting})}
\mtx N'
= \mtx Y_{\msf B} \mtx N'.
\label{eq:nullformula}
\end{equation}
The cost of post-processing sketches using block nullification is $\mathcal{O}(s^3)$, where $s$ depends on the number of points in target box and the near field. Because the nullspace basis $\mtx N'$ has orthonormal columns, block nullification preserves the spectral properties of the original test matrix.

\subsection{Block Extraction}\label{sec:block_extract}

In addition to sketching low-rank factors, we also need to extract the interactions between a target box and its near-field neighbors. In the factorization described in Section~\ref{sec:skel_algos}, the exact indices that need to be extracted are not known a priori---they are determined dynamically over the course of the algorithm depending on which indices are selected as skeletons.
Designing structured test matrices to extract these subblocks directly would both increase the sample complexity and require additional sketching operations as the algorithm proceeds. Instead, block extraction is similar in spirit to block nullification and applies a linear transformation to a dense Gaussian test matrix $\mtx \Omega$ to introduce structure into the test matrix \textit{as needed} during post-processing.

Suppose that the matrix $\mtx A$ is partitioned according to $[\msf B, \msf N, \msf F]$, and that for a subset $\msf R \subseteq \msf B$, the interactions with the far field satisfy ${\mtx A_{\msf R \msf F}} \approx \mtx 0$. Our goal is to recover the remaining dense subblocks $\mtx A_{\msf R \msf B}$ and $\mtx A_{\msf R \msf N}$.
A natural approach is to construct a \textit{structured} test matrix that directly enables extracting these subblocks. Define $\mtx \Omega'$ and its corresponding sketch $\mtx Y'$ as:
\begin{equation}
\underset{N \times l}{\mtx Y'} = \underset{N \times N}{\mtx A} \underset{N \times l}{\mtx \Omega'}, \qquad \text{where} \quad
\mtx \Omega' =
\begin{pmatrix}
\mtx I_{\msf B} & \\[0.5em]
& \mtx I_{\msf N}\\[0.5em]
\multicolumn{2}{c}{\mtx \Omega'_{\msf F}}
\end{pmatrix}, \quad l = |\msf B| + |\msf N|.
\label{eq:Yprime_pinv}
\end{equation}
Extracting $\mtx Y'_{\msf R}$ from this sketch directly yields the desired subblock.

Rather than constructing $\mtx \Omega'$ explicitly, we can compute it implicitly by applying a pseudoinverse transformation to the subblocks of a dense Gaussian matrix $\mtx \Omega \in \mathbb{R}^{N \times s}$. Specifically:
\begin{equation}
\underset{N \times l}{\mtx \Omega'} = \underset{N \times s}{\mtx \Omega}\ \underset{s \times l}{\mtx P'}, \qquad \text{where} \quad
\mtx P' = \begin{pmatrix} \mtx \Omega_{\msf B}\\ \mtx \Omega_{\msf N} \end{pmatrix}^{\dagger},
\label{eq:Omegaprime_pinv}
\end{equation}
for $s > |\msf B| + |\msf N|$.
Since we have already drawn a sketch of $\mtx A$ as in (\ref{eq:rand_sample_setting}), we can apply $\mtx P'$ directly to a small extracted subblock without forming $\mtx \Omega'$ explicitly by
following a similar procedure to (\ref{eq:nullformula}):
\begin{equation}
\begin{pmatrix} \mtx A_{\msf R \msf B} & \mtx A_{\msf R \msf N} \end{pmatrix} =
\underbrace{\mtx A_{\msf R,:} \mtx \Omega'}_{\text{extract from } (\ref{eq:Yprime_pinv})} =
\mtx A_{\msf R,:} \underbrace{\mtx \Omega \mtx P'}_{\text{using } (\ref{eq:Omegaprime_pinv})} =
\underbrace{\mtx A_{\msf R,:} \mtx \Omega}_{\text{extract from } (\ref{eq:rand_sample_setting})} \mtx P' =
\mtx Y_{\msf R} \mtx P'.
\label{eq:pinvformula}
\end{equation}
The cost of post-processing the sketches using block extraction is $\mathcal{O}(s^3)$. Although the resulting entries $\mtx \Omega'_{\msf F}$ are no longer Gaussian, this does not substantially affect the accuracy of extracting sparse blocks, since these entries are multiplied by components of $\mtx A$ that are (by construction) numerically small.

\subsection{Randomized Strong Skeletonization}
\label{sub:randomized_strong_skeletonization}

In Section~\ref{sec:strong_skel}, we described how to construct diagonalization matrices that decouple a subset of box indices $\msf R \subseteq \msf B$ from their far field $\msf F$, in settings where matrix entries are easily accessible. We now extend this approach to the setting where the matrix is only accessible through its action on vectors, using the techniques of block nullification and block extraction introduced in Sections~\ref{sec:block_null} and~\ref{sec:block_extract}.
Suppose we have drawn sketches of $\mtx A$ and its adjoint as in (\ref{eq:rand_sample_setting}), using $s$ Gaussian random test vectors.

The first step is to construct sparsification matrices $\mtx E$ and $\mtx F$ by computing an interpolation matrix $\mtx T$ and partitioning indices $\msf R \cup \msf S = \msf B$, such that (\ref{eq:skel_lowrank_simult}) holds. This is accomplished by computing an interpolative decomposition (ID) of the combined sketch:
\begin{equation*}
\left[ \msf R \cup \msf S,\ \mtx T \right] = \texttt{id} \left( \mtx Y_{\msf B}' + \mtx Z_{\msf B}' \right), \qquad \text{where}\qquad \mtx Y_{\msf B}' = \mtx A_{\msf B \msf F} \mtx \Omega'_{\msf F},\ \mtx Z_{\msf B}' = \mtx A_{\msf F \msf B}^* \mtx \Psi'_{\msf F}
\end{equation*} with $\mtx \Omega'_{\msf F}, \mtx \Psi'_{\msf F}$ Gaussian random matrices, and \texttt{id} defined in \eqref{eq:id_def}.
Rather than explicitly constructing these structured test matrices, block nullification allows us to compute the required sketches by post-processing the original sketches from (\ref{eq:rand_sample_setting}):
\begin{equation}
\mtx Y_{\msf B}' = \mtx Y_{\msf B}\ \texttt{null} \begin{pmatrix} \mtx \Omega_{\msf B}\\ \mtx \Omega_{\msf N} \end{pmatrix}, \qquad \mtx Z_{\msf B}' = \mtx Z_{\msf B}\ \texttt{null} \begin{pmatrix} \mtx \Psi_{\msf B}\\ \mtx \Psi_{\msf N}\end{pmatrix},
\label{eq:nullYZ}
\end{equation}
for $s \ge | \msf B| + | \msf N | + k + p$. 

Next, to compute the sparse elimination matrices $\mtx L$ and $\mtx U$, we need to extract sparse subblocks, such as $\mtx X_{\msf R \msf N}$ and $\mtx X_{\msf N \msf R}$. In order to accomplish this efficiently, we would like to sketch the sparsified system (\ref{eq:orth_sparsifying}), where the interaction between $\msf R$ and $\msf F$ is approximately zero. We aim to generate sketches of the sparsified matrix by reusing the original sketches of $\mtx A$.
This is achieved using the fact that $\mtx E$ and $\mtx F$ are invertible. Specifically, for an appropriate permutation:
\begin{equation*}
\mtx A \approx \mtx E^{-1} \left(\begin{array}{cc|cc}
{\mtx X}_{\msf R \msf R} & {\mtx X}_{\msf R \msf S} & {\mtx X}_{\msf R \msf N} & \\[1mm]
{\mtx X}_{\msf S \msf R} & {\mtx A}_{\msf S \msf S} & {\mtx A}_{\msf S \msf N} &  {\mtx A}_{\msf S \msf F}\\ \hline
{\mtx X}_{\msf N \msf R} & {\mtx A}_{\msf N \msf S} & {\mtx A}_{\msf N \msf N} & {\mtx A}_{\msf N \msf F}\\[1mm]
& {\mtx A}_{\msf F \msf S} & {\mtx A}_{\msf F \msf N} & {\mtx A}_{\msf F \msf F}
\end{array}\right) \mtx F^{-1} := \mtx E^{-1} \hat{\mtx A} \mtx F^{-1}.
\end{equation*}
Thus, we can obtain sketches of $\hat{\mtx A}$ from the original sketches
(\ref{eq:rand_sample_setting}) by applying the transformations
\[
\hat{\mtx Y} := \mtx E\,\mtx Y 
= \mtx E\,\mtx A\,\mtx \Omega 
= \hat{\mtx A}\,(\mtx F^{-1}\mtx \Omega).
\]
We then define the modified test matrix
\begin{equation*}
\hat{\mtx \Omega} := \mtx F^{-1}\mtx \Omega,
\qquad \text{so that} \qquad
\hat{\mtx Y} = \hat{\mtx A}\,\hat{\mtx \Omega}.
\end{equation*}
Likewise, the sketch of the adjoint can be updated as
$$\hat{\mtx Z} := \mtx F^* \mtx Z, \qquad \hat{\mtx \Psi} := \mtx E^{-*} \mtx \Psi.$$ Although the test matrices are modified and are no longer Gaussian, the modifications are sparse, and the matrices $\mtx E, \mtx F$ are  well-conditioned and do not significantly bias the test matrices.
Now that sketches of $\hat{\mtx A}$ have been obtained, block extraction can be used to compute
\begin{equation}
\begin{pmatrix}{\mtx X}_{\msf R \msf R} & {\mtx X}_{\msf R \msf S} & {\mtx X}_{\msf R \msf N} \end{pmatrix} \approx 
\hat{\mtx Y}_{\msf R} \begin{pmatrix} \hat{\mtx \Omega}_{\msf B}\\\hat{\mtx \Omega}_{\msf N} \end{pmatrix}^{\dagger}, \qquad 
\begin{pmatrix}{\mtx X}^*_{\msf R \msf R} & {\mtx X}^*_{\msf S \msf R} & {\mtx X}^*_{\msf N \msf R} \end{pmatrix}
\approx \hat{\mtx Z}_{\msf R} \begin{pmatrix}\hat{\mtx \Psi}_{\msf B}\\ \hat{\mtx \Psi}_{\msf N} \end{pmatrix}^{\dagger}.
\label{eq:pinvYZ}
\end{equation}
Here, approximate equalities are introduced because entries $\hat {\mtx A}_{\msf R \msf F}$ are only approximately zero,
and their contribution is absorbed into the extracted subblocks.

Finally, we update the original sketch and test matrices to maintain sketches of the
diagonalized matrix $\widetilde {\mtx A}$ for the subsequent steps of the algorithm, using the 
formula
\begin{equation*}
\mtx A \approx \mtx V
\left(\begin{array}{cc|cc}
{\mtx X}_{\msf R \msf R} &  \\[1mm]
 & {\mtx X}_{\msf S \msf S} & {\mtx X}_{\msf S \msf N} &  {\mtx A}_{\msf S \msf F}\\ \hline
& {\mtx X}_{\msf N \msf S} & {\mtx X}_{\msf N \msf N} & {\mtx A}_{\msf N \msf F}\\[1mm]
& {\mtx A}_{\msf F \msf S}  & {\mtx A}_{\msf F \msf N} & {\mtx A}_{\msf F \msf F}
\end{array}\right) \mtx W := \mtx V\ \widetilde{\mtx A}\ \mtx W.
\end{equation*}
Then, the following transformations yield a sketch of $\widetilde{\mtx A}$:
\begin{equation}\label{eq:updateVW_sketch}
\widetilde{\mtx Y} := \mtx V^{-1}\mtx Y
= \mtx V^{-1}\mtx A\mtx \Omega
= \widetilde{\mtx A}\,(\mtx W\mtx \Omega).
\end{equation}
We define the modified test matrix as
\begin{equation}\label{eq:updateVW_test}
\widetilde{\mtx \Omega} := \mtx W\mtx \Omega,
\qquad \text{so that} \qquad
\widetilde{\mtx Y} = \widetilde{\mtx A}\,\widetilde{\mtx \Omega}.
\end{equation}
Likewise, the sketch of the adjoint can be updated as 
\begin{equation}
\widetilde{\mtx Z}:= \mtx W^{-*}\ \mtx Z \qquad \widetilde{\mtx \Psi} := \mtx V^*\ \mtx \Psi.
\label{eq:updateVW_adj}
\end{equation}
While the matrices $\mtx V$ and $\mtx W$ are neither unitary nor independent of $\mtx \Omega$ and $\mtx \Psi$, empirical evidence indicates that the updated test matrices retain sufficient randomness for accurate low-rank approximation throughout the factorization process. Although in principle, certain directions could be preferentially amplified or suppressed, extensive numerical experiments suggest that this effect is minimal in practice.

\subsection{Randomized Strong Recursive Skeletonization}
\label{sec:rsrs_multistep}

In this section, we describe how the randomized procedure introduced in Section \ref{sub:randomized_strong_skeletonization} can be recursively applied to compute an approximate factorization of $\mtx A$ using only matrix--vector products. The sketches drawn from $\mtx A$ are updated throughout the algorithm, and structure can be introduced into the test
matrices as needed to sketch subblocks of the sparsified matrix at each stage.
Algorithm \ref{a:precompute} provides pseudocode. At the finest level of the tree (level $L$), all boxes are leaves, and the randomized diagonalization is applied box-by-box exactly as in Section~\ref{sec:strong_rec_skel}, except that all quantities are obtained from sketches rather than explicit matrix entries.

The randomized algorithm proceeds by repeatedly updating the sketches so that, at each box, they represent the current sparsified system rather than the original matrix. 
After diagonalizing the redundant indices for leaf box $B_1$, we can update the sketches to instead maintain sketches of the sparsified matrix
\begin{equation*}
\widetilde{\mtx Y}\ =\ \widetilde{\mtx A}(\mtx A; B_1)\ \widetilde{\mtx \Omega}, \qquad 
\widetilde{\mtx Z}\ =\ \widetilde{\mtx A}(\mtx A; B_1)^*\ \widetilde{\mtx \Psi}
\end{equation*}
using formulas (\ref{eq:updateVW_sketch}--\ref{eq:updateVW_adj}) to update initially drawn sketches of $\mtx A$.
As mentioned in Remark \ref{remark:additive_lowrank}, the diagonalization matrices
may introduce additive terms in the far-field of adjacent boxes that must be recompressed
at later stages of the algorithm. As such, to compute matrices
$\mtx V_2^{-1}, \mtx W_2^{-1}$ for the next leaf box $B_2$,
we use block nullification and extraction techniques on the \textit{updated} sketches
of $\{\widetilde{\mtx Y}, \widetilde{\mtx \Omega}, \widetilde{\mtx Z}, \widetilde{\mtx \Psi}\}$.
In order to maintain sketches of diagonalized matrix $\widetilde{\mtx A}(\mtx A; B_1, B_2)$,
we again apply update formulas (\ref{eq:updateVW_sketch}--\ref{eq:updateVW_adj}).

Let $n_L$ denote the number of boxes on the finest level $L$.
Then, after diagonalizing redundant indices $\msf R_1, \dots, \msf R_{n_L}$,
the accumulated updates of the original matrix $\mtx A$ take the form
\begin{equation*}
{\widetilde {\mtx Y}} := \mtx V^{-1}_{n_L}\ \cdots\ \mtx V^{-1}_{1}\ \mtx Y, \qquad \widetilde{\mtx \Omega}\ :=\ \mtx W_{n_L}\ \cdots\ \mtx W_{1}\ \mtx \Omega,
\end{equation*}
and likewise for the adjoint sketch and test matrices.
The test matrix update involves the sparsification matrices $\mtx E, \mtx F$, which are also typically well-conditioned, and the block elimination matrices $\mtx L, \mtx U$. These elimination factors are often well-conditioned as well, since they are unit-triangular with modest off-diagonal entries. As a result, the majority of the ill-conditioning in $\mtx A$ is isolated in the diagonal factor $\mtx D$ of the factorization in (\ref{eq:Afact}). 

The sparsification matrices $\mtx E, \mtx F$ and elimination matrices $\mtx L, \mtx U$ are constructed from the sketches of $\mtx A$, and therefore depend implicitly on the choice of test matrices $\mtx \Omega$ and $\mtx \Psi$. 
In principle, this dependence introduces a potential risk: the sketching matrices used to approximate $\mtx A$ are also used to define the transformations that update those sketches. In practice, however, this feedback effect appears to be benign. The sparsification and elimination matrices are sparse and well-conditioned, and as a result, the updated sketches $(\widetilde{\mtx Y}, \widetilde{\mtx \Omega})$ remain sufficiently rich to capture the essential subspace information of the sparsified matrix $\hat{\mtx A}$, without needing to redraw independent sketches at each stage of the algorithm.

\begin{algorithm}[htbp!]
\caption{\texttt{Randomized strong recursive skeletonization}}
\label{a:precompute}
\begin{algorithmic}[1]
\Require{Sketch and test matrices $\{\mtx Y, \mtx Z, \mtx \Omega, \mtx \Psi\}$ of size $N \times s$ satisfying (\ref{eq:rand_sample_setting}).}
\Ensure An invertible factorization $\mtx A_{\rm approx} \approx \mtx A$.

\State Initialize all indices as active:
$$\msf {active} = [1, \dots, N].$$
\For{level $\ell = L$ down to $1$}  \Comment{Upward traversal of the tree.}
  \For{each box $B$ at level $\ell$}
        \State Get relevant active indices:
        $$\msf B_{\rm active} = \msf B \cap \msf{active}, \qquad \msf N_{\rm active} = \msf N \cap \msf{active}.$$

        \State Compute far-field sketches using block nullification:
        \[
        \mtx Y'_{\msf B_{\rm active}} = \mtx Y_{\msf B_{\rm active}} \ \texttt{null} \begin{pmatrix} \mtx \Omega_{\msf B_{\rm active}} \\ \mtx \Omega_{\msf N_{\rm active}} \end{pmatrix},
        \]
        \[
        \mtx Z'_{\msf B_{\rm active}} = \mtx Z_{\msf B_{\rm active}} \ \texttt{null} \begin{pmatrix} \mtx \Psi_{\msf B_{\rm active}} \\ \mtx \Psi_{\msf N_{\rm active}} \end{pmatrix}.
        \]

        \State Compute the row-ID:
        $$
        [\msf R \cup \msf S, \mtx T] = \texttt{id}(\mtx Y'_{\msf B_{\rm active}} + \mtx Z'_{\msf B_{\rm active}}).$$

        \State Form sparsification matrices $\mtx E, \mtx F$ and update sketches:
		\begin{alignat*}{2}
		\hat{\mtx Y}     &\gets \mtx E\ \mtx Y,     &\quad \hat{\mtx \Omega} &\gets \mtx F^{-1}\ \mtx \Omega, \\
		\hat{\mtx Z}     &\gets \mtx F^*\ \mtx Z,   &\quad \hat{\mtx \Psi}   &\gets \mtx E^{-*}\ \mtx \Psi.
		\end{alignat*}
        \State Extract sparse blocks using block extraction:
        \[
        (\mtx X_{\msf R \msf B_{\rm active}}\quad  \mtx X_{\msf R \msf N_{\rm active}}) \approx \mtx Y_{\msf R} \begin{pmatrix} \mtx \Omega_{\msf B_{\rm active}} \\ \mtx \Omega_{\msf N_{\rm active}} \end{pmatrix}^{\dagger},
        \]
        \[
        (\mtx X^*_{\msf B_{\rm active} \msf R}\quad  \mtx X^*_{\msf N_{\rm active} \msf R}) \approx \mtx Z_{\msf R} \begin{pmatrix} \mtx \Psi_{\msf B_{\rm active}} \\ \mtx \Psi_{\msf N_{\rm active}} \end{pmatrix}^{\dagger}.
        \]

        \State Compute sparse elimination matrices $\mtx L, \mtx U$ and update sketches:
		\begin{alignat*}{2}
		\widetilde{\mtx Y}     &\gets \mtx L\ \hat{\mtx Y}, \quad & \widetilde{\mtx \Omega} &\gets \mtx U^{-1}\ \hat{\mtx \Omega}, \\
		\widetilde{\mtx Z}     &\gets \mtx U^* \ \hat{\mtx Z}, \quad & \widetilde{\mtx \Psi} &\gets \mtx L^{-*}\ \hat{\mtx \Psi}.
		\end{alignat*}
        \State Remove $\msf R$ from active index set:
        $$\msf {active} = \msf {active} \setminus \msf R.$$
        \State Assign updated sketch and test matrices:
        \[
        \mtx Y\ \gets\ \widetilde{\mtx Y}, \quad \mtx \Omega\ \gets\ \widetilde{\mtx \Omega}, \quad \mtx Z\ \gets\ \widetilde{\mtx Z}, \quad \mtx \Psi\ \gets\ \widetilde{\mtx \Psi}.
        \]

    \EndFor
\EndFor

\State At the root level, extract the final remaining submatrix:
$$\widetilde{\mtx A}_{\msf B_t \msf B_t} \approx \widetilde{\mtx Y}_{\msf B_t} \ \widetilde{\mtx \Omega}_{\msf B_t}^{\dagger},$$
where $\msf B_t$ is an index vector which denotes the remaining active indices in the domain.

\end{algorithmic}
\end{algorithm}
 
To generalize the procedure to a multilevel setting, we maintain an active set of indices as defined in (\ref{eq:active}) and modify the block nullification and extraction formulas accordingly. When skeletonizing a box $B$, we use the \textit{active} box and neighbor indices, $\msf B_{\text{active}}$ and $\msf N_{\text{active}}$, as defined in (\ref{eq:BNactive}), instead of the original full index sets.
In particular, the null space computation used to form the basis matrices in (\ref{eq:nullYZ}) becomes
\begin{equation*}
\widetilde{\mtx Y}_{\msf B_{\rm active}}' = \widetilde{\mtx Y}_{\msf B_{\rm active}}\ \texttt{null} \begin{pmatrix} \widetilde{\mtx \Omega}_{\msf B_{\rm active}}\\ \widetilde{\mtx \Omega}_{\msf N_{\rm active}} \end{pmatrix}
\end{equation*}
for modified test matrices $\widetilde{\mtx Y}, \widetilde{\mtx \Omega}$.
The block extraction computation in \eqref{eq:pinvYZ} is likewise restricted to these active index sets.

At coarser levels of the hierarchy, this computation remains efficient because the number of active near-field indices is relatively small. However, the null space and block extraction operations implicitly define a modified test matrix $\widetilde{\mtx \Omega}'$ that has support on \textit{inactive} points. Although the corresponding entries of the diagonalized matrix $\widetilde{\mtx A}$ are small (since inactive points have already been eliminated), they are not exactly zero. Consequently, the sketches formed at coarser levels incorporate these small inactive contributions, and the resulting residual errors accumulate as the algorithm progresses up the tree.

One possible approach to mitigating this error accumulation would be to redraw structured test matrices at each level, explicitly controlling for inactive contributions. However, in this work, we deliberately focus on the setting where all samples are drawn once and reused throughout the factorization process. Determining when and how to optimally redraw sketches to control error propagation across levels remains an important direction for future research. 
\subsection{Algorithm Complexity}
\label{sec:complexity}

In this section, we analyze both the sample and computational complexity of Algorithm~\ref{a:precompute}.
We consider a point distribution consisting of $N$ total points in $d$ dimensions ($d=2,3)$ organized hierarchically in multilevel tree with $L$ levels, with at most $m$ points per leaf box. The total number of boxes is approximately
\begin{equation}
n_{\rm boxes} = n_{\rm leaf} + n_{\rm tree}, \qquad n_{\rm leaf} \approx \frac{N}{m}, \quad n_{\rm tree} \approx \frac{1}{2^d - 1} \cdot \frac{N}{m}.
\label{eq:nboxes}
\end{equation}
We denote by $c_{\rm nei}(B)$ the number of neighboring boxes for a given target box $B$ (including the box itself). This constant depends on the geometry and remain uniform across boxes in regular grids. For example, for uniform point distributions in $d$ dimensions, $c_{\rm nei} = 3^d$ for all boxes.

For clarity, we first analyze the case of uniform point distributions in Section~\ref{ssub:uniform_point_distributions}, and then consider the more general case of nonuniform distributions in Section~\ref{ssub:nonuniform_point_distributions}. Throughout, we assume that far-field interactions can be compressed to a fixed rank $k$, specified by the user. We denote by $p$ the oversampling parameter used in randomized compression (cf. Section \ref{sec:randomized_lr}).

\subsubsection{Uniform Point Distributions}
\label{ssub:uniform_point_distributions}

At the finest level of the tree (i.e., level $L$), the dominant costs arise from block nullification and block extraction. For instance, block nullification involves computing the null space of a matrix of size $c_{\rm nei}\ m \times s$. To produce a structured test matrix $\mtx \Omega'$ with $(k + p)$ columns, the number of samples required is
\begin{equation}
s_L \geq c_{\rm nei}\ m + k + p, \qquad \text{for a box on level $L$}.
\label{eq:s_L}
\end{equation}
The computational complexity of the null-space operation per box is then
$\mathcal{O}\left(c_{\rm nei}^2\ m^2\ s\right) = \mathcal{O}\left(c_{\rm nei}^3\ m^3\right)$, because $s_L = \mathcal{O}\left(c_{\rm nei}\ m\right)$.

At coarser levels of the tree, only the active degrees of freedom associated with the box and its neighbors are involved; 
therefore, there are $2^d \cdot c_{\rm nei} k$ active points per tree box, because a tree box has $2^d$ children. The number of samples required for a tree box is
\begin{equation}
s_B \geq 2^d\ c_{\rm nei}\ \ k + k + p, \qquad \text{for tree box $B$}.
\label{eq:s_B}
\end{equation}
Because the samples drawn initially are reused for \textit{all} boxes, we choose $m$ so that
\eqref{eq:s_L} does not dominate \eqref{eq:s_B} as
\begin{equation}
m = 2^d\ k
\label{eq:m_choice_uniform}
\end{equation}
which leads to overall sampling requirement
\begin{equation}
s =\ \max_B s_B =\ 2^d\ c_{\rm nei}\ k + k + p = (6^d + 1)\ k + p.
\label{eq:s_uniform}
\end{equation}
For uniform point distributions, the product $2^d \cdot 3^d = 6^d$ in \eqref{eq:s_uniform} arises frequently in fast multipole methods and $\mathcal{H}^2$ matrices under strong admissibility conditions, where it is associated with the size of the interaction list.

The total post-processing time to reconstruct the factorization accounts for processing all boxes in the hierarchy.
Using (\ref{eq:nboxes}), the total cost is
\begin{equation}
T_{\rm rec} = \mathcal{O}\!\left(c_{\rm nei}^3\ m^3\ \frac{N}{m}\right)
            = \mathcal{O}\!\left(c_{\rm nei}^3\ (2^d k)^2\ N\right)
            = \mathcal{O}\!\left(108^d k^2 N\right).
\label{eq:Trec_uniform}
\end{equation}
The constant $108^d$ arises from the product $(3^d)^3 (2^d)^2 = 108^d$. 
This complexity is linear in the problem size, albeit with large constants that depend on both dimensionality
and the chosen rank parameter.
The cost of applying the computed factorization to a vector (i.e., a matrix--vector solve) scales linearly with $N$ as
\begin{equation}
T_{\rm sol} = \mathcal{O}\left(6^d\ k\ N\right).
\label{eq:Tsol_uniform}
\end{equation}
A key advantage of Algorithm~\ref{a:precompute} is that it requires only a single set of random test matrices, which are reused across all levels via post-processing. This sample efficiency makes it particularly effective for problems where matrix entries are expensive or inaccessible. While the post-processing steps, such as block nullification and block extraction, introduce a large constant prefactor to the runtime, this overhead is often offset by the substantial savings in sample complexity. 

\subsubsection{Nonuniform Point Distributions}
\label{ssub:nonuniform_point_distributions}

For nonuniform point distributions, the number of active degrees of freedom per box can vary significantly with the local point density and the geometry of the domain. Consequently, the sample complexity is no longer uniform across all boxes and must adapt to the local number of active degrees of freedom. We set the leaf-box capacity as
\[
m = c_{\rm cap}\, k, \qquad c_{\rm cap} \le 2^d
\]
where $c_{\rm cap}$ is a user-defined constant (e.g.\ we choose $c_{\rm cap} = 6$ in Section \ref{sec:numerical_results} for surface geometries in three dimensions, whereas for volumetric point distributions we set $c_{\rm cap} = 8$, consistent with \eqref{eq:m_choice_uniform}).

For a given box $B$, the number of samples $s_B$ required is determined by the active degrees of freedom in neighboring boxes. Under the 2:1 balance condition, a neighboring box $B'$ can be either a tree box on the same level or a leaf box one level above. Denoting by $\msf B'$ the index vector for the neighboring box $B'$, we may write schematically
\[
s_B = (k + p)
      + \sum_{\substack{\text{tree box } B' \\ \text{neighboring } B}}
        \min\ \bigl(k,\, \underbrace{\left |\msf B' \cap \msf{active} \right|}_{\substack{\text{number of } \\ \text{active points in } B'}}\bigr)
      + \sum_{\substack{\text{leaf box } B' \\ \text{neighboring } B}}\
        \underbrace{|\msf B'|}_{\le m}
\]
with the $\msf{active}$ index vector defined in \eqref{eq:active}.
This yields an upper bound on the total number of samples needed required:
\begin{equation*}
s = \max_B s_B = (c_{\rm geom}\, + 1)\ k + p,
\end{equation*}
where $c_{\rm geom}$ is a constant that depends only on the geometry, neighbor structure, and the choice of $c_{\rm cap}$ (and independent of $N$).
Likewise, $$T_{\rm sol} = \mathcal O\left(c_{\rm geom}\ k\ N\right),$$ analogous to \eqref{eq:Tsol_uniform}.
The total cost of post-processing the samples across all boxes is
\[
T_{\rm rec}
= \mathcal{O}\!\left(s^3 \frac{N}{m}\right)
= \mathcal{O}\!\left(\frac{(c_{\rm geom} k)^3}{c_{\rm cap} k}\, N\right)
= \mathcal{O}\!\left(\frac {c_{\rm geom}^3}{c_{\rm cap}}\, k^2 N\right).
\]
These complexity estimates are consistent with estimate for the uniform case in (\ref{eq:s_uniform}--\ref{eq:Tsol_uniform}) with $c_{\rm geom} = 6^d$ and $c_{\rm cap} = 2^d$.
In many practical settings—particularly when points lie on a lower-dimensional manifold (e.g., a surface in 3D)—the effective constant $c_{\rm geom}$ is modest, and the nonuniform case can be much cheaper than the uniform $d$-dimensional case. 
\section{Numerical Results}
\label{sec:numerical_results}

In this section, we report the performance of the proposed \texttt{RSRS} algorithm on several model problems in two and three dimensions. 
The algorithm depends on two key parameters: the target rank~$k$, which determines the approximation accuracy, and the oversampling parameter~$p$ (cf. Section \ref{sec:randomized_lr}). 
Throughout all experiments, we set $p = k$, providing a conservative oversampling margin that yields consistent accuracy across a range of problems.
The algorithm operates on a hierarchical tree~$\mathcal{T}$, whose leaves contain at most~$m$ points, chosen as $m= c_{\rm cap}\ k$
for a geometry-dependent constant $c_{\rm cap}$ (cf. Section \ref{ssub:nonuniform_point_distributions}).

The accuracy of the computed factorization~$\mtx A_{\rm approx} \approx \mtx A$ is evaluated using two standard metrics, which are estimated by power iteration:
\begin{equation}
\label{eq:error}
\mathrm{relerr}
= \frac{\|\mtx A - \mtx A_{\rm approx}\|}{\|\mtx A\|},
\qquad
\mathrm{errsolve}
= \|\mtx I - \mtx A_{\rm approx}^{-1}\mtx A\|.
\end{equation}
Here, $\mathrm{relerr}$ quantifies the accuracy of the low-rank compression,
while $\mathrm{errsolve}$ measures the effectiveness of the resulting approximate inverse.
When $\mtx A_{\rm approx}^{-1}$ is used as a direct solver or preconditioner, 
its effectiveness depends on the conditioning of~$\mtx A$. 
If $\mathrm{relerr} \le \epsilon$, then for the exact solution $\mtx x = \mtx A^{-1}\mtx b$ and the approximate solution $\tilde{\mtx x} = \mtx A_{\rm approx}^{-1}\mtx b$, 
the relative error satisfies
\begin{equation}
\frac{\|\mtx x - \tilde{\mtx x}\|}{\|\mtx x\|}
\le
\frac{2\,\epsilon\,\mathrm{cond}(\mtx A)}{1 - \epsilon\,\mathrm{cond}(\mtx A)}.
\label{eq:err_solve_cond}
\end{equation}
Hence, problems with larger condition numbers generally require larger ranks~$k$ to maintain accuracy in the solver.

We consider three representative problems that illustrate the performance of \texttt{RSRS} for both integral and differential equations.
These include a first-kind volume integral equation on a uniform 2D grid (Section~\ref{sec:exp1}), 
a second-kind boundary integral equation on a sphere (Section~\ref{sec:exp2}), 
and a Schur complement arising in a sparse direct solver for the discretized 3D Helmholtz equation (Section~\ref{sec:schur_compression}).
We set the leaf parameter $m= c_{\rm cap}k$ with $c_{\rm cap}=4$ for the uniform grid in Section \ref{sec:exp1}
and $c_{\rm cap}=6$ for the surface geometries in Section \ref{sec:exp2} and Section \ref{sec:schur_compression}.

In Experiments 1 and 3, the condition number $\mathrm{cond}(\mtx A)$ increases with problem size~$N$ 
or with the wavenumber, as in the Helmholtz case.
Despite this, the numerical results show that \texttt{RSRS} maintains stable accuracy and efficient scaling even for moderately ill-conditioned or indefinite systems.
In addition to the error metrics above, we report the number of GMRES iterations required to solve $\mtx A\mtx x = \mtx b$, both with and without the preconditioner $\mtx A_{\rm approx}^{-1}$. 

\begin{table}[!htb]
\centering \tablesize
\renewcommand{\arraystretch}{1.2}
\begin{tabular}{@{}ll@{}}
\toprule
Notation & Description \\
\midrule
$N$ & Number of discretization points \\[2pt]
$k,p$ & Rank and oversampling parameters, respectively \\[2pt]
$m$ & Leaf size parameter \\[2pt]
$n_{\mathrm{samples}}$ & Number of matrix–vector samples of $\mathbf A$ and $\mathbf A^*$ \\[2pt] \midrule
$T_{\mathrm{rec}}$ & Wall-clock time (seconds) to reconstruct $\mathbf A_{\mathrm{approx}}^{-1}$  \\[2pt]
$M$ & Memory (GB) to store $\mathbf A_{\mathrm{approx}}^{-1}$ \\[2pt]
$T_{\mathrm{sol}}$ & Wall-clock time (seconds) to apply $\mathbf A_{\mathrm{approx}}^{-1}$ \\[2pt] \midrule
$\mathrm{relerr}$, $\mathrm{errsolve}$ & Defined in equation~(\ref{eq:error}) \\[2pt]
$\mathrm{relerr}_{\mathrm{bvp}}$ & \makecell[l]{Relative error when using \texttt{RSRS}-based solver 
\\to solve BVP, e.g. \eqref{eq:lpbvp} and \eqref{eq:bvp_helm}, and measured as \eqref{eq:relerr_bvp}  } \\[2pt] \midrule
$n_{\rm iter}$ & \makecell[l]{Number of GMRES iterations required \\ to reach relative tolerance $\mathrm{rtol}=10^{-10}$} \\
\bottomrule
\end{tabular}
\caption{Summary of the notations used in the reported numerical results.}
\label{tab:notation_report}
\end{table}

We also measure the accuracy of \texttt{RSRS} when used as a direct solver for boundary-value problems of the form
\[
\mathcal{A} u(\pxx) = f(\pxx) \quad (\pxx \in \Omega),
\qquad
u(\pxx) = g(\pxx) \quad (\pxx \in \partial\Omega),
\]
where $\mathcal{A}$ is an elliptic operator such as those in~\eqref{eq:lpbvp} and~\eqref{eq:bvp_helm}. 
For a computed potential~$\hat{u}$ and an analytic reference solution~$u$, 
the relative boundary-value error is defined as
\begin{equation}
\mathrm{relerr}_{\mathrm{bvp}}
=
\frac{\|\hat{u}(\pxx_t) - u(\pxx_t)\|}{\|u(\pxx_t)\|},
\qquad \pxx_t \in \Omega,
\label{eq:relerr_bvp}
\end{equation}
where the potentials are evaluated at interior target points.

All experiments were performed on a workstation equipped with an Intel Xeon Gold~6254 CPU and~768~GB of RAM. 
Although the processor supports multiple cores, the current implementation employs limited parallelism---each box is processed sequentially.

The notations used are summarized in Table \ref{tab:notation_report}. Throughout the paper, we use $s$ to denote the number of matrix--vector
samples needed for \texttt{RSRS}. In the numerical results, we report this quantity as
$n_{\mathrm{samples}}$ for readability; the two notations refer to the same
quantity.

\subsection{Experiment 1: Integral Equation on a 2D Grid }
\label{sec:exp1}

Consider a volume integral for the Laplace equation discretized on
a unit square $\Omega = [0,1]^2$.  The collocation points $\pxx$ are placed on a uniform grid
of $\sqrt N \times \sqrt N$ points, and matrix entries are given by
\begin{equation}
\begin{aligned}
\mtx A_{ij} &= \frac 1 N \log(\| \pxx_i - \pxx_j\|),\qquad i \neq j\\ 
\mtx A_{ii} &\approx \int_{-h/2}^{h/2} \int_{-h/2}^{h/2} \log(\|\pxx_i - \pxx_j\|)\ d \pxx_i\ d \pxx_j,\qquad \mbox{where}\ h \equiv 1/\sqrt{N},
\end{aligned}
\label{eq:Aentries_exp1}
\end{equation}
and where modified entries on the diagonal are approximated using the Quadpack library \cite{piessens2012quadpack}.
Because the collocation points are on a uniform grid, the matrix $\mtx A$ and its adjoint can be applied
to vectors in $\mathcal O(N\log N)$ time using the FFT \cite{cooley1965algorithm}.
This example appeared previously in \cite{2017_ho_ying_strong_RS}.

\begin{figure}[htb!]
  \centering

\begin{subfigure}{0.33\textwidth}
    \centering
    \begin{tikzpicture}
      \begin{loglogaxis}[
        width=0.78\linewidth, height=0.78\linewidth,
        axis lines=box,
        tick label style={font=\figlabelsize},
        label style={font=\figlabelsize},
        legend style={font=\figlabelsize, draw=none},
        minor tick num=9,
        xlabel={$N$},
        ylabel={$T_{\mathrm{rec}}\ (s)$},
        ylabel style={
        at={(0.22,0.6)},
        anchor=south,
        font=\figlabelsize,
        },
        xlabel style={
        at={(0.5,-0.2)},
        anchor=south,
        font=\figlabelsize,
        }      ]
\addplot+[thick, mark=*, mark size=1.8pt]
          table[x=N, y={Trec_k40}, col sep=comma]{rev_laplace_grid.csv};

        \addplot+[thick, mark=square*, mark size=1.8pt]
          table[x=N, y={Trec_k60}, col sep=comma]{rev_laplace_grid.csv};
      \addONguideRel{0.5}{0.25}{0.3}
      \end{loglogaxis}
    \end{tikzpicture}
  \end{subfigure}\hfill
\begin{subfigure}{0.33\textwidth}
    \centering
    \begin{tikzpicture}
      \begin{loglogaxis}[
        width=0.78\linewidth, height=0.78\linewidth,
        axis lines=box,
        tick label style={font=\figlabelsize},
        label style={font=\figlabelsize},
        legend style={font=\figlabelsize, draw=none},
        minor tick num=9,
        xlabel={$N$},
        ylabel={$M$ (GB)},
        ylabel style={
        at={(0.22,0.6)},
        anchor=south,
        font=\figlabelsize,
        },
        xlabel style={
        at={(0.5,-0.2)},
        anchor=south,
        font=\figlabelsize,
        }
      ]
        \addplot+[thick, mark=*, mark size=1.8pt]
          table[x=N, y={M_GB_k40}, col sep=comma]{rev_laplace_grid.csv};

        \addplot+[thick, mark=square*, mark size=1.8pt]
          table[x=N, y={M_GB_k60}, col sep=comma]{rev_laplace_grid.csv};
      \addONguideRel{0.5}{0.25}{0.3}
      \end{loglogaxis}
    \end{tikzpicture}
  \end{subfigure}\hfill
\begin{subfigure}{0.33\textwidth}
    \centering
    \begin{tikzpicture}
      \begin{loglogaxis}[
        legend to name=sharedlegend_exp1,
        legend columns=2,
        width=0.78\linewidth, height=0.78\linewidth,
        axis lines=box,
        tick label style={font=\figlabelsize},
        label style={font=\figlabelsize},
        legend style={font=\figlabelsize, draw=gray,    column sep=1.2em,
  row sep=0.3em,
  inner xsep=0.8em,
  legend cell align={left}},
        minor tick num=9,
        xlabel={$N$},
        ylabel={$T_{\mathrm{sol}}\ (s)$},
        ylabel style={
        at={(0.22,0.6)},
        anchor=south,
        font=\figlabelsize,
        },
        xlabel style={
        at={(0.5,-0.2)},
        anchor=south,
        font=\figlabelsize,
        }
      ]
        \addplot+[thick, mark=*, mark size=1.8pt]
          table[x=N, y={Tsolve_k40}, col sep=comma]{rev_laplace_grid.csv};
        \addlegendentry{\shortstack{$k=40,\ p=40$\\$n_{\mathrm{samples}}=1520$}}

        \addplot+[thick, mark=square*, mark size=1.8pt]
          table[x=N, y={Tsolve_k60}, col sep=comma]{rev_laplace_grid.csv};
        \addlegendentry{\shortstack{$k=60,\ p=60$\\$n_{\mathrm{samples}}=2280$}}
      \addONguideRel{0.5}{0.25}{0.3}
      \end{loglogaxis}
    \end{tikzpicture}
  \end{subfigure}

\vspace{0.5em}

\centering

\noindent
\hspace*{0.10\textwidth}\begin{minipage}[c]{0.16\textwidth}
\centering
\begingroup
\resizebox{\linewidth}{!}{\begin{tikzpicture}[scale=3]

\def\N{12}            \def\dotsize{0.012}   \def\rim{black!40}    \def\dotcol{black!65} 

\begin{scope}
    \clip (0,0) rectangle (1,1);
\shade[left color=white, right color=black!30] (0,0) rectangle (1,1);
  \end{scope}

\draw[\rim, line width=0.45pt] (0,0) rectangle (1,1);

\foreach \i in {0,...,\N} {
    \pgfmathsetmacro\x{1.0*\i/\N}
    \foreach \j in {0,...,\N} {
      \pgfmathsetmacro\y{1.0*\j/\N}
      \fill[\dotcol] (\x,\y) circle (\dotsize);
    }
  }

\end{tikzpicture}
 }
\endgroup
\end{minipage}\hfill
\begin{minipage}[c]{0.68\textwidth}
\centering
\pgfplotslegendfromname{sharedlegend_exp1}\\[0.4em]

\underline{\textbf{Experiment 1}} \hspace{1em} {\figtitlesize $n_{\mathrm{samples}} = 37\,k + p$}
\end{minipage}

  \caption{Reconstruction time, memory usage, and solve time vs.\ degrees of freedom for matrix $\mtx A$ arising
  from the discretization of a boundary integral equation with entries given in (\ref{eq:Aentries_exp1}).}
  \label{fig:exp1}
\end{figure} 
\begin{table}[htb!]
\centering

\begin{subtable}[t]{\linewidth}
    \centering \tablesize
    \caption{Accuracy and solver effectiveness.}
    \noindent \clearpage{}\begin{tabular}{rc cc cc}
\toprule
\multicolumn{2}{c}{} & \multicolumn{2}{c}{$k=40,\ p = 40$} & \multicolumn{2}{c}{$k=60,\ p = 60$} \\\cmidrule(lr){3-4} \cmidrule(lr){5-6}
\multicolumn{1}{c}{$N$} & $\rm cond(\mtx A)$ & $\mathrm{relerr}$ & $\mathrm{errsolve}$ & $\mathrm{relerr}$ & $\mathrm{errsolve}$ \\
\midrule
80,656 & 1.2e+05 & 2.0e-07 & 9.8e-04 & 1.3e-09 & 1.5e-05 \\
160,801 & 2.5e+05 & 2.0e-07 & 1.1e-03 & 3.0e-09 & 5.6e-05 \\
321,489 & 5.1e+05 & 2.3e-07 & 3.4e-03 & 1.4e-08 & 1.1e-04 \\
641,601 & 1.0e+06 & 7.6e-07 & 7.3e-03 & 2.1e-08 & 3.2e-04 \\
1,281,424 & 2.1e+06 & 1.3e-06 & 1.3e-02 & 6.8e-08 & 4.1e-04 \\
\bottomrule
\end{tabular}

\clearpage{}
    \label{tab:laplace_grid_acc}
\end{subtable}
\\
\begin{subtable}[t]{\linewidth}
    \centering \tablesize
    \caption{Preconditioner performance.}
    \noindent \clearpage{}\begin{tabular}{r r cc cc}
\toprule
 & {No precond} & \multicolumn{2}{c}{$k=40,\ p = 40$} & \multicolumn{2}{c}{$k=60,\ p = 60$} \\ \cmidrule(lr){2-2} \cmidrule(lr){3-4} \cmidrule(lr){5-6}
\multicolumn{1}{c}{$N$} & \multicolumn{1}{c}{$n_{\rm iter}$} & $n_{\rm samples}$ & $n_{\rm iter}$ & $n_{\rm samples}$ & $n_{\rm iter}$ \\
\midrule
80,656 & 4,910 & 1,520 & 3 & 2,280 & 2 \\
160,801 & 9,734 & 1,520 & 3 & 2,280 & 2 \\
321,489 & $>$10,000 & 1,520 & 4 & 2,280 & 3 \\
641,601 & $>$10,000 & 1,520 & 4 & 2,280 & 3 \\
1,281,424 & $>$10,000 & 1,520 & 5 & 2,280 & 3 \\
\bottomrule
\end{tabular}

\clearpage{}
    \label{tab:laplace_grid_niter}
\end{subtable}
\caption{Summary of results for Experiment 1. 
(a): Accuracy of \texttt{RSRS} and its effectiveness as a solver across problem sizes for two rank choices. 
(b): Preconditioner performance, showing that despite worsening conditioning, \texttt{RSRS} remains effective.}
\label{tab:laplace_summary}
\end{table}

This matrix corresponds to a first-kind integral equation and becomes increasingly ill-conditioned as $N$ grows:
the reported condition number rises from approximately $10^5$ to above $10^6$.
Even so, the matrix $\mtx A$ is still highly compressible, and the \texttt{RSRS} factorization maintains a small relative compression error,
with $\mathrm{relerr}$ in the range $10^{-7}$--$10^{-6}$ for $k=40$ 
and improving by one to two orders of magnitude for $k=60$ (Table~\ref{tab:laplace_grid_acc}).

Figure~\ref{fig:exp1} summarizes the measured quantities and sample costs for each rank parameter~$k$, while Table~\ref{tab:laplace_summary} reports the solver and preconditioner performance.
The inverse quality metric $\mathrm{errsolve}$ remains below $10^{-4}$ for $k=60$  across all problem sizes,
indicating that the approximate inverse is sufficiently accurate to support robust performance as a direct solver
despite the ill-conditioning of $\mtx A$. 
The discrepancy between $\mathrm{relerr}$ and $\mathrm{errsolve}$ is consistent with the bound in \eqref{eq:err_solve_cond}.

Without preconditioning, GMRES requires thousands of iterations and fails to converge within the iteration budget
for larger $N$.
In contrast, the \texttt{RSRS}-based preconditioner reduces the iteration count to between two and five iterations
for all problem sizes and both rank choices (Table~\ref{tab:laplace_grid_niter}).

The timings in Figure~\ref{fig:exp1} show that the reconstruction time $T_{\mathrm{rec}}$,
solve time $T_{\mathrm{sol}}$, and memory footprint $M$ grow approximately linearly with $N$,
consistent with the $\mathcal{O}(k^2 N)$ and $\mathcal{O}(k N)$ complexities analyzed in
Section~\ref{ssub:uniform_point_distributions}. The number of samples required is $n_{\rm samples} = 37\,k + p$, which also matches the complexity analysis.
 
\subsection{Experiment 2: Second-Kind Boundary Integral Equation on Sphere Surface}
\label{sec:exp2}

We consider the Dirichlet problem for the Laplace equation inside the unit sphere, which we denote $\Omega$. Given Dirichlet data \( f \) defined on the surface of the unit sphere  $\partial \Omega$, we aim to find the solution $u$ which satisfies
\begin{equation}
\Delta u(\pxx) = 0 \quad \text{for } \pxx \in \Omega, \qquad u(\pxx) = g(\pxx) \quad \text{for } \pxx \in \partial\Omega.
\label{eq:lpbvp}
\end{equation}
To solve (\ref{eq:lpbvp}), we represent the solution using a {double-layer potential}, for some function $\sigma(\pxx)$ defined on the boundary
\begin{equation}
u(\pxx) = \int_{\partial \Omega} \frac{\partial G}{\partial n_y}(\pxx - \pyy) \, \sigma(\pyy) \, d\pyy,\qquad \pxx \in \Omega
\label{eq:ansatz}
\end{equation}
where \( G(\pxx) = \frac{1}{4\pi \|\pxx\|} \) is the free-space Green’s function for the Laplace equation in \( \mathbb{R}^3 \), \( \frac{\partial G}{\partial n_{\pyy}}\) denotes the normal derivative with respect to the source variable \( \pyy \), and \( n_{\pyy} \) is the outward-pointing normal at point \( \pyy \in \partial \Omega \).
Imposing the Dirichlet condition yields a second-kind Fredholm integral equation to be solved for the unknown $\sigma(\pxx)$
\begin{equation}
- \frac{1}{2} \sigma(\pxx) + \int_{\partial \Omega} \frac{\partial G}{\partial n_{\pyy}}(\pxx - \pyy) \, \sigma(\pyy) \, d\pyy = g(\pxx), \qquad \pxx \in \partial \Omega.
\label{eq:bie}
\end{equation}
To discretize (\ref{eq:bie}), we represent \( \partial \Omega \) as a mesh of flat triangular panels. The centroids of these triangles are used as collocation points in the numerical method. 
To accurately evaluate integrals involving near-field interactions, we employ fourth-order tensor-product Gauss–Legendre quadrature.
The resulting linear system is not self-adjoint due to both the asymmetry introduced
by the double-layer kernel and by the Nystr\"om quadrature. The discretized system is applied rapidly to vectors using the FMM \cite{fmm3d}. This example previously appeared in \cite{2017_ho_ying_strong_RS}.

The boundary integral equation (\ref{eq:bie}) has the form of $\mtx I + \mtx G$, where the latter term $\mtx G$ is a compact operator. Because it is a compact perturbation of the identity operator, the resulting discretized system is well conditioned, and as reported in Table \ref{tab:exp2}, the condition number of the system is essentially constant as the number of discretization points $N$ increases. As a result, GMRES converges in only a few iterations, and a preconditioner is not necessary. We include this example as a numerical benchmark to demonstrate the performance of the algorithm without the effects of ill-conditioning.

We also report the effectiveness of using \texttt{RSRS} as a solver for
\eqref{eq:lpbvp} via the metric $\mathrm{relerr}_{\mathrm{bvp}}$.  
Point sources are placed outside the sphere, generating a known harmonic potential $u$ inside the sphere. 
The induced boundary values on~$\partial\Omega$ are used as Dirichlet data in~\eqref{eq:lpbvp}. 
The \texttt{RSRS}-based direct solver is then used to compute the surface
density $\sigma$ in \eqref{eq:bie}. From this density, the interior potential $\hat u$ is reconstructed as (\ref{eq:ansatz})
and compared against the exact analytic potential $u$ inside the sphere for interior target points as in \eqref{eq:relerr_bvp}.

\begin{figure}[htb!]
  \centering

\begin{subfigure}{0.33\textwidth}
    \centering
    \begin{tikzpicture}
      \begin{loglogaxis}[
        width=0.78\linewidth, height=0.78\linewidth,
        axis lines=box,
        tick label style={font=\figlabelsize},
        label style={font=\figlabelsize},
        legend style={font=\figlabelsize, draw=none},
        minor tick num=9,
        xlabel={$N$},
        ylabel={$T_{\mathrm{rec}}\ (s)$},
        ylabel style={
        at={(0.22,0.6)},
        anchor=south,
        font=\figlabelsize,
        },
        xlabel style={
        at={(0.6,-0.2)},
        anchor=south,
        font=\figlabelsize,
        }      ]
\addplot+[thick, mark=*, mark size=1.8pt]
          table[x=N, y={Trec_k10}, col sep=comma]{rev_laplace_sphere.csv};

        \addplot+[thick, mark=square*, mark size=1.8pt]
          table[x=N, y={Trec_k30}, col sep=comma]{rev_laplace_sphere.csv};
      \addONguideRel{0.5}{0.25}{0.3}
      \end{loglogaxis}
    \end{tikzpicture}
  \end{subfigure}\hfill
\begin{subfigure}{0.33\textwidth}
    \centering
    \begin{tikzpicture}
      \begin{loglogaxis}[
        width=0.78\linewidth, height=0.78\linewidth,
        axis lines=box,
        tick label style={font=\figlabelsize},
        label style={font=\figlabelsize},
        legend style={font=\figlabelsize, draw=none},
        minor tick num=9,
        xlabel={$N$},
        ylabel={$M$ (GB)},
        ylabel style={
        at={(0.22,0.6)},
        anchor=south,
        font=\figlabelsize,
        },
        xlabel style={
        at={(0.6,-0.2)},
        anchor=south,
        font=\figlabelsize,
        }
      ]
        \addplot+[thick, mark=*, mark size=1.8pt]
          table[x=N, y={M_GB_k10}, col sep=comma]{rev_laplace_sphere.csv};

        \addplot+[thick, mark=square*, mark size=1.8pt]
          table[x=N, y={M_GB_k30}, col sep=comma]{rev_laplace_sphere.csv};
      \addONguideRel{0.5}{0.25}{0.3}
      \end{loglogaxis}
    \end{tikzpicture}
  \end{subfigure}\hfill
\begin{subfigure}{0.33\textwidth}
    \centering
    \begin{tikzpicture}
      \begin{loglogaxis}[
        legend to name=sharedlegend_exp2,
        legend columns=2,
        width=0.78\linewidth, height=0.78\linewidth,
        axis lines=box,
        tick label style={font=\figlabelsize},
        label style={font=\figlabelsize},
        legend style={font=\figlabelsize, draw=gray,    column sep=1.2em,
  row sep=0.3em,
  inner xsep=0.8em,
  legend cell align={left}},
        minor tick num=9,
        xlabel={$N$},
        ylabel={$T_{\mathrm{sol}}\ (s)$},
        ylabel style={
        at={(0.22,0.6)},
        anchor=south,
        font=\figlabelsize,
        },
        xlabel style={
        at={(0.6,-0.2)},
        anchor=south,
        font=\figlabelsize,
        }
      ]
        \addplot+[thick, mark=*, mark size=1.8pt]
          table[x=N, y={Tsolve_k10}, col sep=comma]{rev_laplace_sphere.csv};
        \addlegendentry{\shortstack{$k=10,\ p=10$\\$n_{\mathrm{samples}}=770$}}

        \addplot+[thick, mark=square*, mark size=1.8pt]
          table[x=N, y={Tsolve_k30}, col sep=comma]{rev_laplace_sphere.csv};
        \addlegendentry{\shortstack{$k=30,\ p=30$\\$n_{\mathrm{samples}}=2424$}}
      \addONguideRel{0.5}{0.25}{0.3}
      \end{loglogaxis}
    \end{tikzpicture}
  \end{subfigure}

\vspace{0.5em}

\centering

\noindent
\hspace*{0.1\textwidth}\begin{minipage}[c]{0.20\textwidth}
\centering
\begingroup
\resizebox{\linewidth}{!}{\begin{tikzpicture}[scale=3]

\def\R{1.0}          \def\Bands{24}       \def\target{0.22}    \def\fdot{0.018}     \def\bdot{0.012}     \def\backop{0.35}    \def\tilt{70}        \def\rot{20}         

\shade[ball color=black!8] (0,0) circle (\R);

\draw[black!30,line width=0.45pt] (0,0) circle (\R);

\foreach \m in {1,...,\Bands} {
    \pgfmathsetmacro\z{-1 + 2*\m/(\Bands+1)}   \ifdim \z pt < 0pt
      \pgfmathsetmacro\rxy{sqrt(1 - \z*\z)}
      \pgfmathsetmacro\Nk{max(4, round(2*pi*\rxy/\target))}
      \pgfmathsetmacro\phiShift{180/\Nk} \foreach \k in {0,...,\Nk} {
        \pgfmathsetmacro\phi{360*\k/\Nk + \phiShift*(mod(\m,2))}
        \pgfmathsetmacro\x{\R*\rxy*cos(\phi)}
        \pgfmathsetmacro\y{\R*\rxy*sin(\phi)}
        \fill[black!60,opacity=\backop] (\x,\y) circle (\bdot);
      }
    \fi
  }
\fill[black!60,opacity=\backop] (0,0) circle (\bdot*0.85);

\foreach \m in {1,...,\Bands} {
    \pgfmathsetmacro\z{-1 + 2*\m/(\Bands+1)}
    \ifdim \z pt > 0pt
      \pgfmathsetmacro\rxy{sqrt(1 - \z*\z)}
      \pgfmathsetmacro\Nk{max(4, round(2*pi*\rxy/\target))}
      \pgfmathsetmacro\phiShift{180/\Nk}
      \foreach \k in {0,...,\Nk} {
        \pgfmathsetmacro\phi{360*\k/\Nk + \phiShift*(mod(\m,2))}
        \pgfmathsetmacro\x{\R*\rxy*cos(\phi)}
        \pgfmathsetmacro\y{\R*\rxy*sin(\phi)}
        \fill[black!70] (\x,\y) circle (\fdot);
      }
    \fi
  }
\fill[black!70] (0,0) circle (\fdot*0.85);

\begin{scope}[rotate=\rot]
    \draw[black!55,dashed,line width=0.4pt]
      (0,0) ellipse ({\R} and {\R*cos(\tilt)});
  \end{scope}
\begin{scope}[rotate=\numexpr\rot+90\relax]
    \draw[black!55,dashdotted,line width=0.4pt]
      (0,0) ellipse ({\R} and {\R*cos(\tilt)});
  \end{scope}

\end{tikzpicture}
 }
\endgroup
\end{minipage}\hfill
\begin{minipage}[c]{0.68\textwidth}
\centering
\pgfplotslegendfromname{sharedlegend_exp2}\\[0.4em]

\underline{\textbf{Experiment 2}} \hspace{1em} {\figtitlesize $n_{\mathrm{samples}} \approx 80\,k + p$}
\end{minipage}

  \caption{Reconstruction time, memory usage, and solve time vs.\ degrees of freedom for the discretization of the second-kind integral equation in \eqref{eq:bie}.}
  \label{fig:exp2}
\end{figure} 
\begin{table}[htb!]
\centering

\begin{subtable}[t]{\linewidth}
    \centering \tablesize
    \caption{Accuracy and solver effectiveness for $k=10,\ p=10$.}
    \noindent \clearpage{}\begin{tabular}{rc c ccc}
\toprule
\multicolumn{1}{c}{$N$} & $\mathrm{cond}(\mathbf A)$ & $n_{\mathrm{samples}}$ & $\mathrm{relerr}$ & $\mathrm{errsolve}$ & $\mathrm{relerr}_{\mathrm{bvp}}$  \\
\midrule
20,480 & 2.0e+00 & 682 & 2.4e-04 & 4.6e-04 & 8.1e-04 \\
81,920 & 2.0e+00 & 731 & 4.9e-04 & 9.7e-04 & 3.4e-04 \\
327,680 & 2.0e+00 & 764 & 4.9e-04 & 9.2e-04 & 1.5e-04 \\
1,310,720 & 2.0e+00 & 770 & 6.7e-04 & 1.3e-03 & 1.4e-04 \\
\bottomrule
\end{tabular}

\clearpage{}
    \label{tab:laplace_acc1}
\end{subtable}
\\
\begin{subtable}[t]{\linewidth}
    \centering \tablesize
    \caption{Accuracy and solver effectiveness for $k=30,\ p=30$.}
    \noindent \clearpage{}\begin{tabular}{rc c ccc}
\toprule
\multicolumn{1}{c}{$N$} & $\mathrm{cond}(\mathbf A)$ & $n_{\mathrm{samples}}$ & relerr & errsolve & $\mathrm{relerr}_{\mathrm{bvp}}$  \\
\midrule
20,480 & 2.0e+00 & 1,945 & 1.3e-06 & 2.4e-06 & 7.2e-04 \\
81,920 & 2.0e+00 & 2,089 & 2.1e-06 & 4.0e-06 & 3.3e-04 \\
327,680 & 2.0e+00 & 2,424 & 1.2e-06 & 2.2e-06 & 1.6e-04 \\
1,310,720 & 2.0e+00 & 2,310 & 3.4e-06 & 6.7e-06 & 7.7e-05 \\
\bottomrule
\end{tabular}

\clearpage{}
    \label{tab:laplace_acc2}
\end{subtable}
\caption{Summary of results for Experiment 2. 
Accuracy of \texttt{RSRS} and its effectiveness as a solver across problem sizes for (a) $k=10, p=10$ and (b) $k=30, p=30$.}
\label{tab:exp2}
\end{table}

The matrix in this example is well conditioned, with $\mathrm{cond}(\mtx A)\approx2$ for all~$N$, and both $\mathrm{relerr}$ and $\mathrm{errsolve}$ remain stable across problem sizes (Table~\ref{tab:exp2}).
For $k=10$, relerr is on the order of $10^{-4}$, while increasing to $k=30$ improves accuracy by roughly two orders of magnitude, with very little degradation as~$N$ increases.
The boundary-value error $\mathrm{relerr}_{\mathrm{bvp}}$ decreases steadily as the mesh is refined, and using \texttt{RSRS} as a direct solver yields errors comparable to the underlying discretization error---indicating that the boundary-value problem is effectively solved to quadrature accuracy.

With the choice of $m=6k$ for the leaf parameter, the number of samples required is $n_{\mathrm{samples}}\approx80\,k+p$, reflecting the smaller geometric constant associated with a surface (rather than volumetric) distribution (cf.~Section~\ref{ssub:nonuniform_point_distributions}).
Reconstruction and solve times scale linearly with~$N$ (Figure~\ref{fig:exp2}), consistent with the $\mathcal{O}(k^2N)$ and $\mathcal{O}(kN)$ complexity estimates.

\subsection{Experiment 3: Schur Complement Arising in a Sparse Direct Solver for PDEs}
\label{sec:schur_compression}
\begin{figure}[htb!]
\centering
\resizebox{0.8\textwidth}{!}{\begin{tikzpicture}[scale=0.22]

\pgfmathsetmacro{\sep}{5}
\pgfmathsetmacro{\n}{2*5+1}

\foreach \x in {0,...,\n}
{
  \draw (\x,0,\n) -- (\x,\n,\n);
  \draw (\x,\n,\n) -- (\x,\n,0);
}
\foreach \x in {0,...,\n}
{
  \draw (\n,\x,\n) -- (\n,\x,0);
  \draw (0,\x,\n) -- (\n,\x,\n);
}
\foreach \x in {0,...,\n}
{
  \draw (\n,0,\x) -- (\n,\n,\x);
  \draw (0,\n,\x) -- (\n,\n,\x);
}

\foreach \x in {0,...,\n}
  \foreach \y in {0,...,\n}
    \filldraw[gray] (\x,\y,\n) circle (4pt);
\foreach \x in {0,...,\n}
  \foreach \z in {0,...,\n}
    \filldraw[gray] (\x,\n,\z) circle (4pt);
\foreach \y in {0,...,\n}
  \foreach \z in {0,...,\n}
    \filldraw[gray] (\n,\y,\z) circle (4pt);

\foreach \x in {0,...,\n}
{
  \filldraw[blue] (\x,\sep,\n) circle (4pt);
  \filldraw[blue] (\x,\n,\sep) circle (4pt);
  \filldraw[blue] (\n,\x,\sep) circle (4pt);
  \filldraw[blue] (\sep,\x,\n) circle (4pt);
  \filldraw[blue] (\sep,\n,\x) circle (4pt);
  \filldraw[blue] (\n,\sep,\x) circle (4pt);
}

\node [below,blue] at (\sep,0,\n) {$I_1$};
\node [below right,gray] at (\n,0,\n) {$I_2, \dots I_9$};

\node[align=center,font=\figlabelsize] at (\n*0.25,-0.75*\n) {Sparse Matrix $\mtx{A}$};

\draw[->, thick] (\n+\sep,\sep) -- (3*\n-\sep,\sep);

\node[above,font=\figlabelsize] at (2*\n,\sep) {Eliminate\ $\textcolor{gray}{I_2, \dots, I_9}$};

\node at (3*\n,0)
      [anchor=south west, inner sep=0, outer sep=0, yshift=-1.2cm]
      {\includegraphics[width=0.25\textwidth]{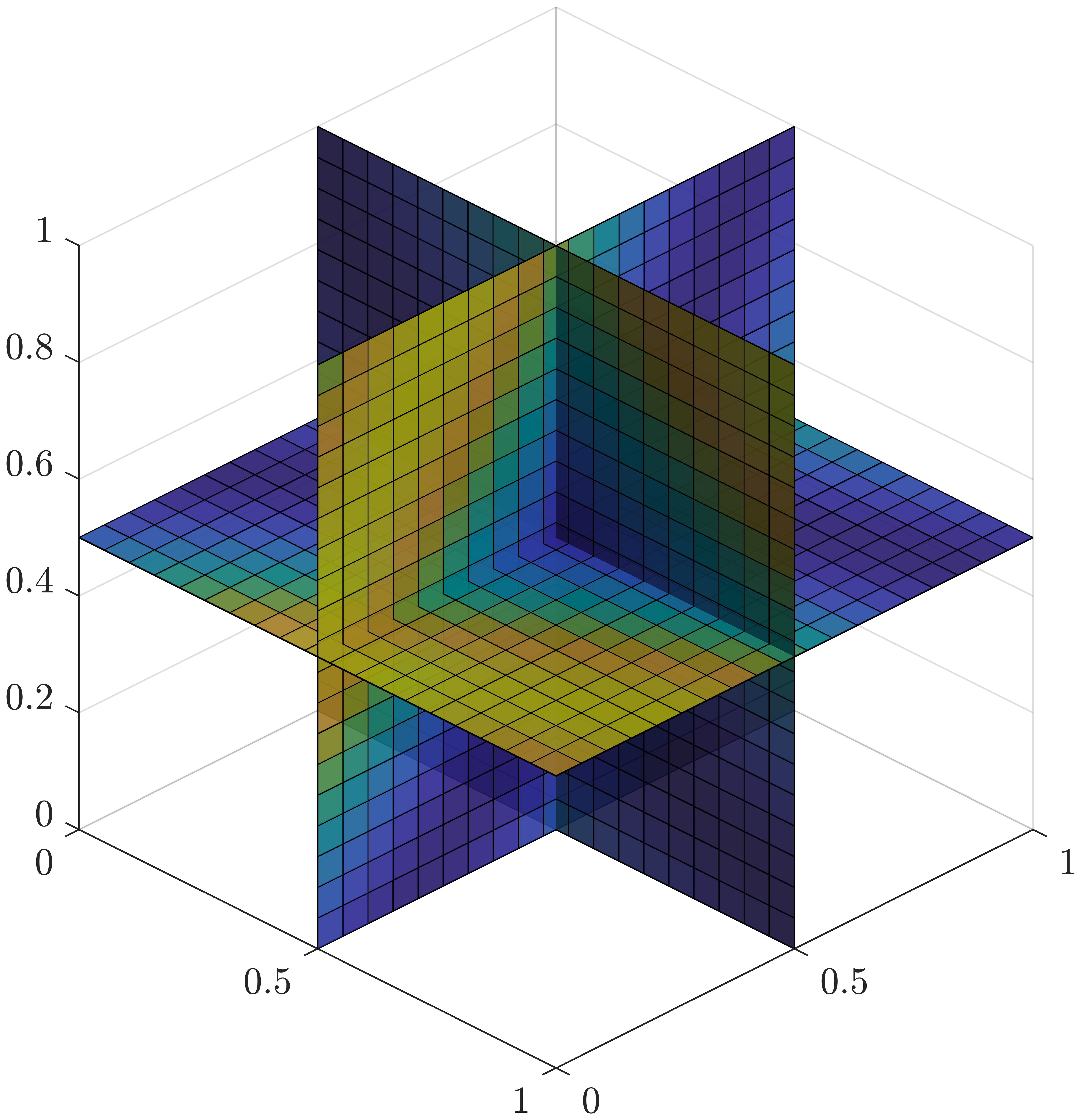}};

\node[below,anchor=west,font=\figlabelsize] at (2.75*\n,-0.75*\n) {Dense Schur Complement $\mtx{S}_{11}$};
\end{tikzpicture} }
\caption{Discretization of the 3D Helmholtz equation using a 7-point stencil produces a sparse matrix $\mtx{A}$ (left). Eliminating the gray interior octants results in a dense Schur complement on the separator $I_1$ (right).}
\label{fig:exp3_geom}
\end{figure}

In this example, we demonstrate how $\mathcal{H}^2$ matrix factorizations can be leveraged to improve the complexity and memory footprint of sparse solvers, building upon prior work~\cite{2013_xia_randomized,2009_xia_superfast,yesypenko2024slablu}.
We consider the 3D Helmholtz equation,
\begin{equation}
(-\Delta - \kappa^2) u(\pxx) = f(\pxx), \quad \pxx \in \Omega = [0,1]^3, 
\qquad 
u(\pxx) = g(\pxx), \quad \pxx \in \partial\Omega.
\label{eq:bvp_helm}
\end{equation}
discretized using a standard 7-point finite difference stencil on a uniform grid. The resulting sparse matrix $\mtx{A}$ has the familiar 7-point stencil sparsity pattern, illustrated in Figure~\ref{fig:exp3_geom}. We set the wavenumber as $\kappa=13.6$, so that there are about two wavelengths across the unit domain.

To factorize $\mtx A$, we partition the domain into eight interior gray octants and a top-level separator $I_1$, as shown in Figure~\ref{fig:exp3_geom}. The interior octants are eliminated recursively first, and the separator $I_1$ is eliminated last. This ordering minimizes fill-in and isolates dense Schur complements to smaller separators.
The sparsity pattern of $\mtx{A}$ after this reordering is:
\begin{equation}
\mtx A = \begin{bmatrix}
\mtx A_{99} & &&&\mtx A_{91}\\
& \mtx A_{88} & && \mtx A_{81}\\
&& \ddots &&\vdots\\
&&& \mtx A_{22}& \mtx A_{21}\\
\mtx A_{19}& \mtx A_{18} & \dots & \mtx A_{12} & \mtx A_{11}
\end{bmatrix},
\label{eq:Astiffness}
\end{equation}
with an associated invertible factorization
$\mtx A = \mtx L\ \mtx D\ \mtx U,$
where $\mtx L$ and $\mtx U$ are sparse triangular factors, and $\mtx D$ is block diagonal with the Schur complement $\mtx S_{11}$ on $I_1$, defined by:
\begin{equation}
\mtx S_{11} = \mtx A_{11} - \sum_{j=2}^9 \mtx A_{1j} \mtx A_{jj}^{-1} \mtx A_{j1}.
\label{eq:S11}
\end{equation}
Each factor $\mtx A_{jj}^{-1} \mtx A_{j1}$ is dense, leading to a dense Schur complement $\mtx S_{11}$ in the sparse factorization of \eqref{eq:Astiffness}. The matrices $\mtx A_{jj}$ correspond to sparse subdomains that can be efficiently factorized using standard sparse direct solvers. Applying $\mtx A_{jj}^{-1}$ to a vector requires $\mathcal{O}(N^{4/3})$ operations, making the application of $\mtx S_{11}$ and its adjoint efficient when amortized over multiple right-hand sides.
Figure~\ref{fig:exp3} summarizes the numerical results for this experiment, while Table~\ref{tab:exp3} reports the iteration counts for solving the Schur complement system using GMRES, both with and without the \texttt{RSRS}-based solver.

By replacing the dense Schur complement $\mtx{S}_{11}$ with an approximate invertible factorization, 
we can efficiently solve the discretized Helmholtz problem~\eqref{eq:bvp_helm}.
The Dirichlet boundary data are prescribed from an analytic solution corresponding to the free-space Helmholtz Green’s function,
\[
u(\pxx) = \frac{\sin(\kappa \|\pxx - \pyy\|)}{\|\pxx - \pyy\|},
\]
where the source point~$\pyy$ is located outside the domain~$\Omega$.
The relative error between the analytic solution~$u$ and the computed solution~$\hat{u}$ 
inside~$\Omega$ is reported in Table~\ref{tab:exp3} as $\mathrm{relerr}_{\mathrm{bvp}}$, defined in~\eqref{eq:relerr_bvp}.

\begin{figure}[htb!]
  \centering

\begin{subfigure}{0.33\textwidth}
    \centering
    \begin{tikzpicture}
      \begin{loglogaxis}[
        width=0.78\linewidth, height=0.78\linewidth,
        axis lines=box,
        tick label style={font=\figlabelsize},
        label style={font=\figlabelsize},
        legend style={font=\figlabelsize, draw=none},
        minor tick num=9,
        xlabel={$N$},
        ylabel={$T_{\mathrm{rec}}\ (s)$},
        ylabel style={
        at={(0.22,0.6)},
        anchor=south,
        font=\figlabelsize,
        },
        xlabel style={
        at={(0.65,-0.2)},
        anchor=south,
        font=\figlabelsize,
        }      ]
\addplot+[thick, mark=*, mark size=1.8pt]
          table[x=N, y={Trec_k60}, col sep=comma]{rev_helm_schur.csv};

        \addplot+[thick, mark=square*, mark size=1.8pt]
          table[x=N, y={Trec_k80}, col sep=comma]{rev_helm_schur.csv};
      \addONguideRel{0.5}{0.2}{0.3}
      \end{loglogaxis}
    \end{tikzpicture}
  \end{subfigure}\hfill
\begin{subfigure}{0.33\textwidth}
    \centering
    \begin{tikzpicture}
      \begin{loglogaxis}[
        width=0.78\linewidth, height=0.78\linewidth,
        axis lines=box,
        tick label style={font=\figlabelsize},
        label style={font=\figlabelsize},
        legend style={font=\figlabelsize, draw=none},
        minor tick num=9,
        xlabel={$N$},
        ylabel={$M$ (GB)},
        ylabel style={
        at={(0.22,0.6)},
        anchor=south,
        font=\figlabelsize,
        },
        xlabel style={
        at={(0.65,-0.2)},
        anchor=south,
        font=\figlabelsize,
        }
      ]
        \addplot+[thick, mark=*, mark size=1.8pt]
          table[x=N, y={M_GB_k60}, col sep=comma]{rev_helm_schur.csv};

        \addplot+[thick, mark=square*, mark size=1.8pt]
          table[x=N, y={M_GB_k80}, col sep=comma]{rev_helm_schur.csv};
      \addONguideRel{0.5}{0.2}{0.3}
      \end{loglogaxis}
    \end{tikzpicture}
  \end{subfigure}\hfill
\begin{subfigure}{0.33\textwidth}
    \centering
    \begin{tikzpicture}
      \begin{loglogaxis}[
        legend to name=sharedlegend_exp3,
        legend columns=2,
        width=0.78\linewidth, height=0.78\linewidth,
        axis lines=box,
        tick label style={font=\figlabelsize},
        label style={font=\figlabelsize},
        legend style={font=\figlabelsize, draw=gray,    column sep=1.2em,
  row sep=0.3em,
  inner xsep=0.8em,
  legend cell align={left}},
        minor tick num=9,
        xlabel={$N$},
        ylabel={$T_{\mathrm{sol}}\ (s)$},
        ylabel style={
        at={(0.22,0.6)},
        anchor=south,
        font=\figlabelsize,
        },
        xlabel style={
        at={(0.65,-0.2)},
        anchor=south,
        font=\figlabelsize,
        }
      ]
        \addplot+[thick, mark=*, mark size=1.8pt]
          table[x=N, y={Tsolve_k60}, col sep=comma]{rev_helm_schur.csv};
        \addlegendentry{\shortstack{$k=60,\ p=60$\\$n_{\mathrm{samples}}=5580$}}

        \addplot+[thick, mark=square*, mark size=1.8pt]
          table[x=N, y={Tsolve_k80}, col sep=comma]{rev_helm_schur.csv};
        \addlegendentry{\shortstack{$k=80,\ p=80$\\$n_{\mathrm{samples}}=8643$}}
      \addONguideRel{0.5}{0.2}{0.3}
      \end{loglogaxis}
    \end{tikzpicture}
  \end{subfigure}

\vspace{0.5em}

\centering

\noindent
\hspace*{0.05\textwidth}\begin{minipage}[c]{0.25\textwidth}
\centering
\begingroup
\resizebox{\linewidth}{!}{\tdplotsetmaincoords{70}{120}

\begin{tikzpicture}[tdplot_main_coords, scale=4/3, line cap=round]
\def\L{1.0}      \def\N{14}       \def\ax{0.50}    \def\ay{0.50}    \def\az{0.50}    

\colorlet{planeX}{black!14}   \colorlet{planeY}{black!12}
  \colorlet{planeZ}{black!18}
  \def\op{0.65}                 \tikzset{mesh/.style={black!55, very thin}}

\newcommand{\DrawXpatch}[4]{\begin{scope}
    \clip (\ax,#1,#3) -- (\ax,#2,#3) -- (\ax,#2,#4) -- (\ax,#1,#4) -- cycle;
    \fill[planeX,opacity=\op] (\ax,0,0) -- (\ax,\L,0) -- (\ax,\L,\L) -- (\ax,0,\L) -- cycle;
    \foreach \k in {0,...,\N} {\pgfmathsetmacro\t{\k/\N*\L}
      \draw[mesh] (\ax,\t,0) -- (\ax,\t,\L);
      \draw[mesh] (\ax,0,\t) -- (\ax,\L,\t);
    }
  \end{scope}}

\newcommand{\DrawYpatch}[4]{\begin{scope}
    \clip (#1,\ay,#3) -- (#2,\ay,#3) -- (#2,\ay,#4) -- (#1,\ay,#4) -- cycle;
    \fill[planeY,opacity=\op] (0,\ay,0) -- (\L,\ay,0) -- (\L,\ay,\L) -- (0,\ay,\L) -- cycle;
    \foreach \k in {0,...,\N} {\pgfmathsetmacro\t{\k/\N*\L}
      \draw[mesh] (0,\ay,\t) -- (\L,\ay,\t);
      \draw[mesh] (\t,\ay,0) -- (\t,\ay,\L);
    }
  \end{scope}}

\newcommand{\DrawZpatch}[4]{\begin{scope}
    \clip (#1,#3,\az) -- (#2,#3,\az) -- (#2,#4,\az) -- (#1,#4,\az) -- cycle;
    \fill[planeZ,opacity=\op] (0,0,\az) -- (\L,0,\az) -- (\L,\L,\az) -- (0,\L,\az) -- cycle;
    \foreach \k in {0,...,\N} {\pgfmathsetmacro\t{\k/\N*\L}
      \draw[mesh] (\t,0,\az) -- (\t,\L,\az);
      \draw[mesh] (0,\t,\az) -- (\L,\t,\az);
    }
  \end{scope}}

\DrawXpatch{0}{\ay}{0}{\az}
\DrawYpatch{\ax}{\L}{0}{\az}
\DrawYpatch{0}{\ax}{0}{\az}
\DrawXpatch{\ay}{\L}{0}{\az}

\DrawZpatch{\ax}{\L}{0}{\ay}
\DrawZpatch{0}{\ax}{0}{\ay}
\DrawZpatch{\ax}{\L}{\ay}{\L}
\DrawZpatch{0}{\ax}{\ay}{\L}

\DrawXpatch{0}{\ay}{\az}{\L}
\DrawYpatch{\ax}{\L}{\az}{\L}
\DrawYpatch{0}{\ax}{\az}{\L}
\DrawXpatch{\ay}{\L}{\az}{\L}

\draw[black!40, line width=0.45pt]
    (0,0,0) -- (\L,0,0) -- (\L,\L,0) -- (0,\L,0) -- cycle
    (0,0,\L) -- (\L,0,\L) -- (\L,\L,\L) -- (0,\L,\L) -- cycle
    (0,0,0) -- (0,0,\L)
    (\L,0,0) -- (\L,0,\L)
    (\L,\L,0) -- (\L,\L,\L)
    (0,\L,0) -- (0,\L,\L);

\end{tikzpicture}
 }
\endgroup
\end{minipage}\hfill
\begin{minipage}[c]{0.68\textwidth}
\centering
\pgfplotslegendfromname{sharedlegend_exp3}\\[0.4em]

\underline{\textbf{Experiment 3}} \hspace{1em} {\figtitlesize $n_{\mathrm{samples}} \approx 107\,k + p$}
\end{minipage}

  \caption{Reconstruction time, memory usage, and solve time vs.\ degrees of freedom for the Schur complement $\mtx S_{11}$ in \eqref{eq:S11} arising during the course of the sparse factorization 
  of \eqref{eq:Astiffness}.}
  \label{fig:exp3}
\end{figure} 
\begin{table}[htb!]
\centering

\begin{subtable}[t]{\linewidth}
    \centering \tablesize
    \caption{Accuracy and solver effectiveness.}
    \noindent \clearpage{}\begin{tabular}{rc ccc ccc}
\toprule
\multicolumn{2}{c}{} & \multicolumn{3}{c}{$k=60,\ p = 60$} & \multicolumn{3}{c}{$k=80,\ p = 80$} \\ \cmidrule(lr){3-5} \cmidrule(lr){6-8} 
\multicolumn{1}{c}{$N$} & $\rm cond(\mtx A)$ & $\mathrm{relerr}$ & $\mathrm{errsolve}$ & $\mathrm{relerr}_{\mathrm {bvp}}$& $\mathrm{relerr}$ & $\mathrm{errsolve}$ & $\mathrm{relerr}_{\mathrm {bvp}}$ \\
\midrule
40,717 & 2.7e+03 & 6.3e-05 & 2.4e-03 & 4.1e-03 & 9.0e-06 & 2.8e-04 & 3.9e-03 \\
81,181 & 3.8e+03 & 1.9e-04 & 6.9e-03 & 7.8e-03 & 3.8e-05 & 1.2e-03 & 2.1e-03 \\
159,391 & 5.2e+03 & 2.0e-03 & 5.8e-02 & 2.2e-02 & 9.7e-05 & 4.3e-03 & 3.1e-03 \\
319,807 & 7.4e+03 & 8.0e-04 & 1.3e-01 & 4.1e-02 & 3.9e-05 & 3.0e-03 & 8.3e-04 \\
\bottomrule
\end{tabular}

\clearpage{}
    \label{tab:helm_acc}
\end{subtable}
\\
\begin{subtable}[t]{\linewidth}
    \centering \tablesize
    \caption{Preconditioner performance.}
    \noindent \clearpage{}\begin{tabular}{r c cc cc}
\toprule
 & {No precond} & \multicolumn{2}{c}{$k=60,\ p = 60$} & \multicolumn{2}{c}{$k=80,\ p = 80$} \\
\cmidrule(lr){2-2} \cmidrule(lr){3-4} \cmidrule(lr){5-6}
\multicolumn{1}{c}{$N$} & $n_{\rm iter}$ & $n_{\rm samples}$ & $n_{\rm iter}$ & $n_{\rm samples}$ & $n_{\rm iter}$ \\
\midrule
40,717 & $>$10,000 & 5,580 & 4 & 7,440 & 3 \\
81,181 & $>$10,000 & 5,580 & 5 & 8,643 & 4 \\
159,391 & $>$10,000 & 5,580 & 7 & 7,440 & 5 \\
319,807 & $>$10,000 & 5,580 & 9 & 8,561 & 4 \\
\bottomrule
\end{tabular}

\clearpage{}
    \label{tab:helm_niter}
\end{subtable}
\caption{Summary of results for Experiment 3. 
(a): Accuracy of the \texttt{RSRS} factorization and resulting direct solver across mesh refinements for two rank choices.
(b): GMRES iteration counts with and without \texttt{RSRS} preconditioning, demonstrating robust performance.}
\label{tab:exp3}
\end{table}

The 3D Helmholtz Schur complement is an indefinite matrix, with the conditioning growing moderately with $N$.
Across all problem sizes, the \texttt{RSRS} factorization produces accurate compressions, 
and increasing the rank parameter from $k=60$ to $k=80$ yields roughly one to two orders of magnitude improvement in
$\mathrm{relerr}$ and $\mathrm{errsolve}$. 
The boundary-value error $\mathrm{relerr}_{\mathrm{bvp}}$ decreases as the mesh is refined and
drops below $10^{-3}$ at the finest resolution for $k=80$.

As shown in Table~\ref{tab:exp3}, the unpreconditioned Schur system fails to converge 
within $10^4$ GMRES iterations, whereas the \texttt{RSRS}-based preconditioner reduces 
the iteration count to only $4$--$9$ iterations for $k=60$ and $3$--$5$ iterations for 
$k=80$. With $m=6k$ as the leaf parameter, the observed sample costs are
$n_{\mathrm{samples}} \approx 107\,k + p$. Together with the linear trends in
$T_{\mathrm{rec}}$, $T_{\mathrm{sol}}$, and memory shown in Figure~\ref{fig:exp3}, these results
are consistent with the complexity estimates in Section~\ref{ssub:nonuniform_point_distributions}.

\section{Conclusions}

This manuscript introduces \texttt{RSRS}, an algorithm for simultaneously compressing and inverting $\mathcal{H}^2$ matrices in matrix-free settings. By leveraging novel randomized sketching techniques, \texttt{RSRS} constructs an approximate invertible factorization using random sketches of the matrix and its adjoint. Dense Gaussian test matrices are used for sketching, and the necessary structure is introduced via linear algebraic post-processing techniques, eliminating the need for carefully structured test matrices and avoiding excessive sample costs.

The algorithm builds on and extends recently proposed LU factorization methods for $\mathcal{H}^2$ matrices~\cite{2017_ho_ying_strong_RS,sushnikova2023fmm}, adapting them to settings where matrix entries are inaccessible. Its effectiveness is demonstrated across a range of applications, including both integral and differential equations. For ill-conditioned problems where iterative solvers fail to converge, \texttt{RSRS} computes a highly effective solver using only a modest number of matrix and adjoint-vector products. The algorithm is especially advantageous for multiple solves, where the sampling costs of \texttt{RSRS} are naturally amortized, though our experiments demonstrate that these amortization benefits can already be observed for a single solve.

\texttt{RSRS} opens several promising avenues for future research. These include the analysis of sketch accuracy under recursive transformations, high-performance implementations in distributed environments, and hybrid strategies that combine weak and strong admissibility to balance computational cost and accuracy. Beyond PDEs and integral equations, the method has potential applications in inverse problems, uncertainty quantification, and PDE-constrained optimization---particularly in settings where the matrix is only accessible via its action, and robust fast solvers are essential for tackling large-scale problems.

\subsection*{Acknowledgments}

The work reported was supported by the Office of Naval Research (N00014-18-1-2354), by the National Science Foundation (DMS-1952735 and DMS-2313434), and by the Department of Energy ASCR (DE-SC0022251).

\subsection*{Data Availability} Not applicable.

\section*{Declarations} 
\subsection*{Conflict of Interest} The authors have no relevant financial or non-financial interests to disclose. 

\bibliographystyle{plain}
\bibliography{main_bib,main_bib_extra}

\end{document}